\newif\ifdraft
\drafttrue
 
\documentclass[3p, times, review ]{elsarticle}

\usepackage[utf8]{inputenc}
\usepackage[english]{babel}
  \usepackage[]{algorithm2e}
\usepackage[utf8]{inputenc}
\usepackage[T1]{fontenc}
\usepackage{array}
\usepackage{amssymb}
\usepackage{amsfonts}
\usepackage{amsmath}
\usepackage{cases}
\usepackage{hyperref}
\newcolumntype{\xx}{\bm{x}}

\newcommand{\var}{\texttt}
\usepackage{amsthm}
\usepackage{graphicx}
\usepackage{subfig} 
\usepackage{textcomp}
\usepackage{bm}
\usepackage{booktabs}
\usepackage{mathtools}

\makeatletter

\newcommand{\Rmnum}[1]{\expandafter\@slowromancap\romannumeral #1@}
\makeatother

\theoremstyle{definition}

\newtheorem{Remark}{Remark}[section]

\newtheorem{form}{Formulation}[section]

\newcommand{\R}{\mathbb{R}}

\newcommand{\blue}[1]{\textcolor{black}{#1}}

\newcommand{\unit}[2]{#1\,\mathrm{#2}}

\usepackage{color}

\ifdraft
\usepackage{color}
\usepackage[draft,disable]{todonotes}

\else
\newcommand{\todo}[1]{}

\fi
\newcommand{\noii}[1]{\textcolor{black}{#1}}

\newcommand{\new}[1]{\textcolor{black}{#1}}

\newcommand{\xx}{\textbf{x}}

\newcommand{\RN}[1]{%
  \textup{\uppercase\expandafter{\romannumeral#1}}%
}

\usepackage[figuresright]{rotating}

\begin{document}

\begin{frontmatter}


\title{A Bayesian estimation method for variational phase-field fracture problems}
\author[tuwien,Hannover]{Amirreza Khodadadian\corref{cor1}}
\ead{khodadadian@ifam.uni-hannover.de}
\author[Hannover]{Nima Noii}
\ead{noii@ifam.uni-hannover.de}
\author[tuwien]{Maryam Parvizi}
\ead{maryam.parvizi@tuwien.ac.at}
\author[Amir]{Mostafa Abbaszadeh}
\ead{m.abbaszadeh@aut.ac.ir}
\author[Hannover]{Thomas Wick}
\ead{thomas.wick@ifam.uni-hannover.de}
\author[tuwien,ASU]{Clemens Heitzinger}
\ead{clemens.heitzinger@tuwien.ac.at}

\cortext[cor1]{Corresponding author}

\address[tuwien]{Institute of Analysis and Scientific Computing, 
  Vienna University of Technology (TU Wien),
  Wiedner Hauptstraße 8--10, 
  1040 Vienna, Austria}
\address[Amir]{Faculty of Mathematics and Computer Sciences,
  Amirkabir University of Technology, No.~424, Hafez Ave., 15914, Tehran, Iran}
\address[Hannover]{Institute of Applied Mathematics, Leibniz University Hannover, Welfengarten 1, 30167
  Hanover, Germany}
\address[ASU]{School of Mathematical and Statistical Sciences, Arizona State University, Tempe, AZ 85287, USA}

\date{}

\begin{abstract}
In this work, we propose a parameter estimation framework 
for fracture propagation problems. The fracture problem is 
described by a phase-field method. Parameter estimation 
is realized with a Bayesian approach. Here, the focus 
is on uncertainties arising in the solid material parameters
and the critical energy release rate. 
A reference value (obtained on a sufficiently {refined} mesh) 
as the replacement of measurement {data}
will be chosen, and their posterior distribution is obtained. 
Due to time- and mesh dependenc{ies} of the
problem, the computational costs can be high. 
Using Bayesian inversion, we solve the problem 
on a relatively coarse mesh and fit the parameters. 
{In several numerical examples 
our proposed framework is substantiated and
the obtained 
load-displacement curves, that are usually the target functions,
are matched with the reference values.}
\end{abstract}


\end{frontmatter}

\section{Introduction}\label{intro}
This work is devoted to parameter identifications in 
fracture failure problems. To formulate fracture phenomena, {a} 
phase-field formulation for quasi-brittle fracture is used. 
The variational phase-field formulation is 
a thermodynamically consistent framework to 
compute the fracture failure process. 
This formulation emanates from the regularized 
version of the sharp crack surface function, which was first 
modeled in a variational framework in \cite{FraMar98}.
Regularized fracture phenomena are described with an 
additional auxiliary smooth indicator function \cite{BourFraMar00}, 
which is denoted as 
crack phase-field (here indicated by $d$).  Along with a mechanical field  (denoted by $\bm{u}$), 
a minimization problem for the multi-field problem $(\bm{u},d)$ can be formulated. The main feature of such a variational formulation 
is to approximate the discontinuities in $\bm{u}$ across the lower-dimensional 
crack topology with the phase-field function $d$.

The resulting, regularized formulation leads to a diffusive transition zone between two phases in the solid, which corresponds to the fractured phase (i.e., $d=0$) and intact phase (i.e., $d=1$), respectively. \new{The transition zone is determined by the phase-field regularization parameter $\ell$, 
also well-known as the length-scale parameter.  
The parameter $\ell$ is related to the element size $h$ and specifically $h\leq \ell$ (e.g., $\ell=2h$).
Therefore, sufficiently small length-scales are computationally demanding.}
To date, the focus in such cases was on local mesh adaptivity and parallel
computing in order to reduce the computational cost significantly; see
for instance \cite{noii2019phase,MaWi19,HeWheWi15,HeiWi18_pamm,BuOrSue13,BuOrSue10,ArFoMiPe15,Wi16_dwr_pff,MaWaWiWo19}. Another recent approach is a global-local technique in which
parts of the domain are solved with a simplified approach \cite{GeNoiiAllLo18,NoAlWiWr19} that
also aims to reduce the computational cost.

\new{Generally, material parameters fluctuate randomly in space. In fact, the mechanical material parameters are spatially variable and, therefore, the uncertainty related to spatially varying properties can be represented by random fields. For instance, the material stiffness property has spatial variability. 
In fact, there are several sources of uncertainty including the class of extensometer or strain gauge resolution,
	uncertainty in the dimensional measurements, the classification and resolution of the load cell, misalignment of the specimen or strain measurement device, temperature effects, operator-dependent factors, data fitting routines and analysis methods, etc \cite{lord2006measurement}. Therefore, in order to provide a reliable model, the uncertainty effect must be taken into account.}


The main goal in this work is to identify such \new{uncertain parameters} for phase-field 
fracture problems. 
The underlying framework of parameter estimation using Bayesian inference is described 
in the following.
Bayesian inference  is a probabilistic method used to estimate the unknown parameters according to the prior knowledge. The observations (experimental or synthetic measurements) can be used to update the prior data and provide the posterior estimation. The distribution provides useful information about the possible range of parameters and their variations and mean. Markov chain Monte Carlo (MCMC) \cite{hoang2013complexity} is a common computational approach for extracting information of the inverse problem (posterior distribution). Metropolis-Hastings (MH) algorithm \cite{hastings1970monte} is the most popular MCMC method to generate a Markov chain employing a proposal distribution for new steps. In practice, a reliable estimation of influential parameters is not possible or needs significant efforts.
In \cite{khodadadian2019bayesian,mirsian2019new}, the authors used the Metropolis-Hastings algorithm to estimate the unknown parameters in field-effect sensors. It enables authors to estimate probe-target density of the target molecules which can not be experimentally estimated. {We refer interested readers} to \cite{minson2013bayesian,cardiff2009bayesian} for more applications of Bayesian 
estimation. In the same line, other optimization approach{es} can be used to determine intrinsic material properties of the specimen
from experimental load-displacement curve{s}, 
see e.g., \cite{noii2019characterization}.

\new{As previously mentioned, we consider fractures in elastic solids in this work. 
	The principal material parameters are the shear modulus $\mu$  and the effective bulk modulus, $K=\lambda+ \frac{2\mu}{3}$ (here $\lambda$ denotes Lamé's first parameter) and Griffith's critical energy release rate $G_c$. Using Bayesian inversion, the objective is to determine the unknown elasticity parameters.
}

\new{For a homogeneous material, the stability requires positive-definiteness of the elasticity tensor. For an externally unconstrained homogeneous solid, the conditions of structural stability needs that the fourth-order stiffness tensor is positive-definite. The condition for an isotropic, linear elastic medium gives rise to the shear modulus $\mu$  and the effective bulk modulus, $K$ be strictly positive  \cite{kochmann2012analytical}.
  Regarding $\lambda$, the bound $\lambda>-\frac{2\mu}{3}$ may relate it to the shear modulus. Also, for the isotropic materials (as used in this paper) Poisson's ratio $\nu$ satisfies the condition $-1 < \nu < \frac{1}{2}$ \cite{greaves2011poisson}. These two elasticity parameters ($\lambda$ and $\nu$) 
are not well-suited for the estimation due to their bounds and dependency. Therefore, for the elasticity parameterization, we chose the eigenvalues, i.e., $K$ and $\mu$, and strive to estimate the joint probability, being updated jointly using MCMC.}

	\new{Griffith's theory describes that crack propagation occurs if a certain reduction of the potential energy due to the change of surface energy associated with incremental crack extension reaches to its critical value \cite{Zehnder2013}. Here, Griffith's critical energy rate  $G_c$ measures the amount of energy dissipated in a localized fracture state \cite{comi2001fracture}, thus has units of energy-per-unit-area. 
\blue{In case $G_c$ is unknown, one possibility is to employ the Bayesian setting for its identification.}
		Physically speaking, there is a direct relation between $G_c$ and material stiffness, which means that in stiffer materials more energy is needed for the crack initiation. Computationally speaking, this value is independent of the elasticity parameters.  
	Finally, since we should deal with three positive values ($\mu$, $K$) and $G_c$, in order to remove the positivity constraints, we transfer these parameters and estimate the transfered values $\mu^*=\log(\mu), K^*=\log(K)$, and $G_c^*=\log(G_c)$. } 
 
\new{In our Bayesian framework}, a reference value (obtained on a sufficiently refined mesh which termed here to the $\textit{virtual observation}$) as the 
replacement of measurement will be chosen. \new{Then, the posterior density of the elasticity parameters (joint probability) and the critical energy release rate} is obtained. 
The computational costs can be high, specifically when an appropriate 
estimation is required {inside multi-physics frameworks}, 
see e.g. \cite{miehe2016phase, heider2017phase, LeeWheWi16, noii2019phase}. Using Bayesian inversion; we strive to solve {such problems} 
with a coarser mesh and fit the parameters. The obtained load-displacement 
curve (as an important characteristic output) is matched with the reference value.  
In spite of using coarser meshes and therefore significantly lower computational 
costs (in terms of CPU timings), the accuracy of the solution is reliable 
(crack initiation and material 
fracture time estimated precisely).

The paper is organized as follows:
In Section \ref{section2}, we describe the variational isotropic phase-field formulation for the brittle fracture that is  a  thermodynamically  consistent  framework  to  compute
the  fracture  failure  process. In section \ref{section3}, the Bayesian inference is explained. We describe how the MH algorithm will be used to estimate the unknown parameters in phase-field fracture. Also we point out the critical points in the load-displacement curve, which must be estimated precisely with the Bayesian approach. 
In Section \ref{sec_Bay_for_PFF}, the Bayesian framework is adopted to estimate 
unknown parameters in the phase-field fracture approach.
In Section \ref{section4}, three specific numerical examples with different parameters and geometry will be given. We will use two proposal distributions (uniform and normal distribution) to sample the candidates and estimate the unknown parameters with different mesh sizes. Finally, in \ref{conclusions} we will draw paper conclusions and explain our future planes for employing Bayesian inversion in heterogeneous materials.

\section{Variational isotropic phase-field brittle fracture}
\label{section2}

\subsection{The primary fields for the variational phase-field formulation}
\label{sec_noniso}
\begin{figure}[!ht]
	\centering
	\includegraphics[width=14cm,height=6cm]{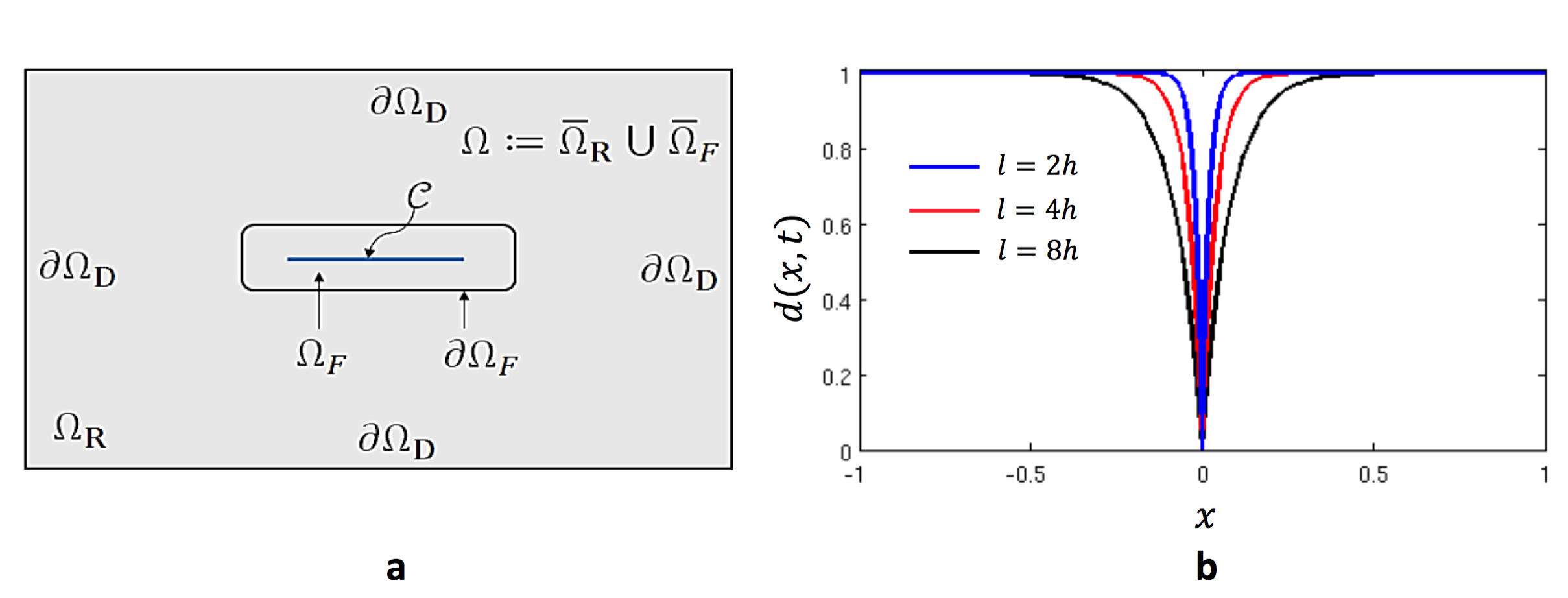}
	\caption{(a) Geometric setup: the intact region indicated by
		$\Omega_R$ and $\mathcal{C}$ is the crack phase-field
		surface. The entire domain is denoted by $\Omega$. The
		crack
		phase-field is approximated in the domain $\Omega_F$. The fracture boundary is $\partial
		\Omega_F$ and the outer boundary of the domain is $\partial
		\Omega$. $\Omega_F$ is represented by means of
		$d$ such that the transition area is $0< d< 1$ with
		thickness $2\ell$. (b) Regularized crack phase-field
		profile for a different length scale. A smaller value
		for the length scale lets the crack phase-field
		profile converge to a delta distribution.}
	\label{Figure1}
\end{figure}

We consider a smooth, open and bounded domain $\Omega \subset \R^\delta~(\text{here}~ \delta=2)$. 
In this computational domain, a lower dimensional fracture can be indicated by $\mathcal{C}\subset \R^{\delta-1}$. 
In the following, Dirichlet boundaries conditions indicated as $\partial\Omega_D :=
\partial\Omega$, and Neumann boundaries conditions are given on $\partial_N \Omega := \Gamma_N \cup \partial\mathcal{C}$, where $\Gamma_N$ 
is the outer boundary of $\Omega$ and $\partial\mathcal{C}$ is the
crack boundary. The geometric setup including notations is illustrated in Figure \ref{Figure1}a. 
The surface fracture $\mathcal{C}$ is estimated in
$\Omega_F\subset\Omega \subset \mathbb{R}^\delta$. A region without
any fracture (i.e., an intact region) is indicated by $\Omega_R:=\Omega \backslash \Omega_F\subset\Omega \subset \mathbb{R}^\delta$ such that $\Omega_R\cup\Omega_F=\Omega$ and $\Omega_R\cap\Omega_F=\varnothing$.

The variational phase-field formulation is a thermodynamically consistent framework to compute
the fracture process. Due to the presence of the crack surface, we
formulate the fracture problem as a two-field problem including the
displacement field $\bm u(\bm{x})\colon \Omega\to\mathbb{R}^\delta$
and the crack phase-field $d(\bm{x})\colon\Omega\to[0,1]$. The crack
phase-field function $d(\bm{x})$ interpolates between $d=1$, which
indicates undamaged material, and $d=0$, which indicates a fully broken material phase.

For stating the variational formulations, the spaces
\begin{align}\label{space1}
\bm{V}:=& ~\{ {\bf H}^1(\Omega)^\delta:\bm u=\bar{\bm u}\; \mathrm{on}
\; \partial\Omega_D  \},\\
W:=& ~\text{H}^1(\Omega) ,\\
W_{in} :=&~ \{ d \in \text{H}^1(\Omega)^{\delta-1} : 0 \leq
d\leq d^{old} \}
\end{align}
are used.
Herein, $W_{in}$ denotes a closed, non-empty and convex set which is a
subset of the linear function space $W=\text{H}^1(\Omega)$ (see e.g.,
\cite{KiStam00}).

\subsection{Variational formulation for the \noii{isotropic} mechanical contribution}\label{Section21}

In the following, a variational setting for quasi-brittle fracture in
bulk materials with small deformations is formulated. To formulate the
bulk free energy stored in the material, we define the first and second invariants as
\begin{equation}\label{eq1}
I_1(\bm{\varepsilon})=tr(\bm{\varepsilon}),
\qquad I_2(\bm{\varepsilon})=tr(\bm{\varepsilon}^2),
\end{equation}
with the second-order infinitesimal small strain tensor defined as
\begin{equation}\label{eq2}
\bm{{\varepsilon}} = \nabla_s \bm{u} = sym[ \nabla \bm{u} ].
\end{equation}
\noii{The isotropic scalar valued free-energy function reads
	\begin{equation}\label{eq6}
	\widetilde{\Psi}\big(I_1(\bm{\varepsilon}),I_2(\bm{\varepsilon})\big)
	:=(\frac{K}{2})I^2_1(\bm{\varepsilon})
	-\mu~\Big(\frac{I^2_1(\bm{\varepsilon})}{3}- I_2(\bm{\varepsilon})\Big)
	\quad\text{with}\quad K>0\quad\text{and}\quad\mu>0,
	\end{equation}
	where $K= \lambda+\frac{2}{3}\mu$ is the bulk modulus}. 
A stress-free condition for the bulk energy-density function
requires $\widetilde{\Psi}\big(I_1(\bm{0}),I_2(\bm{0})\big)=0$. Hence,
the bulk free-energy functional including the stored internal energy
and the imposed external energy is
\begin{equation}
\begin{aligned}
\mathcal{E}_{bulk}(\bm{u})=\int_{\Omega_C}\widetilde{\Psi}(\bm\varepsilon)
\mathrm{d}{\bm{x}}-\int_{{\partial_N\Omega_C }} {{\bm\tau}} \cdot \bm u\,\mathrm{d}s
\label{eq8}
\end{aligned}
\end{equation}
where $\bm\tau$ is the imposed traction traction vector on 
${\partial_N\Omega_C := \Gamma_N \cup \mathcal{C}}$ and the body-force is neglected. 

Following \cite{FraMar98}, we define the total energetic functional
which includes the stored bulk-energy functional and fracture
dissipation as
\begin{equation}
\mathcal{E}(\bm{u},\mathcal{C})= \mathcal{E}_{bulk}(\bm{u})
+ G_c \mathcal{H}^{\delta-1}  (\mathcal{C}),
\label{eq9}
\end{equation}
where $G_c$ is the so called the Griffith's critical elastic-energy
release rate. Also, $\mathcal{H}^{\delta-1}$ refers to the
$(\delta-1)$-dimensional Hausdorff measure (see e.g.\
\cite{BourFraMar00}). Following \cite{BourFraMar00},
$\mathcal{H}^{\delta-1}$ is regularized (i.e.\ approximated) by the
crack phase-field $d(\bm{x})$ (see e.g.\ \cite{BourFraMar00}). Doing
so, a second-order variational phase-field formulation is employed;
see Section \ref{Section22}. Additional to that, a second-order stress
degradation state function (\textit{intacted-fractured transition}
formulation) is used as a monotonically decreasing function which is
lower semi-continuous order; see Section \ref{Section233}.

\subsection{Crack phase-field formulation in a regularized setting}\label{Section22}
Let us represent a regularized (i.e.,\ approximated) crack surface for
the sharp-crack topology (which is a Kronecker delta function) thorough the exponential function $d(\bm{x})=1-\exp^{ - \vert \bm{x} \vert /l }$, which satisfies $d(\bm{x})=0$ at $\bm{x}=0$ as a Dirichlet boundary condition and $d(\bm{x})=1$ as $\bm{x}\rightarrow \pm \infty$. This is explicitly shown in Figure \ref{Figure1}b for different length scales. \noii{Here, $\bm{x}$ is a position variable in the Cartesian coordinate system, meaning $\bm{u}$ and ${d}$ have a certain value at each position within the geometry. The first observation through the explicit formulation is that, the crack phase-field $d$ constituting a smooth transition zone dependent on the regularization parameter $\ell$. In engineering or physics, $\ell$ is often a so-called characteristic length-scale parameter. This may be justified since this zone weakens the material and is a physical transition zone from the unbroken material to a fully damaged state. In practice, \new{choices such as} $\ell=2h$ or $\ell=4h$ \new{are often employed.}}
 \noii{Following \cite{pham2011gradient,MieWelHof10b}, a regularized \textit{crack surface energy functional} for the second term in Eq. \ref{eq9} reads} 
\begin{equation}\label{eq11}
G_c \mathcal{H}^{\delta-1}  (\mathcal{C}):=G_c\int_\Omega\gamma_\ell(d, \nabla d) \,\mathrm{d}{\bm{x}}
\quad \text{with}\quad
\gamma_{\ell}(d,\nabla d):=  \frac{\ell}{2\ell} {(1-d)^2} + \frac{\ell}{2}
\nabla d \cdot \nabla d
\end{equation}
based on the crack surface density function \noii{$\gamma_{\ell}(d,\nabla d)$} per unit volume of the solid. The equation is so-called AT-2 model because of the quadratic term in PDE.

\begin{figure}[!ht]
	\centering
	\includegraphics[width=15cm,height=5.5cm]{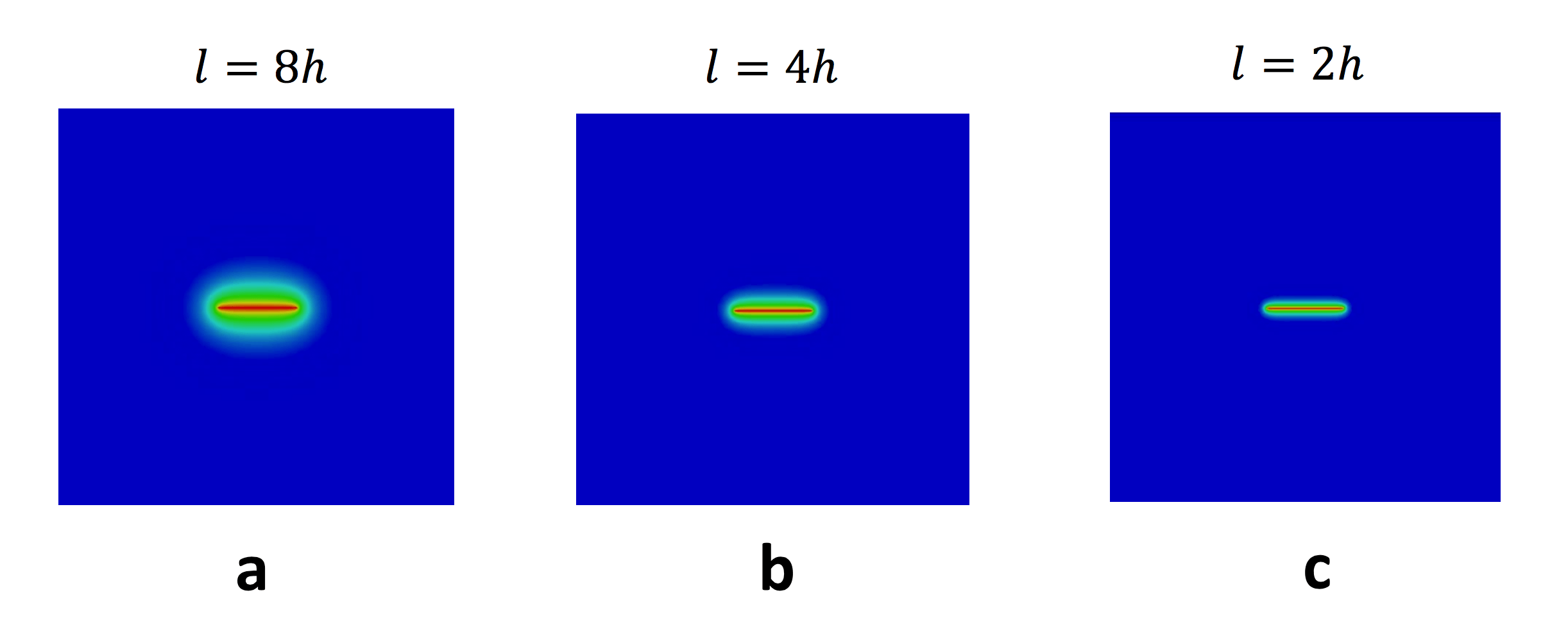}
	\caption{Effect of different length scales on the crack
		phase-field resolution as calculated by the minimization
		problem in Eq.~\ref{min_d} such that $\ell_a>\ell_b>\ell_c$.}
	\label{Figure11}
\end{figure}

We set sharp crack surfaces as Dirichlet boundary conditions in
$\mathcal{C}\subset\Omega$. Hence, the crack phase-field \noii{$d(\bm{x},t)$} is obtained from the minimization of the regularized crack density function as
\begin{align}
\label{min_d}
d(\bm{x})=\underset{d(\bm{x})\in W_{in} \;
	\text{with} \; d(\bm{x})=0 \; \forall \bm{x}\in \mathcal{C}}{\mathrm{argmin}} \int_\Omega\gamma_l(d, \nabla d) \,\mathrm{d}{\bm{x}}.
\end{align} 
Figure \ref{Figure11} gives the numerical solution that arises from
the minimization Eq.~\ref{min_d} and demonstrates the effect of
different regularized length scales on the numerical
solution. Clearly, a smaller length scale leads to a narrower
transition zone (see Figure \ref{Figure11}c). That is also in
agreement with the crack phase-field profile shown in Figure
\ref{Figure1}b.

\subsection{Strain-energy decomposition for the bulk free-energy}\label{Section23}

Fracture mechanics is the process which results in the compression
free state. As a result, a fracture process behaves differently in the
\textit{positive phase} and in \textit{negative phase}, see e.g. \cite{aldakheel2014towards}. In the following, an additive split for the strain energy density function to distinguish  the \noii{positive} and \noii{negative}  phases is used. Instead of dealing with a full linearized strain tensor $\bm\varepsilon(\bm u)$, the
additive decomposition
\[
\bm\varepsilon(\bm u)=\bm\varepsilon^{+}(\bm u)+\bm\varepsilon^{-}(\bm
u)
\quad \text{with}\quad
\bm\varepsilon^{\pm}(\bm
u):=\sum_{i=1}^{\delta} \langle\varepsilon_i\rangle^{\pm} {\textbf{N}_i} \otimes {\textbf{N}_i},
\]
of the strain tensor based on its eigenvalues is used \cite{MieWelHof10b,HeWheWi15}.  Herein,
$\langle x \rangle_{\pm} := \frac{ x {\pm} |x|}{2}$ refers to the a
Macaulay brackets for $x \in \R^\pm$. Furthermore,
$\bm\varepsilon^{+}$ and $\bm\varepsilon^{-}$ refer to the positive and
negative parts of the strain, respectively. The $\{\varepsilon_i\}$
are the principal strains (i.e., the eigenvalues of the
$\bm\varepsilon(\bm u)$) and the $\{\textbf{N}_i\}$ are the principal
strain directions (i.e., the eigenvectors of the $\bm\varepsilon(\bm
u)$).  To determine the positive and negative parts of total strain
${\bm {\varepsilon}}$, a positive-negative fourth-order projection
tensor is
\begin{equation}\label{eq16}
\mathbb{P}^\pm_{\bm {\varepsilon}}:=\frac{\partial \bm {\varepsilon}^\pm}{\partial \bm {\varepsilon}}=\frac{\partial \big(\displaystyle\sum_{i=1}^{\delta} \langle\varepsilon_i\rangle^{\pm}  {\textbf{N}_i} \otimes {\textbf{N}_i}\big)}{\partial \bm {\varepsilon}},
\end{equation}
such that the fourth-order projection tensor $\mathbb{P}^\pm_{\bm
	{\varepsilon}}$ projects the total linearized strain ${\bm
	{\varepsilon}}$ onto its positive-negative counterparts, i.e.,
$\bm {\varepsilon}^{\pm}=\mathbb{P}^\pm_{\bm {\varepsilon}}:\bm
{\varepsilon}$. Hence an additive formulation of the strain-energy
density function consisting of the positive and the negative parts reads
\begin{equation}\label{eq17}
{\Psi}\big(I_1(\bm{\varepsilon}),I_2(\bm{\varepsilon})\big):=\underbrace{\widetilde{\Psi}^{+}\big(I^{+}_1(\bm{\varepsilon}),I^{+}_2(\bm{\varepsilon})\big)}_{\text{tension term}}+\underbrace{\widetilde{\Psi}^{-}\big(I^{-}_1(\bm{\varepsilon}),I^{-}_2(\bm{\varepsilon})\big)}_{\text{compression term}}.
\end{equation}
Here, the scalar valued principal invariants in the positive and negative modes are
\begin{equation}\label{eq18}
I_1^{\pm}(\bm{\varepsilon}):=\langle{I_1(\bm{\varepsilon})}\rangle_{\pm}, \quad I^{\pm}_2(\bm{\varepsilon}):=I_2(\bm{\varepsilon}^{\pm}).
\end{equation}

\noii{Here, the first positive/negative invariant $I_1(\bm{\varepsilon})$ is strictly related to the tension/compression mode, respectively, meaning that if $tr(\bm{\varepsilon})>0$ requires that we are in tension mode otherwise we are in compression state. The second invariant $I_2(\bm{\varepsilon})$ is mainly due to the positive and negative eigenvalue of the strain tensor, where its positive value requires that we are either in shear or in tension mode otherwise it is in compression. Thus, we distinguished between tension/compression and also a isochoric mode of our constitutive model, and only the positive part of the energy is degraded.}

\subsection{Energy functional for the isotropic crack topology}\label{Section233}

Due to the  physical response of the fracture process, it is assumed
that the degradation of the bulk material due to the crack propagation
depends only on the tensile and isochoric counterpart of the stored
bulk energy density function. Thus, there is no degradation of the
bulk material in negative mode, see \cite{MieWelHof10b}. Hence, the
degradation function denoted as $g(d_+)$ acts only on the positive
part of bulk energy given in Eq.~\ref{eq17}, i.e.,
\begin{equation}
g(d_+):=d_+^2 ,\qquad g\colon [0,1]\rightarrow [0,1].
\label{eq20}
\end{equation}
This function results in degradation of the solid during the evolving crack phase-field
parameter $d$. Due to the transition between the intact region and the
fractured phase, the degradation function has flowing properties, i.e.,
\begin{equation}
g(0)=0,\quad g(1)=1,\quad g(d)>0 \;\; \text{for} \;\;d>0, \quad g'(0)=0,\quad g'(1)>0.
\label{eq200}
\end{equation}

Following \cite{MieWelHof10b}, the small residual scalar
$0<\kappa\ll1$ is introduced to prevent numerical instabilities. It is
imposed on the degradation function, which now reads
\begin{equation}
g(d_+):=(1-\kappa)d_+^2 + \kappa,\qquad g\colon [0,1]\rightarrow [0,1).
\label{eq2000}
\end{equation}

The stored bulk density function is denoted as $w_{bulk}$.  Together
with the fracture density function $w_{frac}$, it gives the the total density function
\begin{equation}
\label{total_free_energy}
w(\bm{\varepsilon}, d, \nabla d) = w_{bulk} (\bm{\varepsilon}, d) +  w_{frac} (d, \nabla d),
\end{equation}
with
\begin{equation}
\label{eq116}
w_{bulk} (\bm{\varepsilon}, d)=g(d_+)\widetilde{\Psi}^{+}\big(I^{+}_1(\bm{\varepsilon}),I^{+}_2(\bm{\varepsilon})\big)+\widetilde{\Psi}^{-}\big(I^{-}_1(\bm{\varepsilon}),I^{-}_2(\bm{\varepsilon})\big),
\end{equation}
\begin{equation*}
w_{frac}(d, \nabla d) = G_c \gamma_l(d, \nabla d).
\end{equation*}

\begin{form}[Energy functional for isotropic crack topology]
	\textit{	\label{form_2}
		We assume that \new{$K$} and $\mu$ are given as well as initial conditions $\bm u_0=\bm u(\bm{x},0)$ and  $d_0=d(\bm{x},0)$. For the loading increments 
		$n \in \{1,2,\ldots, N\}$, find $\bm u:=\bm u^n\in V$ and
		$d:=d^n\in W_{in}$ such that the functional
		\begin{align*}
		\mathcal{E} ( {\bm u},d)&=\mathcal{E}_{bulk}(\bm u,d_+,\chi
		)+\; \mathcal{E}_{frac}(d)+\;\mathcal{E}_{ext}(\bm u)\\[1mm] 
		&=\underbrace{\int_\Omega g(d_+)\;\widetilde{\Psi}^{+}(I^{+}_1,I^{+}_2) 
			+\widetilde{\Psi}^{-}(I^{-}_1,I^{-}_2)\; \mathrm{d}{\bm{x}}}_{\text{bulk term}}
		+ G_c\underbrace{\int_\Omega\gamma_{l}(d, \nabla d)\mathrm{d}{\bm{x}}}_{\text{fracture term}}
		-\underbrace{\int_{\partial_N\Omega } {\bm {\bar\tau}} \cdot \bm u\,\mathrm{d}s}_{\text{external load}},
		\end{align*}
		is minimized.}
\end{form}

Herein, to make sure that phase-field quantity~$d$ lies in the
interval $[0, 1]$, we define $d_+$ to map negative values of $d$ to
positive values. In Formulation \ref{form_2}, the stationary points of
the energy functional are determined by the first-order necessary
conditions, namely the Euler-Lagrange equations, which can be found by
differentiation with respect to ${\bm u}$ and $d$.
\begin{form}[Euler-Lagrange equations]
	\label{form_3}
	\textit{	Let \noii{$K>0$, $\mu>0$} be given as well as the initial
		conditions $\bm u_0=\bm u(\bm{x},0)$ and
		$d_0=d(\bm{x},0)$. For the loading increments $n \in
		\{1,2,\ldots, N\}$, find $\bm u:=\bm u^n\in V$ and $d:=d^n\in
		W_{in}$ such that
		\begin{equation}
		\begin{aligned}
		&{\mathcal E}_{\bm u}(\bm u, d;\delta\bm u)=
		\int_\Omega g(d_+) {\bm {\widetilde{\sigma}}^{iso,+}}(\bm u): {\bm \varepsilon}(\delta\bm u)\mathrm{d}{\bm{x}}+ 
		\int_\Omega {\bm {\widetilde{\sigma}}^{iso,-}(\bm u)}: {\bm \varepsilon}(\delta\bm u)\mathrm{d}{\bm{x}}-\int_{\partial_N\Omega } {\bm {\bar\tau}} \cdot {\delta\bm u}\,\mathrm{d}s= 0 \qquad\forall {\delta\bm u}\in V, \\
		&{\mathcal E}_d(\bm u,d;\delta d-d) =(1-\kappa)\int_\Omega 2d_+ \widetilde{D}. (\delta d-d)\mathrm{d} {\bm{x}}\\
		&\qquad \qquad  \qquad  \qquad+ G_c \int_\Omega \left( \frac{1}{\ell}(d-1)\cdot(\delta d-d)
		+ \ell \nabla d\cdot \nabla(\delta d-d)\right)\mathrm{d}{\bm{x}}\geq 0 \qquad\forall
		\delta d\in W \cap L^{\infty}.
		\end{aligned}
		\end{equation}}
\end{form}

Herein, ${\mathcal E}_{\bm u}$ and ${\mathcal E}_d$ are the first
directional derivatives of the energy functional ${\mathcal E}$ given
in Formulation \ref{form_2} with respect to the two fields, i.e., $\bm u$ and $d$, respectively. Also, $\widetilde{D}$ is a crack driving state function which
depends on a state array of strain- or stress like quantities and $\delta\bm u\in\{ {\bf H}^1(\Omega)^2: \delta\bm u=\bm 0 \; \mathrm{on} \; \partial\Omega_D \}$ is the deformation test function and $\delta d\in H^1(\Omega)$ is the phase-field test function. 

Furthermore, the second-order constitutive stress tensor with respect to Eq. \ref{eq116} reads
\begin{equation}
\begin{aligned}
&{\bm \sigma}(\bm{\varepsilon},d ):=\frac{\partial w_{bulk}(\bm{\varepsilon}, d)}{\partial {\bm \varepsilon}} =g(d_+)\frac{\partial \widetilde{\Psi}^{+}}{\partial \bm \varepsilon}+\frac{\partial \widetilde{\Psi}^{-}}{\partial \bm \varepsilon}
=g(d_+)~{\bm {\widetilde{\sigma}}^{+}}+{\bm {\widetilde{\sigma}}^{-}},
\label{eq21}
\end{aligned}
\end{equation}
with 
\noii{\begin{equation}
	\bm {\widetilde{\sigma}}^{\pm}(\bm{\varepsilon})
	:=K~I_1^{\pm}(\bm{\varepsilon}) -2\mu~\Big(\frac{1}{3}I_1^{\pm}(\bm{\varepsilon})  {\textbf{I}} 
	- 2 \bm\varepsilon_\pm\Big)
	\quad\text{with}\quad K>0\quad\text{and}\quad\mu>0.
	\label{eq24}
	\end{equation}}

\subsection{Crack driving forces for brittle failure}\label{Section25}

Following \cite{aldakheel2018,teichtmeister2017}, we determine the crack driving state function to couple between two PDEs. Hence, crack driving state function acts as a right hand side for the phase-field equation. To formulate the crack driving state function, we consider the crack
irreversibility condition, which is the inequality constraint $\dot d
\le 0$ imposed on our variational formulation. The first variation of the total pseudo-energy density with respect to the crack phase-field given in (\ref{total_free_energy}) reads
\begin{equation}\label{29_H1}
- \delta_{d} w(\bm{\varepsilon}, d, \nabla d)=(\kappa-1)2d_+ \big[\widetilde{\Psi}^{+}
\big] - G_c \delta_d\gamma_\ell(d, \nabla d)
\geq 0.
\end{equation}
Herein, the functional derivative of $\gamma_l(d,\nabla d)$ with respect to $d$ is 
\begin{equation}\label{eq29_H5}
\int_\Omega \delta_{d} \gamma_{\ell}(d,\nabla d) \mathrm{d}{\bm{x}}=\int_\Omega \frac{1}{\ell}[(d-1)-\ell^2 \Delta d] \mathrm{d}{\bm{x}}.
\end{equation}
Maximization the inequality given in Eq. \ref{29_H1} with respect to the time history $s\in [0,t_n]$ reads
\begin{equation}\label{29_H2}
(\kappa-1)2d_+\max_{s \in [0,t_n]}
\big[\widetilde{\Psi}^{+}
\big] =G_c \delta_d\gamma_\ell(d, \nabla d).
\end{equation} 
We multiply  Eq.\ \ref{29_H2} by $\frac{l}{G_c}$.  Then Eq.\
\ref{29_H2} can be restated as
\begin{equation}\label{eq29_H3}
(\kappa-1)2d_+\mathcal{H} =\ell \delta_d\gamma_l \quad
\text{if} \quad \mathcal{H}:= \max_{s\in [0,t_n]}\widetilde{D}
\quad \text{with} \quad \widetilde{D}:=\frac{\ell{\widetilde{\Psi}}^{+}}{G_c}.
\end{equation}

Here, $\mathcal{H}:=\mathcal{H}(\bm{\varepsilon},t)$ denotes a
positive crack driving force that is used as a history field from
initial time up to the current time. Note that the crack driving state
function $\widetilde{D}$ is affected by the length-scale parameter $\ell$
and hence depends on the regularization parameter. 

\begin{algorithm}[ht!]
	{\bf Input:} \hspace{0.2cm}$\bullet$ loading data $(\bar{\bm u}_{n},\bar{\bm t}_{n})$ on $\partial \Omega_D\subset\partial\Omega$, \\ 
	\hspace{1.26cm}$\bullet$ solution $(\bm u_{n-1},d_{n-1})$ from time step $n-1$. \\ [0.1cm]  
	
	Initialization of alternate minimization scheme ($k=1$):  \hspace{0.2cm}$\bullet$ set \var{FLAG}:=true
	
	\While{\var{FLAG}}{		
		\quad  1.  given $\bm u^k$, solve ${\mathcal E}_d(\bm u^k,d;\delta d)=0$ for $d$, set $d:=d^k$,\\
		\quad 2. given $d^{k-1}$, solve ${\mathcal E}_{\bm u}(\bm u,d^{k-1};\delta\bm u)=0$ for $\bm u$, set $\bm u:=\bm u^k$,\\
		\quad 3. define alternate minimization residual for the obtained pair $(\bm u^k,d^k)$
		\begin{align}
		\mathrm{Res}_\mathrm{Stag}^k:=|{\mathcal E}_{d}(\bm u^k,d^k;\delta d)|+|{\mathcal E}_{\bm u}(\bm u^k,d^k;\delta\bm u)| \quad \forall \; \delta\bm u\in{\bf V}\quad\delta d\in{W},
		\end{align}

		\quad 4.$~\textbf{if}~$ $\mathrm{Res}_\mathrm{Stag}^k\leq\texttt{TOL}_\mathrm{Stag}$ \textbf{then}\\ 
		
		\qquad\quad $\bullet$ set $(\bm u^k,d^k):=(\bm u_n,d_n)$\\
		\qquad\quad $\bullet$  $\var{FLAG}$:=false\\
		
		\qquad \textbf{else}\\
		
		\qquad\quad $k+1\rightarrow k$\\
		
		\qquad \textbf{end if}	
	}
	\caption{Alternate minimization scheme for Formulation (\ref{form_3_final}) at a fixed loading step $n$.} 
	\label{TableStag}	
\end{algorithm}

\begin{form}[Final Euler-Lagrange equations]
	\textit{	\label{form_3_final}	
		Let us assume that \noii{$K>0$, $\mu>0$} are given as well as the
		initial condition $\bm u_0=\bm u(\bm{x},0)$ and
		$d_0=d(\bm{x},0)$. For the loading increments $n \in
		\{1,2,\ldots, N\}$, find $\bm u:=\bm u^n\in V$ and $d:=d^n\in
		W$ such that
		\begin{equation}\label{AT-2}
		\begin{aligned}
		&{\mathcal E}_{\bm u}(\bm u, d_+;\delta\bm u)=
		\int_\Omega g(d_+) {\bm {\widetilde{\sigma}}^{+}_{\bm
				\varepsilon}}(\bm u): {\bm \varepsilon}(\delta\bm u)\,\mathrm{d}{\bm{x}}+ 
		\int_\Omega {\bm {\widetilde{\sigma}}^{-}_{\bm \varepsilon}(\bm u)}: {\bm \varepsilon}(\delta\bm u)\,\mathrm{d}{\bm{x}}-\int_{\partial_N\Omega } {\bm {\bar\tau}} \cdot {\delta\bm u}\,\mathrm{d}s= 0 \qquad\forall {\delta\bm u}\in \boldsymbol{V}, \\
		&{\mathcal E}_d(\bm u,d;\delta d) =(1-\kappa)\int_\Omega 2d_+ \mathcal{H} \delta d \,\mathrm{d}{\bm{x}}+ \int_\Omega \Big((d-1) \delta d
		+ \ell^2 \nabla d\cdot\nabla\delta d\Big)\mathrm{d}{\textbf{x}}= 0 \qquad\forall
		\delta d\in W.
		\end{aligned}
		\end{equation}}
\end{form}

The multi-field problem given in Formulation (\ref{form_3_final}) depending on $\bm{u}$ and $d$ implies alternately fixing $\bm u$ and $d$, which is a so called alternate minimization scheme, and then solving the corresponding equations until convergence. The alternate minimization scheme applied to the Formulation (\ref{form_3_final}) is summarized in Algorithm \ref{TableStag}.

\subsection{\noii{The influence of the $\kappa$ on the stress-strain curve}}\label{Section26}

\begin{figure}[!t]
	\centering
	\subfloat{\includegraphics[clip,trim=7cm 0cm 4cm 0cm, width=8.6cm,height=5.5cm]{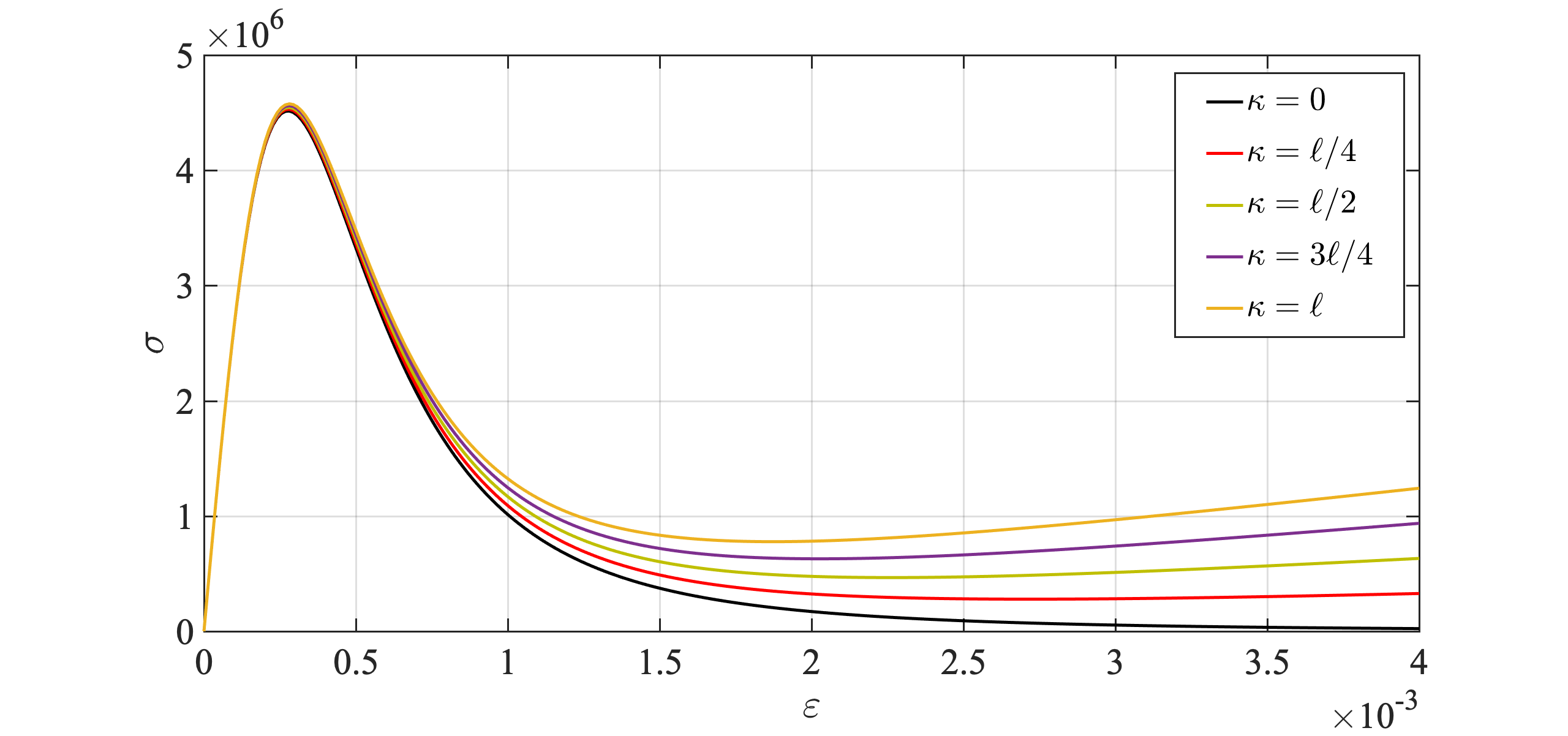}}
	\subfloat{\includegraphics[clip,trim=7cm 0cm 4cm 0cm, width=8.6cm,height=5.5cm]{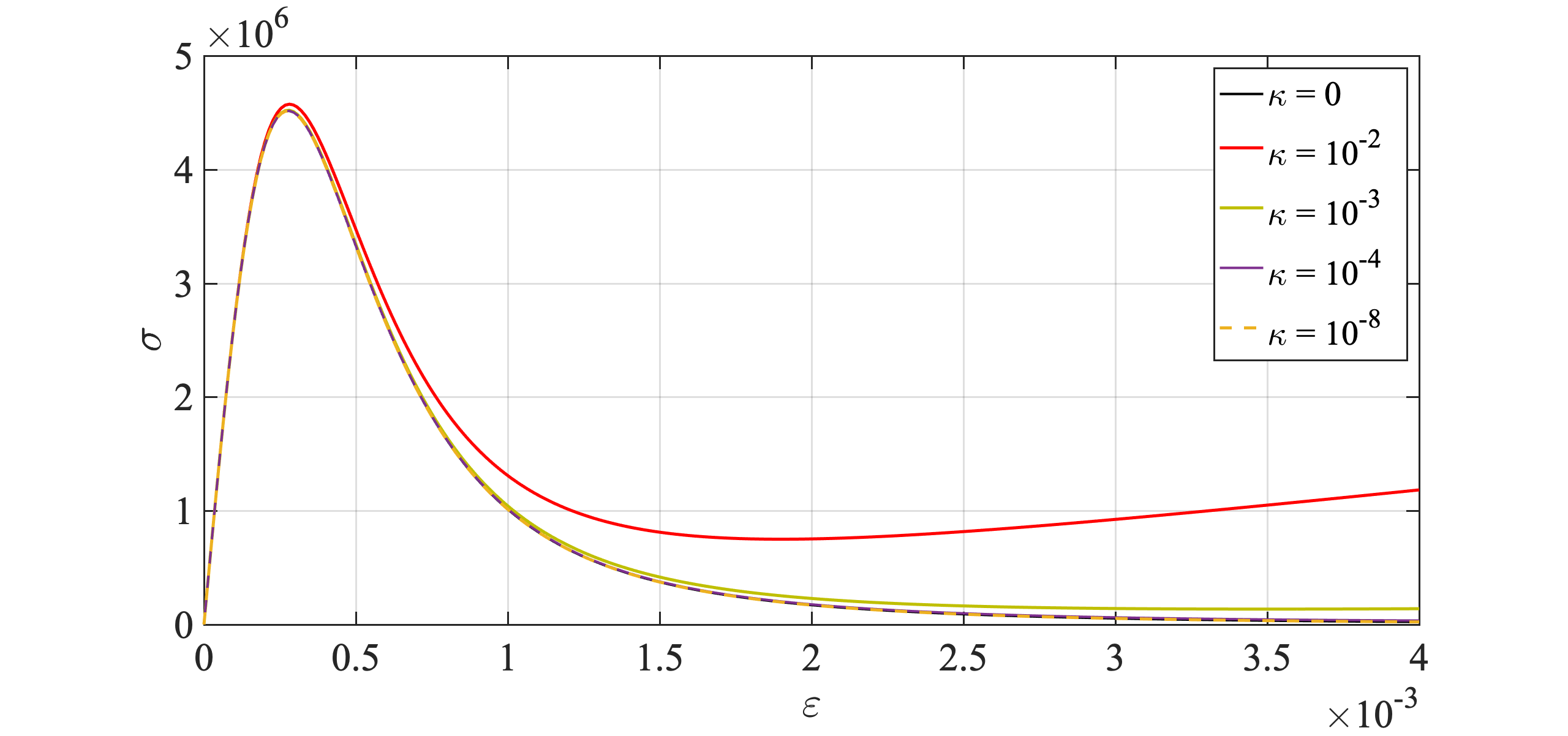}}
	\caption{The influence of the $\kappa$ on the stress-strain curve; Left plot represent $\kappa=\kappa(\ell$ and right plot presents $\kappa$ as a numerical parameter which is sufficiently small.}
	\label{Figure111}
\end{figure}

\new{In this part, the influence of the $\kappa$ on the stress-strain curve is taken into account. Following \cite{miehe2016phase}, the homogeneous solution at the quasi-static stationery state of the phase-field partial differential equation in the loading case takes the following form 
\begin{equation}\label{eqA31}
d_{homo}=\frac{1}{1+2(1-\kappa)\widetilde{D}} \;\in[0,1],
\end{equation} 
which results from the free Laplacian operator $\Delta (\bullet)=0$ assumption in Eq. \ref {29_H2} without any source terms (zero left-hand sides). Here, the crack driving state function $\widetilde{D}$ is given in Eq. \ref{eq29_H3}. Because, we are in the elastic limit, prior to the onset of fracture, then no split is considered. We now aim to relate a stress state $\sigma$ with the isotopic phase-field formulation. To do so, a non-monotonous function in the one-dimensional setting for the degrading stresses takes the following form by
\begin{equation}\label{eqA32}
\sigma=g(d)\widetilde{\sigma}\\
=\big(\frac{(1-\kappa)}{\big(1+2(1-\kappa)\widetilde{D}\big)^{2}}+\kappa \big)~ E\varepsilon.
\end{equation}
To see the influence of the $\kappa$ in Eq. \ref{eqA32}, the concrete material is considered.  Following, \cite{marigo2016overview} for a concrete material which has a brittle response, a typical values for material parameters reads,
\begin{equation}\label{mat_conc}
	E=29\,GPa,\quad \sigma_c=4.5\,MPa \quad \text{and} \quad G_c=70\,N/m.
\end{equation}
We set $\ell=0.0105~m$. Thus, we can do a plot for the stress-strain curve through Eq. \ref{eqA32} by considering the material set given above. Figure \ref{Figure111} shows the effect of stress state for different strain loading. The black curve represents the stress-strain curve while $\kappa=0$. Evidently, it can be grasped through Figure \ref{Figure111} with $\kappa=0$ the $\sigma_c$ is exactly $\sigma_c=4.5\,MPa$ as it is required for the concrete material, see \cite{marigo2016overview}. If we consider $\kappa\neq0$ as a function of characteristic length-scale, see Figure \ref{Figure111} left, we can observe a good agreement with $\kappa=0$ up to the peak point while after some strain value it becomes different as $\kappa$ changes. Unfortunately, we can not observe any converged response if we consider $\kappa$ as a function of $\ell$. In contrast, if we chose $\kappa$ sufficiently small, see Figure \ref{Figure111} right, as much as $\kappa$ reduced, in here less than $\kappa\le10^{-4}$, we observed a very identical response with $\kappa=0$, thus it behaves as numerical parameters rather than material parameters.} 

\section{Stochastic model for Bayesian inversion}
\label{section3}
In this section, we explain how we use Bayesian inversion to identify parameters.
Then, we introduce a computationally effective numerical technique to estimate the unknown parameters. 

%

In the phase-field model, the uncertaint{ies arise} from the \new{elasticity parameters} including the shear modulus $\mu$ and \new{the bulk modulus $K$} as well as Griffith's critical elastic energy release rate (material stiffness
parameter) $G_c$,  which are assumed to be random fields. Specifically, we represent
the \new{parameters} uncertainty (\blue{spatial variability}) by a spatially-varying \new{log-normal} random field. 

\blue{The Karhunen-Lo\'eve expansion (KLE) expansion method is used  to reduce the dimensionality of the random field.}
The field
$\Theta$ \blue{representing the elasticity parameters and the energy release rate} can be characterized by its expectation  and covariance using the expansion. Considering the probability density function $\mathbb{P}$, the \new{covariance function} is
\begin{align}
\operatorname{Cov}_{\Theta}(\bm{x},\textbf{y})=\int_{\Omega}\left(\Theta(\bm{x},\omega)-\Theta(\bm{x})\right)\left(\Theta(\bm{y},\omega)-\Theta(\bm{y})\right)\,\text{d}\mathbb{P}(\omega),
\end{align}
which leads to the \new{KL-expansion}
\begin{align}\label{KL}
\Theta (\bm{x},\omega)=\bar{\Theta}(\bm{x})+\sum_{n=1}^{\infty}\sqrt{\psi_n }k_n(\bm{x}) \xi_n (\omega).
\end{align}
Here the first term is the mean value, $k_n$ are the orthogonal eigenfunctions, $\psi_n$ are the corresponding eigenvalues of the eigenvalue problem \cite{ghanem2003stochastic}
\begin{eqnarray}
\int_{D}\text{Cov}_{\Theta} (\bm{x},\bm{y})k_n(\bm{y})~ d\bm{y}=\psi_n k_n(\bm{x}),
\end{eqnarray}
and the $\lbrace \xi_n(\omega)\rbrace$ are mutually uncorrelated random variables satisfying
\begin{eqnarray}
\mathbb{E}[\xi_n]=0,\hspace{1em}\mathbb{E}[\xi_n \xi_m]=\delta_{nm},
\end{eqnarray}
where $\mathbb{E}$ indicates the expectation of the random
variables.

The infinite series can be truncated to a finite series
expansion (i.e., an $N_\mathrm{KL}$-term truncation) by \cite{ghanem2003stochastic}
\begin{align}\label{KL1}
\Theta (\bm{x},\omega)=\bar{\Theta}(\bm{x})+\sum_{n=1}^{N_{_{\text{KL}}}}\sqrt{\psi_n}k_n(\bm{x})\xi_n(\omega).
\end{align}
For the Gaussian random field, we employ an exponential covariance kernel as 
\begin{align}
\operatorname{Cov}_{\Theta} (\bm{x},\bm{y})=\sigma^2\exp \left( -\frac{\|\bm{x}-\bm{y}\|}{\zeta}\right),
\end{align}
where $\zeta$ is the correlation length as well as $\sigma$ is the standard deviation.

For a random field, we describe the parameters using a KL-expansion. Considering the Gaussian field $\xi (\bm{x})$, a log-normal random field can be generated by the transformation $\tilde{\xi}(\bm{x})=\exp(\xi(\bm{x}))$. For instance, for the parameter \new{$K$}, the truncated KL-expansion can be written as 
\begin{align}\label{KL2}
\tilde{\xi}_{\new{K}} (\bm{x},\omega)=\exp\left(\bar{\xi_{\new{K}}}(\bm{x},\omega)+\sum_{n=1}^{N}\sqrt{\psi_n}k_n(\bm{x})\xi_n(\omega)\right).
\end{align}

\subsection{Bayesian inference}
We consider  Formulation \eqref{form_3_final} as the forward model
$\textbf{y}=\mathcal{G}(\,\Theta\,(\bm{x}))$, where $\mathcal{G}\colon
L^2(\Omega)\rightarrow L^2(\Omega)$. The forward model explains the
response of the model to different influential parameters $\Theta$
(here $\mu$, \new{$K$}, and $G_c$).  We can write the statistical
model in the form \cite{smith2013uncertainty}
\begin{align}\label{stat}
\mathcal{M}=\mathcal{G}(\,\Theta)+\varepsilon,
\end{align} 
where $\mathcal{M}$ indicates a vector of observations (e.g.,
measurements). The error term $\varepsilon$ arises from uncertainties
such as measurement error due experimental situations. More precisely,
it is due to the modeling and the measurements and is assumed to have
a Gaussian distribution of the form $\mathcal{N}(0,H)$ with known
covariance matrix $H$.  The error is independent and identically
distributed and is independent from the realizations.  Here, for sake of
simplicity, we assume $H=\sigma^2I$ (for a positive constant
$\sigma^2$).

For a realization $\theta$ of the random field $\Theta$ corresponding
to a realization $m$ of the observations $\mathcal{M}$, the posterior
distribution is given by 
\begin{align}
\pi(\theta|m)=\frac{\pi(m|\theta)\pi_0(\theta)}{\pi(m)}=\frac{\pi(m|\theta)\pi_0(\theta)}{\int_{W_m}\pi(m|\theta)\pi_0(\theta)\,d\theta},
\end{align}
where $\pi_0(\theta)$ is the prior density (prior knowledge) and ${W_m}$ is the
space of parameters $m$ (the denominator is a normalization constant)
\cite{stuart2010inverse}. The likelihood function can be defined as \cite{smith2013uncertainty}
\begin{align}
\label{likelihood}
\pi(m|\theta):=\frac{1}{(2\pi \sigma^2)^{\bar{n}/2}}\exp\left(-\sum_{n=1}^{\bar{n}} \frac{(m_n-\mathcal{G}(\theta))^2}{2\sigma^2} \right).
\end{align}

\begin{figure}[!ht]
	\centering
	\includegraphics[width=9cm,height=5.5cm]{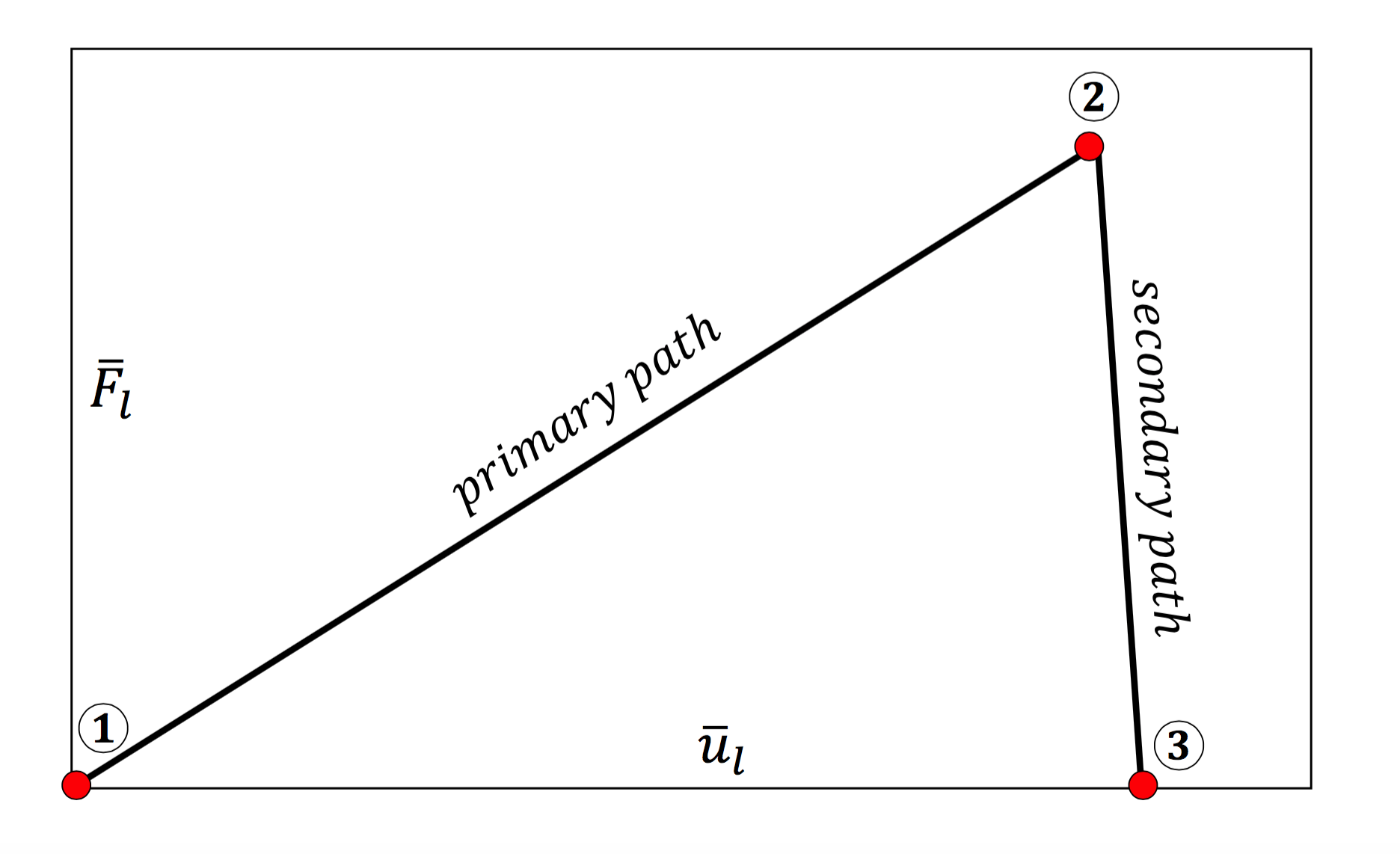}
	\caption{The schematic of load-deflection response for the failure process including primary path (prior to the crack initiation, i.e., between point 1 and 2) and secondary path (during crack propagation, i.e., between point 2 and 3).
	}
	\label{Figure12}
\end{figure}

As an essential characteristic of the phase-field model, the
load-displacement curve (i.e., the global measurement) in addition to
the crack pattern (i.e., the local measurement) are appropriate
quantities to show the crack propagation as a function of time. Figure
\ref{Figure12} indicates the load-displacement curve during the
failure process. Three major points are the following.

\textbf{\raisebox{.5pt}{\textcircled{\raisebox{-.9pt} {1}}}}{\textbf{\textbf{First stable position}.}} This point corresponds to the stationary limit such that we are completely in elastic region ($d(\bm{x},0)=1 \;\forall \bm{x}\in \Omega\backslash\mathcal{C}$).

\textbf{\raisebox{.5pt}{\textcircled{\raisebox{-.9pt} {2}}}}{\textbf{\textbf{First peak point}.}} Prior to this point crack nucleation has occurred and now we have crack initiation. Hence, this peak point corresponds to the critical load quantity such that the new crack surface appears (i.e., there exist some elements which have some support with $d=0$). 

\textbf{\raisebox{.5pt}{\textcircled{\raisebox{-.9pt}
			{3}}}}{\textbf{\textbf{Failure point}.}} At this point, failure
of the structure has occurred and so increasing the load applied to
the material will not change the crack surface anymore.

The interval between point 1 and point 2 in Figure \ref{Figure12}
typically refers to the primary path where we are almost in the
elastic region. The secondary path (sometimes referred to as the
softening damage path) starts with crack initiation occurring at point 2. The
whole process recapitulates the load-deflection curve in the failure
process.

The main aim of solving the inverse problem followed here is to
determine the random field $\Theta$ to satisfy \eqref{stat}.  We
strive to find a posterior distribution of suitable values of the
parameters $\mu$, \new{$K$}, and $G_c$ in order to match the simulated
values (arising from \eqref{AT-2}) with the observations. The
distribution provides all useful statistical information about the
parameter.

\begin{Remark}
	\label{Rem1}
	\textit{ 
		Note that the principal parameters $h$, $\kappa$, and $\ell$ are
		mathematically linked in Formulation  \ref{form_3_final}.
\new{Here, we use $\ell=2h$ and $\kappa$ is sufficiently small which is compatible with Subsection \ref{Section26}. In Section \ref{section4}, the values of $\kappa$ 
in the computations will be specified.
Further, a sufficiently small $h$ is
		chosen to obtain the reference solution. }
	}
	
\end{Remark}


The crack pattern is a time-dependent process (more precisely in a
quasi-static regime, the cracking process is load-dependent), i.e.,
after initiation it is propagated through time. In order to
approximate the parameters precisely, we estimate the likelihood
during all time steps. Therefore, the posterior distribution maximizes
the likelihood function for all time steps, and therefore we have an
exacter curve for all crack nucleation and propagation times.

MCMC is a suitable technique to calculate the posterior
distribution. When the parameters are not strongly correlated, the MH
algorithm \cite{smith1993bayesian} is an efficient computational
technique among MCMC methods. We propose a new candidate (so-called $\theta$, i.e., a value of
$(\mu,\new{K},G_c)$) according to a proposal distribution 
(for instance uniform or normal distributions) and
calculate its acceptance/rejection probability. The ratio indicates
how likely the new proposal is with regard to the current sample. In
other words, by using the likelihood function \eqref{likelihood}, the
ratio determines whether the proposed value is accepted or rejected
with respect to the observation (here the solution of  Formulation \ref{form_3_final} with a very
fine mesh).  As mentioned, fast convergence means that the parameters
are fully correlated. A summary of the MH algorithm is given below.

\begin{algorithm}[ht!]
	\text{Initialization}: set prior data $\theta^0$ and number of samples $N$.\\  
	\For{$i=1:N$}{		
		\quad 1. Propose a new candidate based on the proposal distribution $\theta^{*}\sim\mathcal{K}(\theta^*|\,\theta^{i-1})$.    \\
		\quad  2. Compute the acceptance/rejection probability $\upsilon (\theta^*|\,\theta^{i-1})=\min\left(1, \cfrac{\pi(\theta^*|\,m)}{\pi(\theta^{i-1}|\,m)}\,\cfrac{\mathcal{K}(\theta^{i-1}|\,\theta^*)}{\mathcal{K}(\theta^*|\,\theta^{i-1}))}\right)$.\\
		
		\quad 3. Generate a random number $\mathcal{V}\sim \text{U}\,(0,1)$.\\
		
		\quad 4.$~\textbf{if}~$ $\mathcal{V}<\upsilon~$ \textbf{then}\\ 
		
		\qquad\quad accept the proposed candidate $\theta^*$ and set $\theta^i:=\theta^*$\\
		
		\qquad \textbf{else}\\
		
		\qquad\quad reject the proposed candidate $\theta^*$ and set $\theta^{i}:=\theta^{i-1}$\\
		
		\qquad \textbf{end if}	
	}
	\caption{The Metropolis-Hastings algorithm.} 
	\label{algorithm}	
\end{algorithm}

\section{Bayesian inversion for phase-field fracture}
\label{sec_Bay_for_PFF}
In this key section, we combine the phase-field algorithm from 
Section \ref{section2} with the Bayesian framework presented 
in Section \ref{section3}.

First, we define two sampling strategies as follows:
\begin{itemize}
	\item \textit{One-dimensional}~Bayesian inversion. We first use $N$-samples (according to the proposal distribution) and extract the posterior distribution of the first \new{set e.g., ($\mu^*$, $K^*$)} where other parameter is according to the mean value. Then obtained information is used to estimate the posterior distribution of next unknown (i.e., $G_c^*$). \new{In order to employ the estimated values, the exponential of the estimated parameters is used in the AT-2 model (see Algorithm \ref{algorithm}).}
	
	\item $\textit{Multi-dimensional}$~Bayesian inversion.  A three-dimensional candidate \new{$(\mu^*,K^*,G_c^*)$} is proposed and the algorithm computes its acceptance/rejection probability. 
\end{itemize}
To make the procedures more clarified we explain the multi-dimensional approach in Algorithm \ref{algorithm3}. Clearly, for the one-dimensional setting; for each parameter (e.g., $\theta^*=\new{(\mu^*,K^*)}$), it can be reproduced separately. We will study both techniques in the first example and the more efficient method will be used for other simulations.  

\begin{algorithm}[ht!]
	\For{$i=1:N$}{		
		\vspace{0.1cm}
		\quad 1. Propose  the $i$-th candidate $\theta^*=\new{(\mu^*,K^*,G_c^*)}$ according to the proposal distribution (unform or normal).    \\
		\quad 2. $\bullet$~set \var{FLAG}=true\\ \qquad $\bullet$~set $ n=0$\\
		\While{\var{FLAG}}{
			\vspace{0.1cm}
			\qquad \quad(i) solve the Formulation \eqref{form_3_final} by Algorithm \ref{TableStag} considering   $\texttt{TOL}_\mathrm{Stag}$ and\\
			\qquad \hspace{0.71cm} the proposed candidate $\theta^*$.\\
			\qquad \quad(ii) approximate $(\bm u_n,d_n)$\\
			\qquad \quad(iii) estimate the crack pattern at the loading stage $n$ by 
			\begin{align}
			\bar{F}_n=\int_{\partial \Omega_D}  \bm{n}\cdot{\bm \sigma}\cdot\bm{n}\,\mathrm{d}{\bm{x}}
			\end{align}
			
			\qquad \quad(iv) \textbf{if}~ $\Bigg\{\exists \,d=0~\text{in}~\Omega\backslash\mathcal{C}\Bigg\}~\&~$$\Bigg\{ \|\bar{F}_n \|<\texttt{TOL}_\mathrm{Load} \Bigg\}$\quad \text{or}\quad $n<n_{max}$ ~\textbf{then}\\\
			\qquad \qquad \qquad $\bullet$ set \var{FLAG}=false\\
			\qquad \qquad~	
			\textbf{else}\\
			\qquad \qquad \qquad $\bullet$ set~ $n=n+1$\\
			\qquad \qquad~	
			\textbf{end if}\\
		}				
		\quad 3. Calculate the likelihood function \eqref{likelihood} for $\bar{F}$ (during all $n$-steps, until $\bar{n}$) with respect to  $\theta^*$
		where $m_n$ \\
		\qquad	indicates  the reference value	 at  the $n$-th loading step.\\
		\quad 4. Compute the acceptance/rejection probability $\nu(\theta^*|\,\theta^{i-1})$.\\
		\quad 5. Use Algorithm \ref{algorithm} to determine $\theta^i$  (i.e., $\theta^*$ is accepted/rejected).\\ 
	}
	\caption{The multi-dimensional Bayesian inversion for phase-field fracture.} 
	\label{algorithm3}	
\end{algorithm}

Here, $n_{max}$ is the sufficiently large value that is set by the user. Also, $\texttt{TOL}_\mathrm{Load}$ is a sufficiently small value to guarantee that the crack phase-field model reached to the material failure time. Note, in part (iv) for the while-loop step, the criteria $\|\bar{F}_n \|<\texttt{TOL}_\mathrm{Load}$ in the secondary path (i.e., during crack propagation state) guarantees that reaction force under imposed Dirichlet boundary surface is almost zero. Hence, 
no more force exists to produce a fractured state. We now term this as a complete failure point. But, in some cases, e.g., shear test as reported in \cite{MieWelHof10b}, by increasing the monotonic displacement load, $\bar{F}_n$ is not reached to zero. For this type of problem, if $n<n_{max}$ holds, then phase-field step (i.e., while-loop step) in Algorithm \ref{algorithm3} will terminate. 

The physical aim of using Bayesian inversion in phase-field fracture is adjusting the effective parameters to fit the solution with the reference values (see Remark \ref{Rem1}). 
With (future) experiments (experimental load-displacement until the failure point), these can be used as reliable reference values. 

\section{Numerical examples}
\label{section4}
In this section, we consider {three} numerical test problems 
to determine the unknown parameters using given Bayesian inference. 
Specifically, we {propose:}
\begin{itemize}
	\item Example 1. the single edge notch tension (SENT) test;
	\item Example 2. double edge notch tension (DENT) test;
	\item Example 3. tension test with two voids.
\end{itemize}
The observations can be computed by very fine meshes (here the reference values) as an appropriate replacement of the measurements (see Remark \ref{Rem1}).
Regarding the observational noise, $\sigma^2=1\times10^{-3}$ is assumed. The main aim here is to estimate the effective parameters ($\mu$, $\new{K}$, and $G_c$) in order to match the load-displacement curve with the reference value.  
To characterize the random fields, we can use the KL-expansion with $N_{_{\text{KL}}}=100$ and the correlation length $\zeta=2$ as well.   

\new{In all examples,} the phase-field parameters set by \new{$\kappa = 10^{-8}$}, 
and regularized length scale \new{$\ell = 2{h}$} (respecting the condition $h < l$). The stopping criterion for the iterative Newton method scheme, i.e. the relative
residual norm that is
\begin{align}
\texttt{Residual}:= \| \bm R(\bm x_{k+1}) \|
\leq \texttt{Tol}_\texttt{N-R} \| \bm R(\bm x_{k}) \|,
\end{align}
is chosen to \new{$\texttt{Tol}_\texttt{N-R}=10^{-8}$}. Here, $\bm R$ indicates a discretized setting of weak forms described in Formulation (\ref{form_3_final}). Regarding alternate minimization scheme we set \new{$\texttt{TOL}_\mathrm{Stag}=10^{-4}$} for all numerical examples and \new{$\texttt{TOL}_\mathrm{Load}=10^{-3}$} is chosen to guarantee that we solve the model only until the material failure time. \new{In the examples, the random fields modeled as a log-normal random field.  For the numerical simulations, all variables are discretized by first-order quadrilateral finite elements.}


%

\subsection{Example 1. The single edge notch tension (SENT) test}
\label{SENT}

\begin{figure}[ht!]
	\centering
	\includegraphics[width=15cm,height=7.5cm]{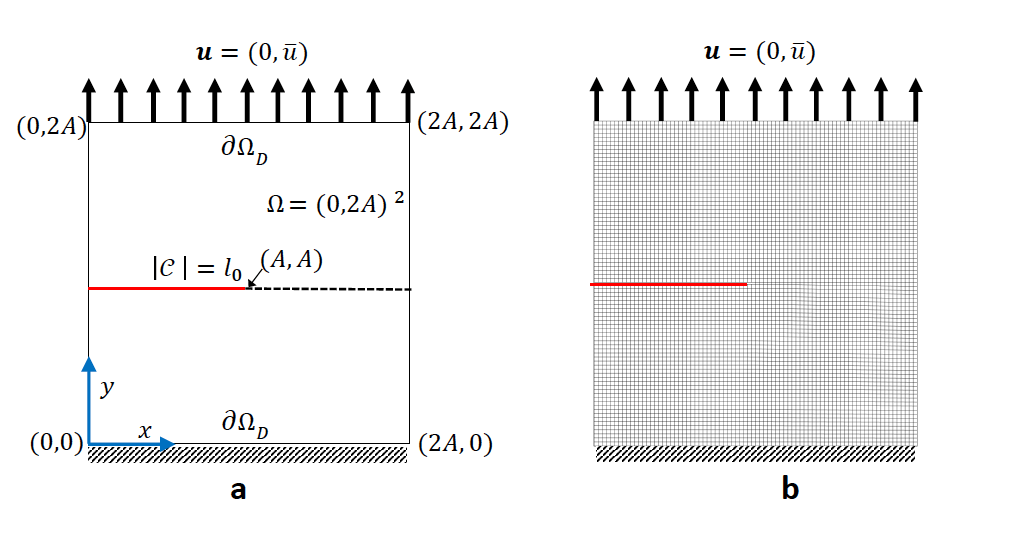}
	\caption{Schematic of SENT (Example 1) (left) and its corresponding mesh with $h=1/80$ (right).}
	\label{schematics1}
\end{figure}

This example considers the single edge notch tension. The specimen is
fixed at the bottom. We have traction-free conditions on both sides. A
non-homogeneous Dirichlet condition is applied at the top. The domain
includes a predefined single notch (as an initial crack state imposed
on the domain) from the left edge to the body center, as shown in
Figure \ref{schematics1}a. We set $A=0.5\;\mathrm{mm}$ hence
$\Omega=(0,1)^2 \mathrm{mm}^2$, hence the predefined notch is in the
$y=A$ plane and is restricted to $0\le|\mathcal{C}|\le A$.  This
numerical example is computed by imposing a monotonic displacement
$\bar{u}=1\times10^{-4}$ at the top surface of the specimen in a
vertical direction.  The finite element discretization corresponding
to $h=1/80$ is indicated in Figure \ref{schematics1}b.

For the shear modulus, we assume the variation range
$(\unit{60}{kN/mm^2},\unit{100}{kN/mm^2})$. Regarding the \new{the bulk modulus}
$K$, the parameter varies between \new{$\unit{140}{kN/mm^2}$ and
	$\unit{200}{kN/mm^2}$}. Finally, we consider the interval between
$\unit{2.1\times 10^{-3}}{kN/mm^2}$ and $\unit{3.3\times
	10^{-3}}{kN/mm^2}$ for $G_c$. \new{Furthermore, we assume that in this example, the variables are spatially constant random variables (they are not random fields).}

We solved the PDE model (Formulation \ref{form_3_final}) with $\mu=\unit{80}{kN/mm^2}$, \new{$K=\unit{170}{kN/mm^2}$}, and  $G_c=\unit{2.7\times 10^{-3}}{kN/mm^2}$ \cite{MieWelHof10b} and the displacement during the time (as the reference solution) with $h=1/320$ was obtained. The main goal is to obtain the suitable values of $\mu,~\new{K}$, and $G_c$ such that the simulations match the reference value.

For this example, we use a uniformly distributed prior distribution
and the \new{uniform} proposal distribution  
\begin{align}
\mathcal{K}(\theta\rightarrow \theta^*):=\frac{1}{\theta_2-\theta_1}\chi[\theta_1,\theta_2](\theta),
\end{align}
where $\chi$ indicates the characteristic function of the interval $[\theta_1,\theta_2]$ \blue{(where $\theta$ denotes a set of parameters)}.

\begin{figure}[t!]
	\centering
	\subfloat{\includegraphics[width=8cm,height=5cm]{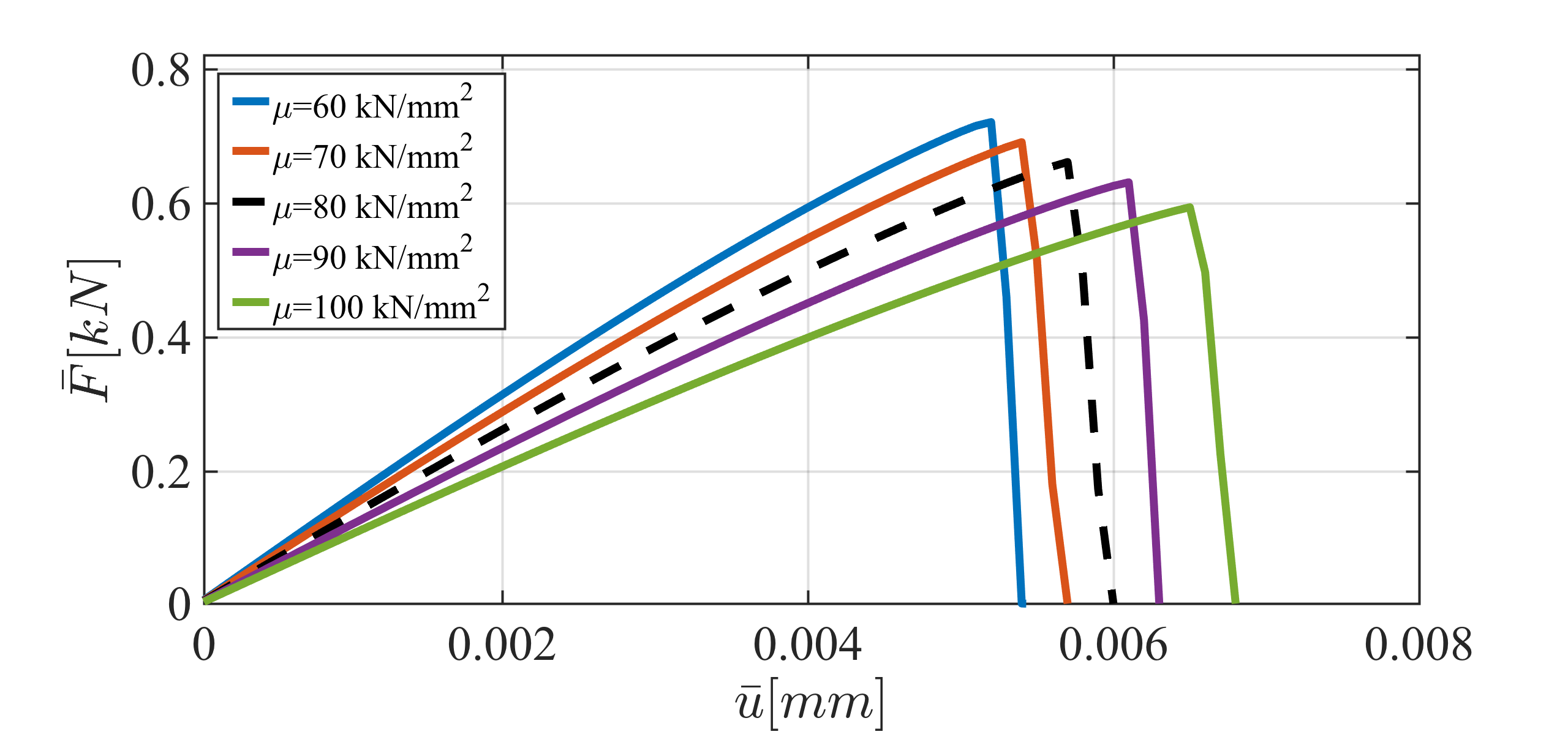}}%
	\hfill 
	\subfloat{\includegraphics[width=8cm,height=5cm]{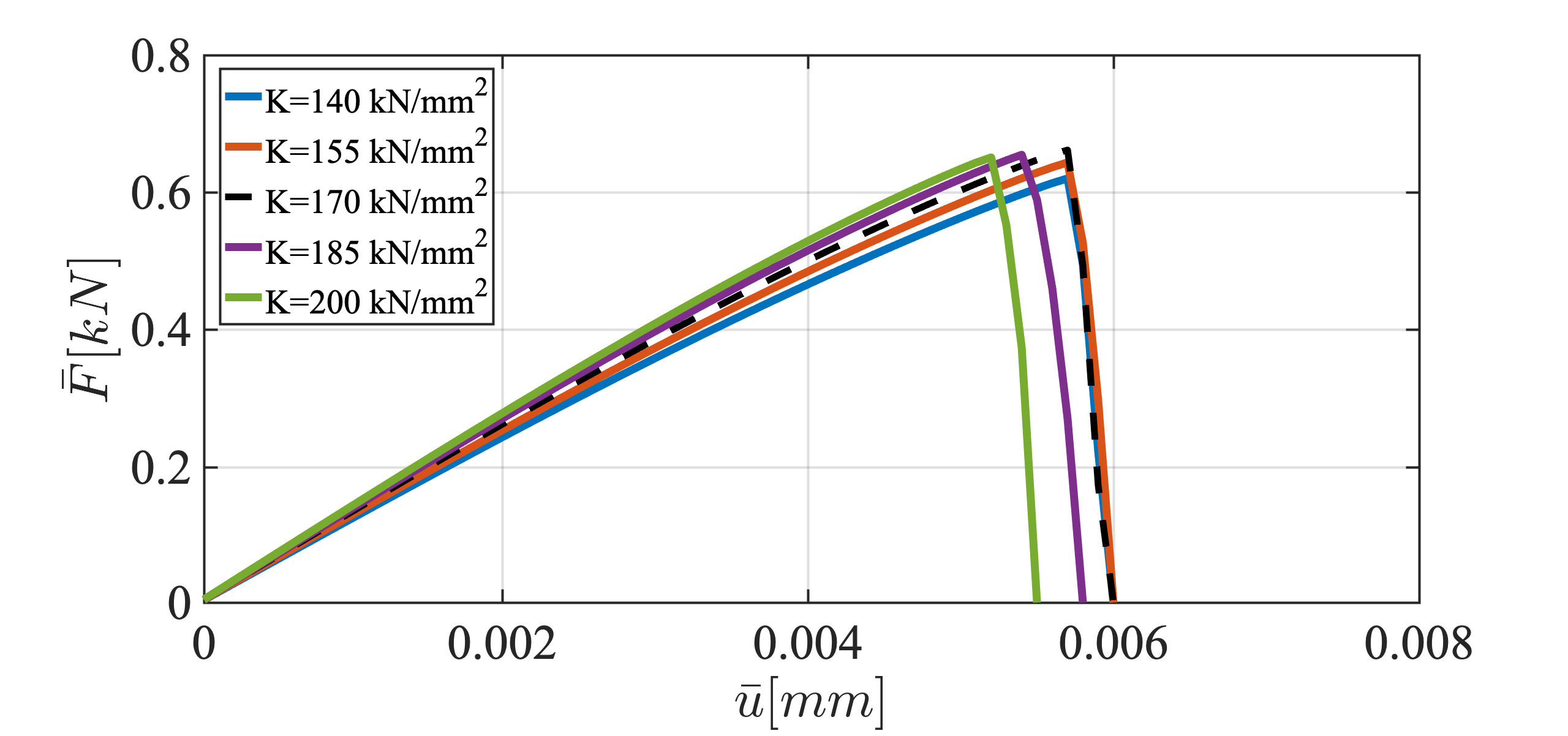}}%
	\newline
	\subfloat{\includegraphics[width=8cm,height=5cm]{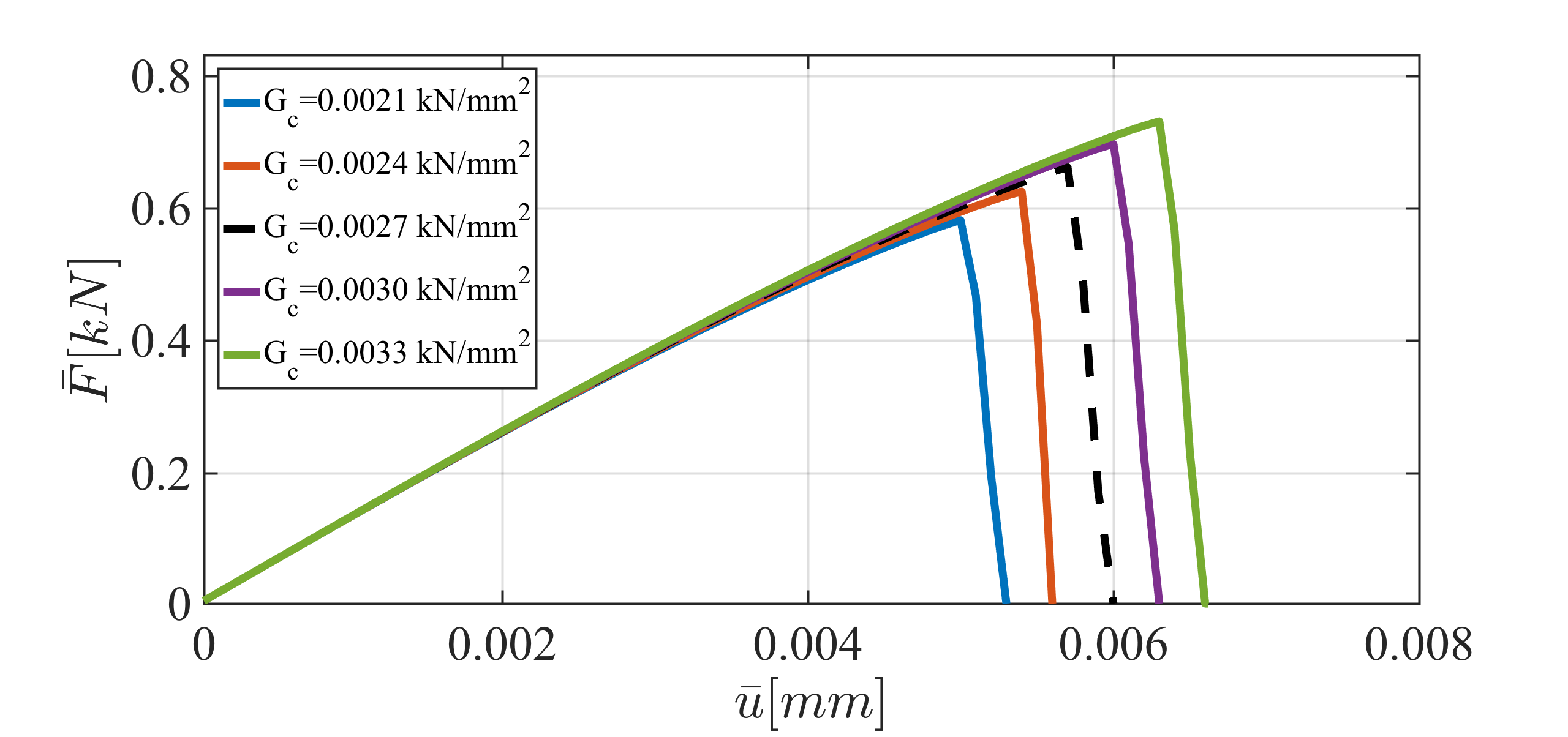}}%
	\caption{The load-displacement curve for different values of
		$\mu$  (top left), \new{$K$} (top right) and $G_c$ (bottom)
		in the SENT example (Example 1).}
	\label{fig:exam1_parameter}
\end{figure}

First, we describe the effect of each parameter on the
displacement. As the \new{elasticity} constants (i.e., $\mu$ and $\new{K}$)
become larger, the material response becomes stiffer; crack initiation
takes longer to occur. Additionally, a larger crack release energy
rate (as an indicator for the material resistance against the crack
driving force) delays crack nucleation and hence crack
dislocation. All these facts are illustrated in Figure
\ref{fig:exam1_parameter}.

\begin{figure}[t!]
	\subfloat{\includegraphics[width=7cm,height=4.5cm]{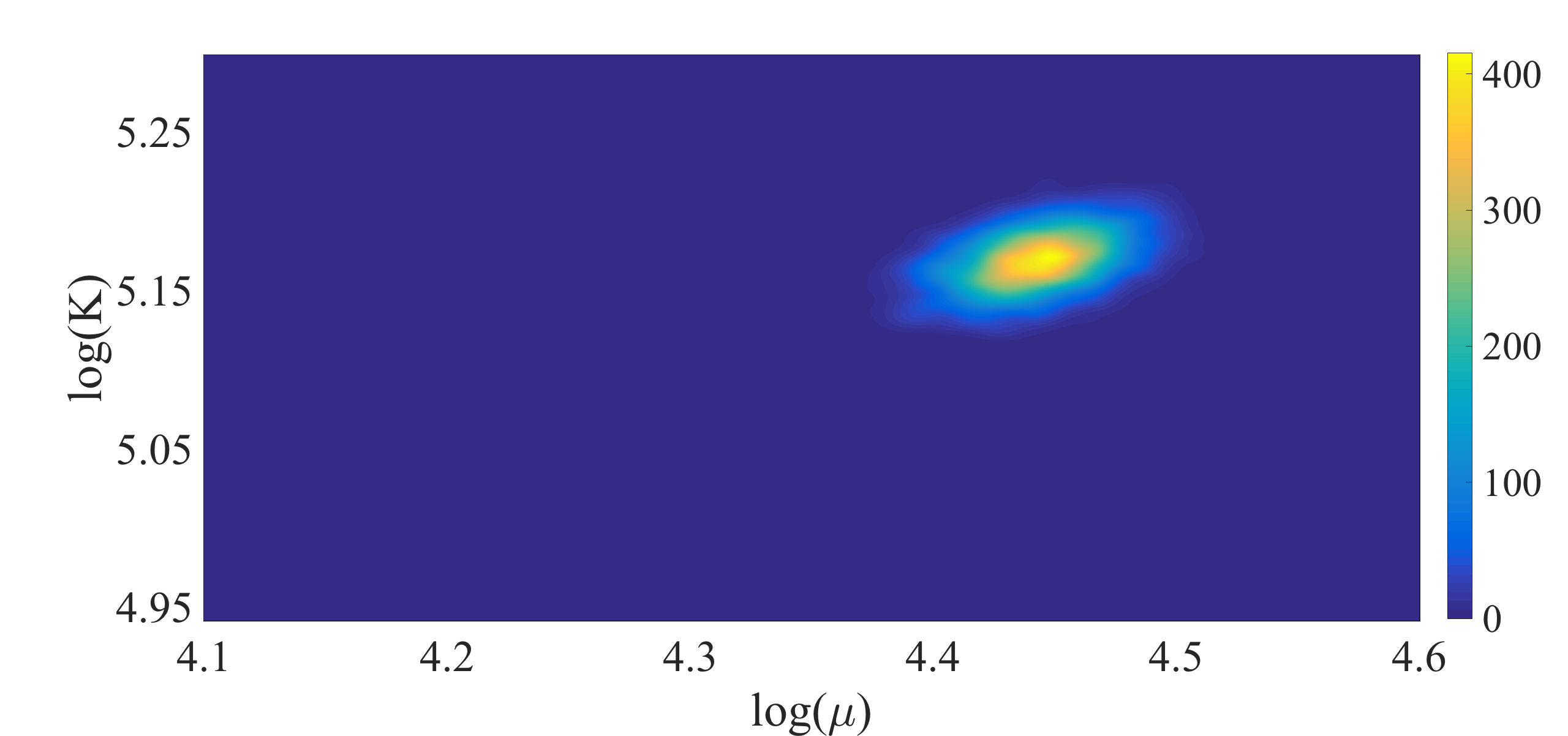}}%
	\subfloat{\includegraphics[width=7cm,height=4.5cm]{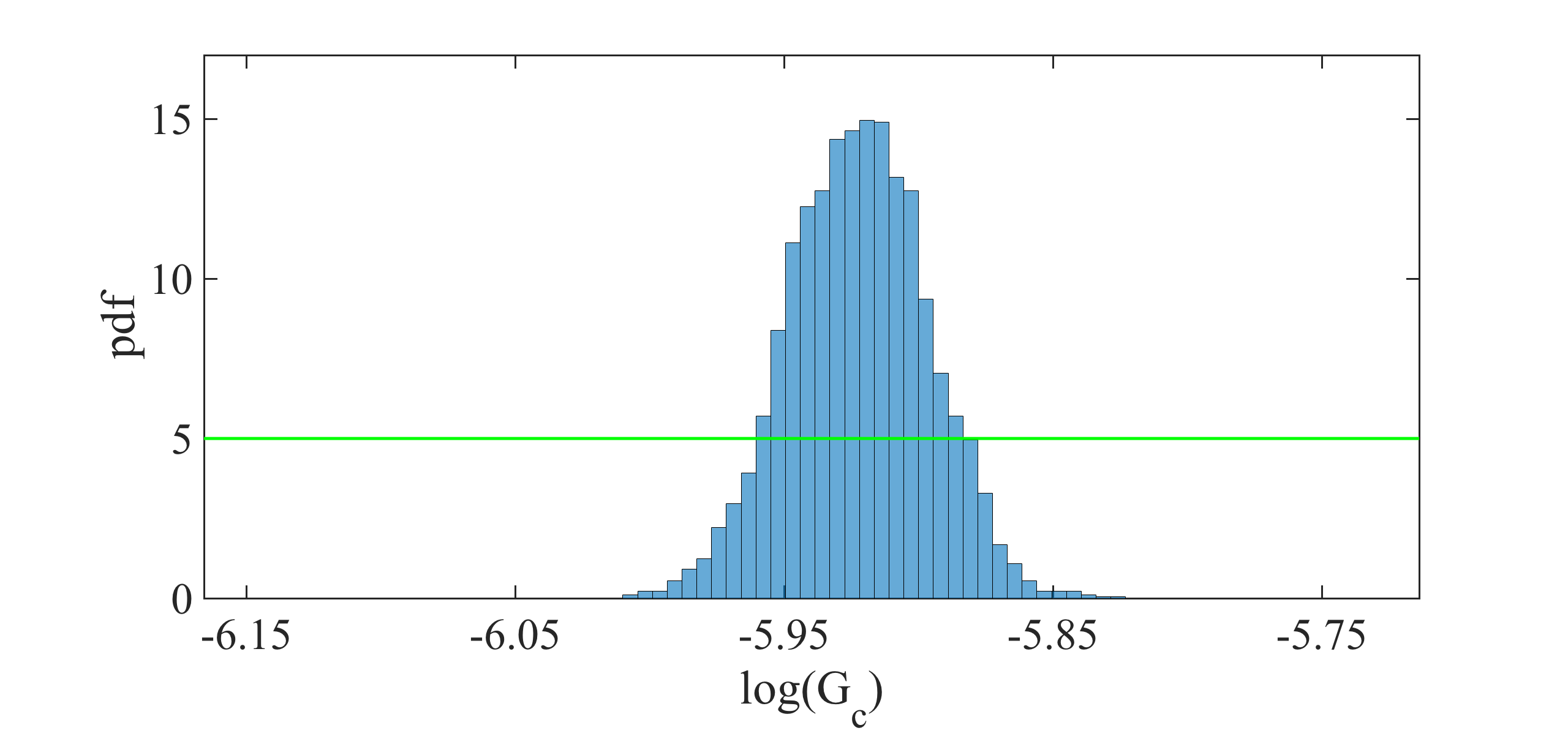}}%
	\hspace{-0.8cm}
	\newline
	\subfloat{\includegraphics[width=7cm,height=4.5cm]{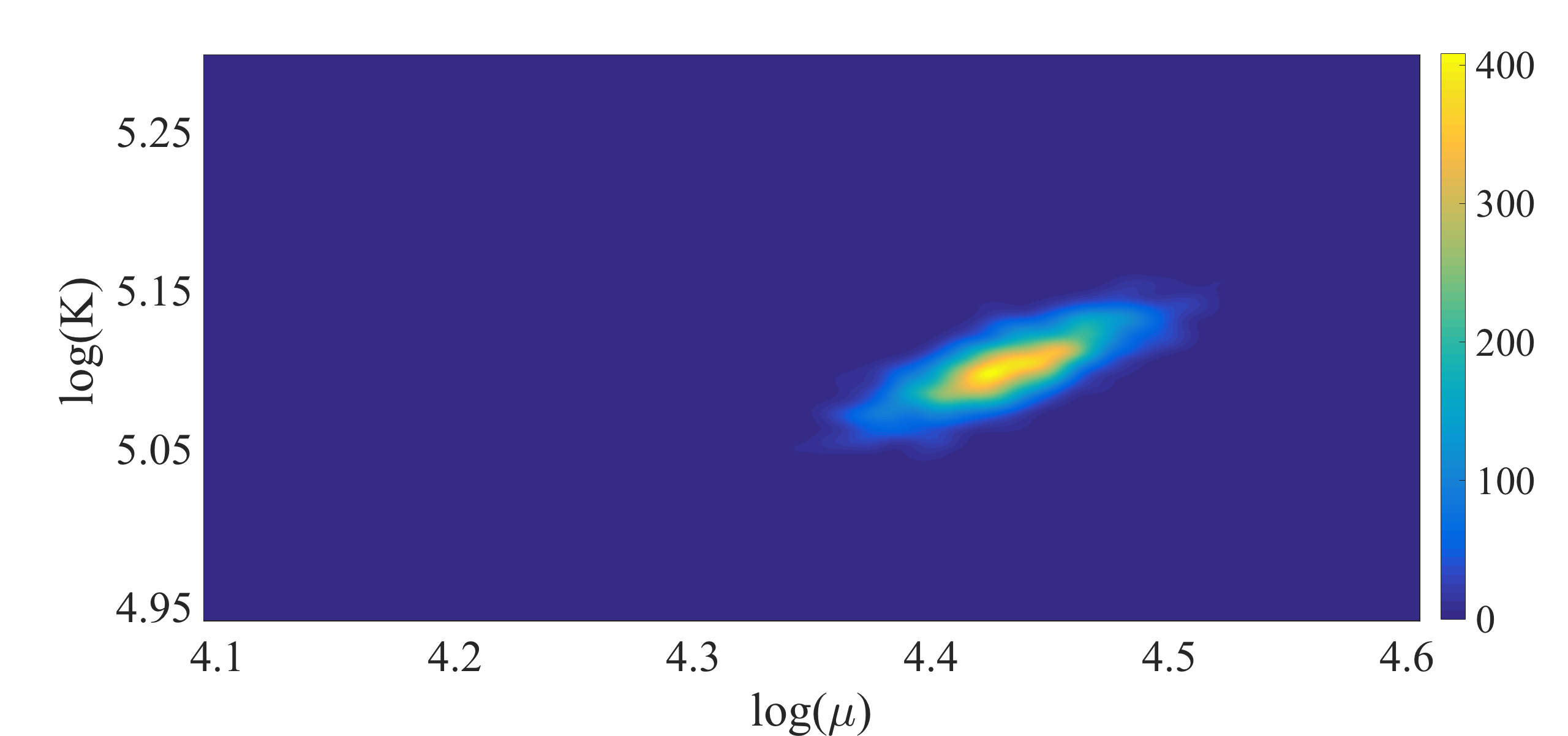}}%
	\subfloat{\includegraphics[width=7cm,height=4.5cm]{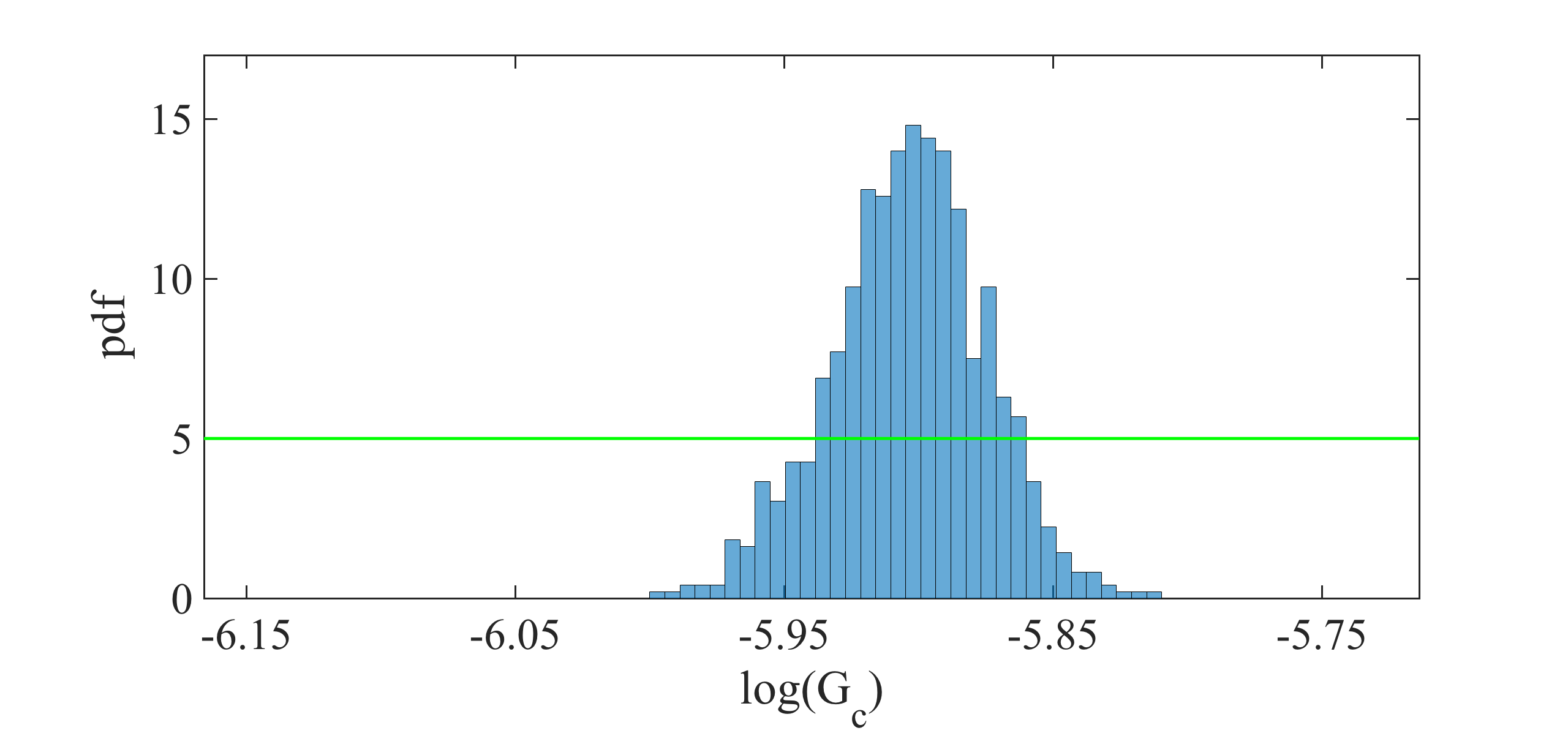}}%
	\caption{\new{Left: the joint probability density of the elasticity parameters. Right:} the  
		prior (green line),
		and the normalized probability density function (pdf) of $G_c$
		for the SENT example. Here we compare the distributions obtained by the one-dimensional Bayesian inversion (first row) and the three-dimensional Bayesian inversion (second row).}
	\label{fig:exam1_hist}
\end{figure}

\begin{table}
	\begin{tabular}{lllllll}
		\hline
		& \text{one-dimensional}&& \text{three-dimensional} &&  \\
		\cline{2-3}\cline{4-5} 
		& mean ($\mathrm{kN/mm^2}$) &   ratio (\%)	 & mean ($\mathrm{kN/mm^2}$)&   ratio (\%)     \\
		\hline
		$\mu$     &   \new{84.2}  & \new{27} &  \new{85.1}   &  \new{28}    \\
		\new{$K$}   &   \new{176.1} &  \new{27}  & 	\new{175.3}    & \new{28}     \\
		$G_c$       & \new{$0.00272$}     & \new{29.1}   & \new{$0.00268$}    & \new{28.6}   \\
		\hline
	\end{tabular}
	\caption{The mean values of the posterior distributions
		obtained by one-dimensional and three-dimensional Bayesian
		inversion in addition to their acceptance ratios for  the SENT (Example 1). The units are in $\mathrm{kN/mm^2}$.} 
	\label{Ex11D3D}
\end{table}

\begin{figure}[ht!]
	\centering
	\includegraphics[width=10cm,height=6.5cm]{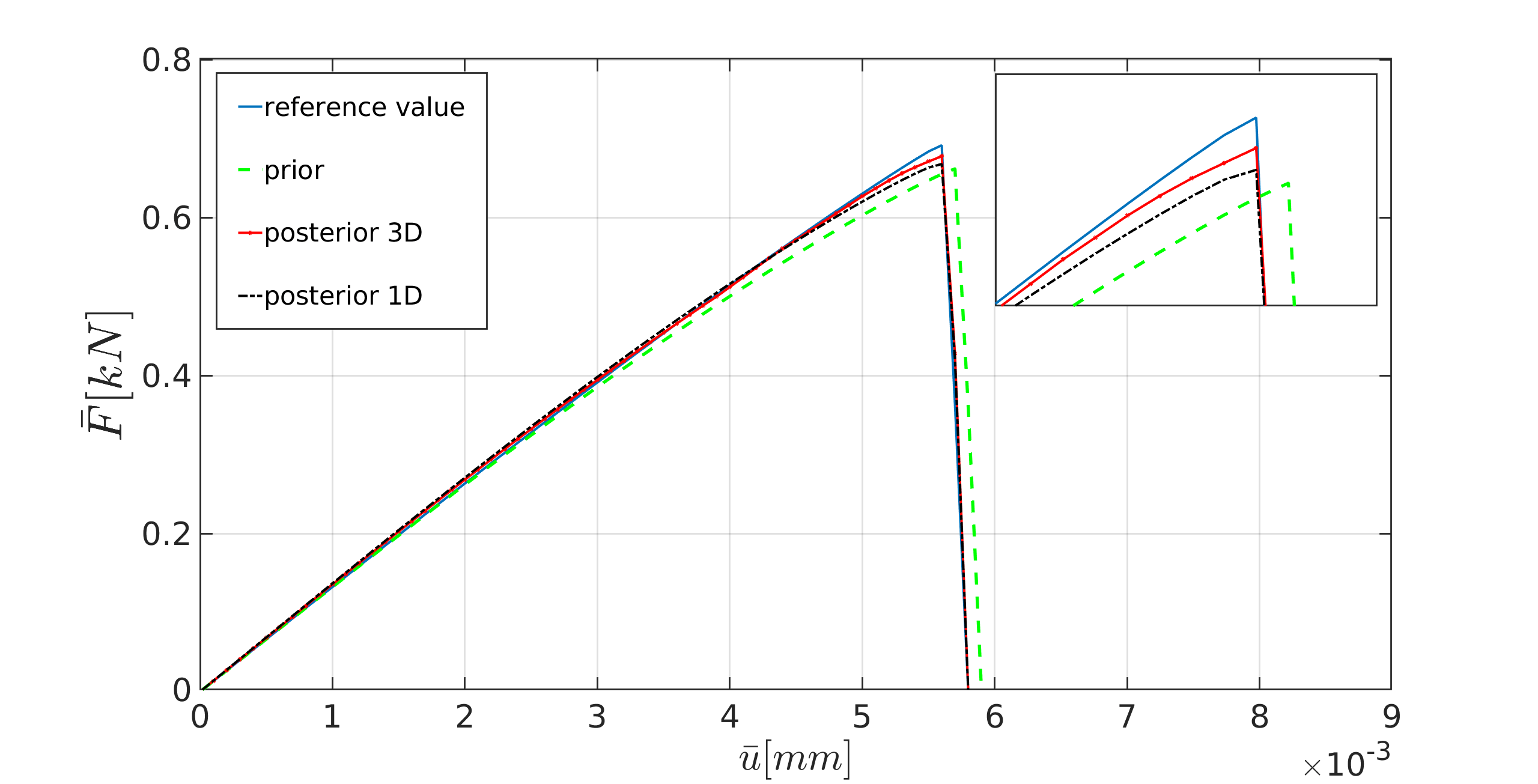}
	\caption{The load-displacement curve for the one-dimensional
		(black) and three-dimensional (red) posterior distributions
		in addition to the ones for the  prior  distribution (green)
		and the reference value (blue) for the SENT example (Example 1) with $h=1/160$.}
	\label{fig:Ex1_posterior}
\end{figure}

\new{The joint probability density of the elasticity parameters and the marginal probability of the posterior are shown in Figure \ref{fig:exam1_hist} including  one- and
	three-dimensional Bayesian inversions. 
	The mean values of the distributions are
	\new{$\mu=\unit{84.2}{kN/mm^2}$} and $\mu=\unit{85.1}{kN/mm^2}$ are obtained for the shear modulus.  
	Here an acceptance rate of 27\% is obtained.  Regarding the
	material stiffness parameter $G_c$, the acceptance rates are near
	29\%. The values are summarized in Table \ref{Ex11D3D}.}

To verify the parameters obtained by the Bayesian approach, we solved
the forward model using the mean values of the posterior
distributions. Figure \ref{fig:Ex1_posterior} shows the
load-displacement diagram according to prior and posterior
distributions. As expected, during the nucleation and propagation
process, using Bayesian inversion results in better agreement
(compared to the prior). Furthermore, a better estimation is achieved
by simultaneous \emph{multi-dimensional} Bayesian inversion. From now
onward, this approach will be used for Bayesian inference.

\begin{figure}[ht!]
	\centering
	\includegraphics[width=10cm,height=6.0cm]{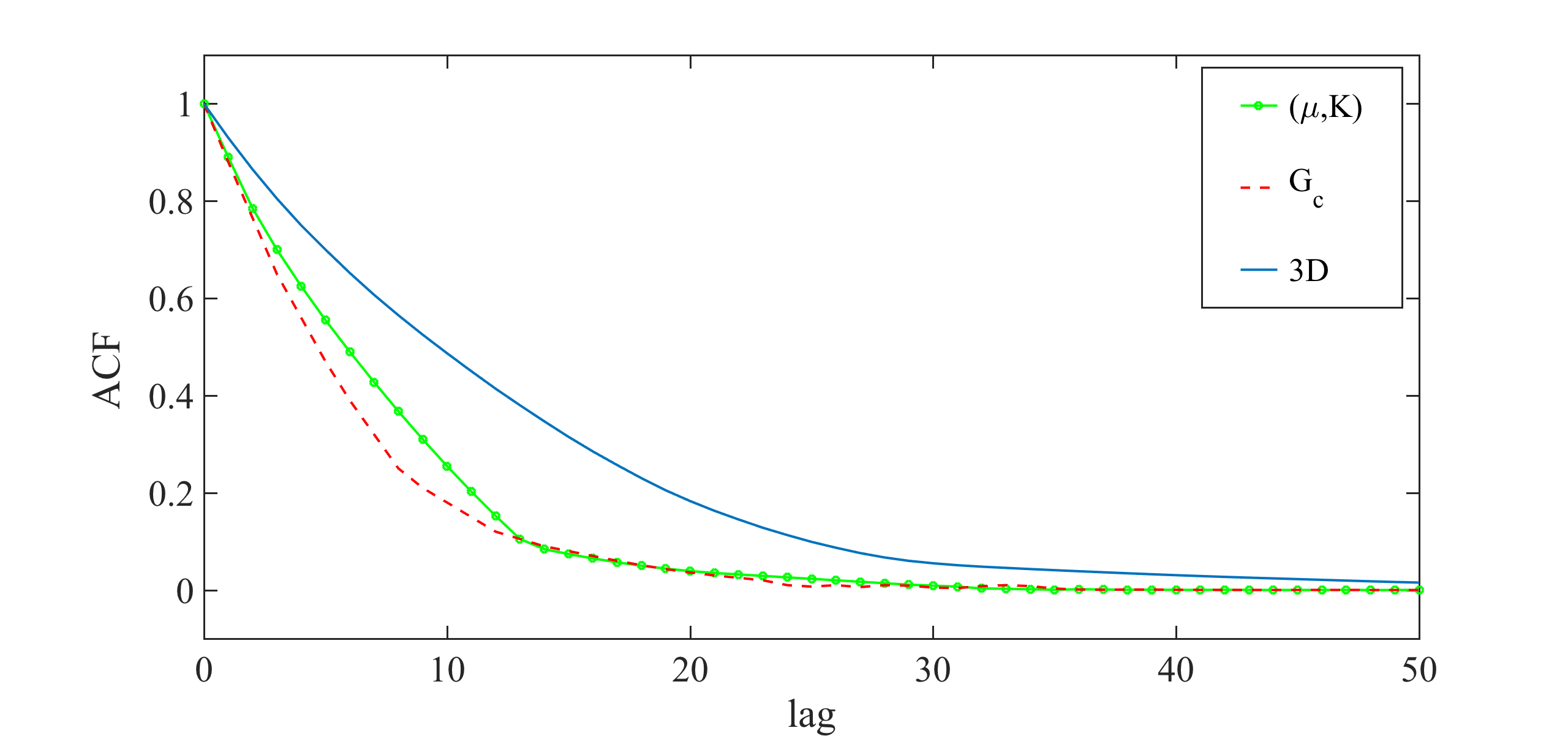}
	\caption{The autocorrelation function for one- and
		multidimensional Bayesian inference in the SENT example.}
	\label{Ex3:ACF}
\end{figure}

\subsubsection{The convergence of MCMC}
A customary method to assess the convergence of the MCMC is the
calculation of its autocorrelation. The lag-$\tau$ autocorrelation
function (ACF) $R\colon\mathbb{N}\rightarrow [-1,1]$ is defined as
\begin{equation*}
R(\tau):=\frac{\sum_{n=1}^{N-\tau} (\theta_n -\bar{\theta})(\theta_{n+\tau} -\bar{\theta}) }{ \sum_{n=1}^{N}\left( \theta_n-\bar{\theta}\right)^2}=\frac{\text{cov}(\theta_n,\theta_{n+\tau})}{\text{var}(\theta_n)},
\end{equation*}
where $\theta_n$ is the $n$-th element of the Markov chain and
$\bar{\theta}$ indicates the mean value. For the Markov chains,
$R(\tau)$ is positive and strictly decreasing. Also, a rapid decay in
the ACF indicates the samples are not fully  correlated and mixing well. Figure
\ref{Ex3:ACF} shows the convergence of the MCMC where the \new{elasticity parameters} and the crirical elastic energy.  Also, we estimated the convergence
observed in the multi-dimensional approach.    As expected,
the multi-dimensional approach converges slower than the
one-dimensional one.

\begin{figure}[t!]
	\centering
	\subfloat{\includegraphics[width=4.5cm,height=4.5cm]{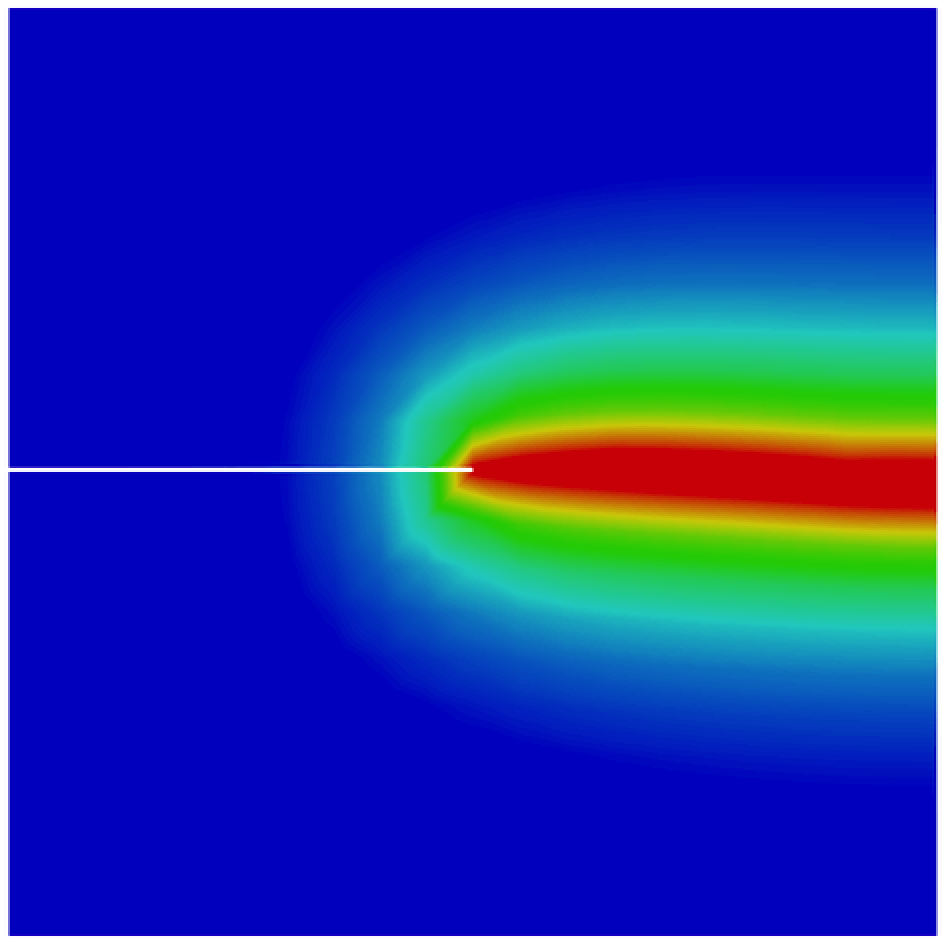}}%
	\hspace{.7cm} 
	\subfloat{\includegraphics[width=4.5cm,height=4.5cm]{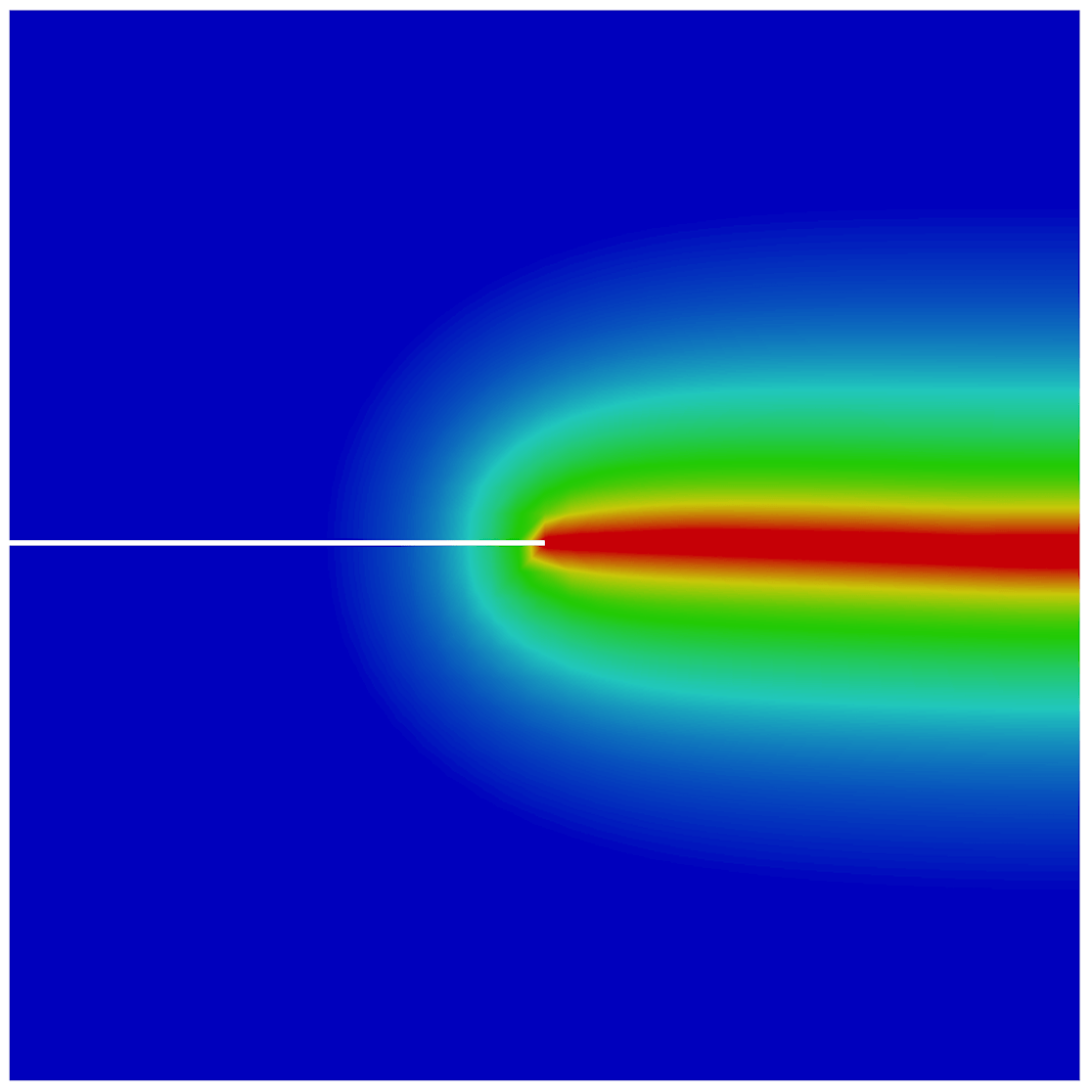}}
	\hspace{.7cm}
	\subfloat{\includegraphics[width=4.5cm,height=4.5cm]{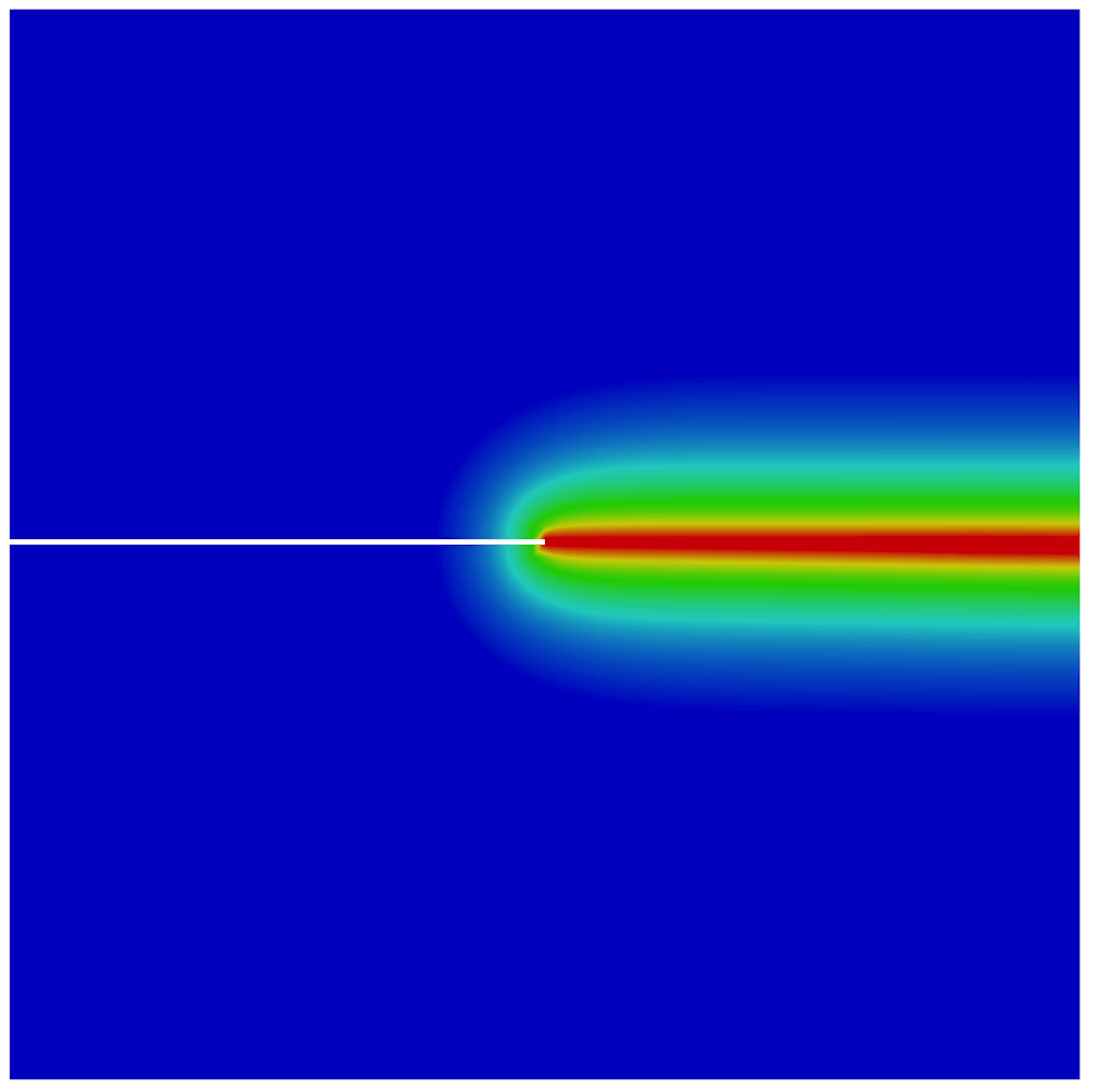}}\newline\\		\subfloat{\includegraphics[width=4.5cm,height=4.5cm]{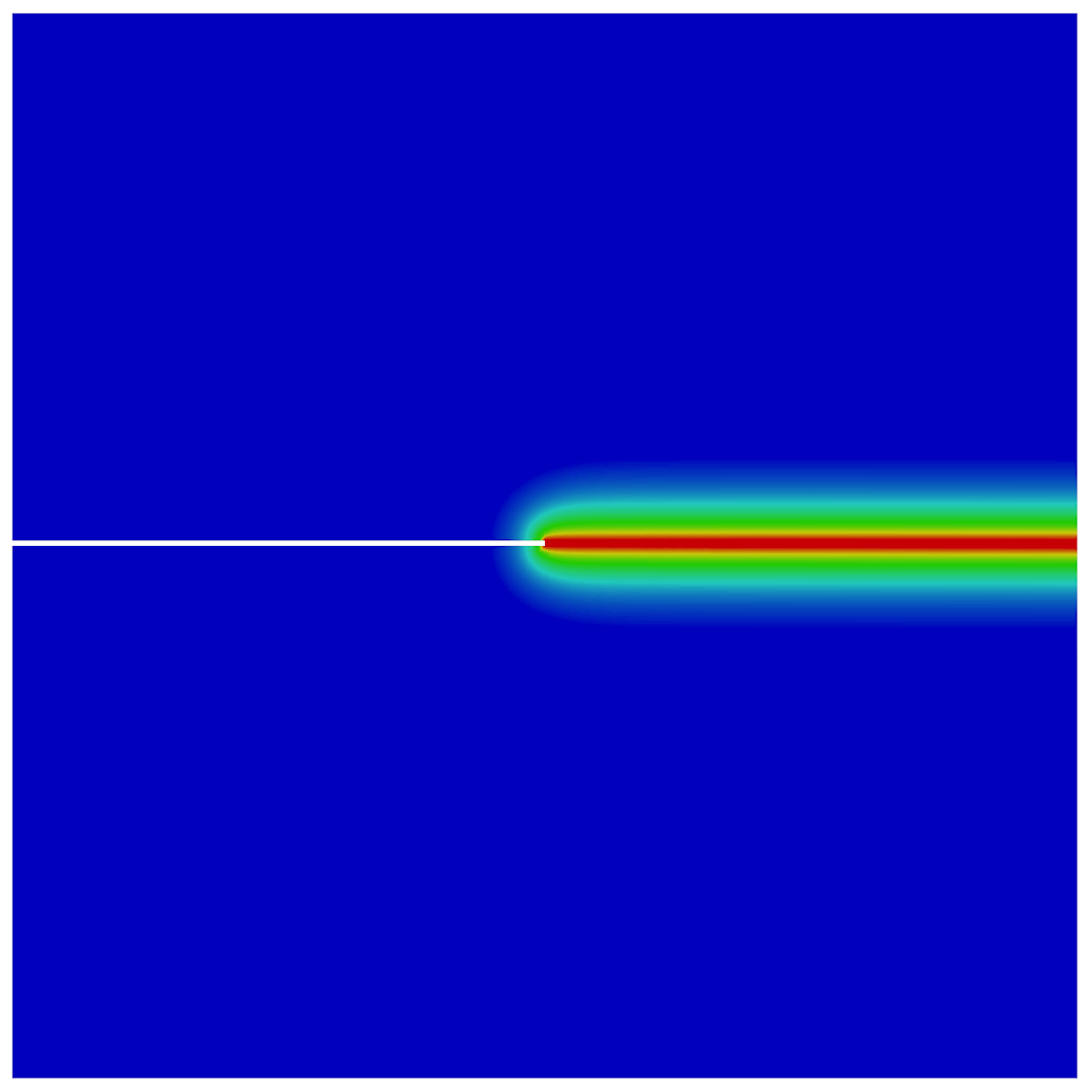}}
	\hspace{.7cm}
	\subfloat{\includegraphics[width=4.5cm,height=4.5cm]{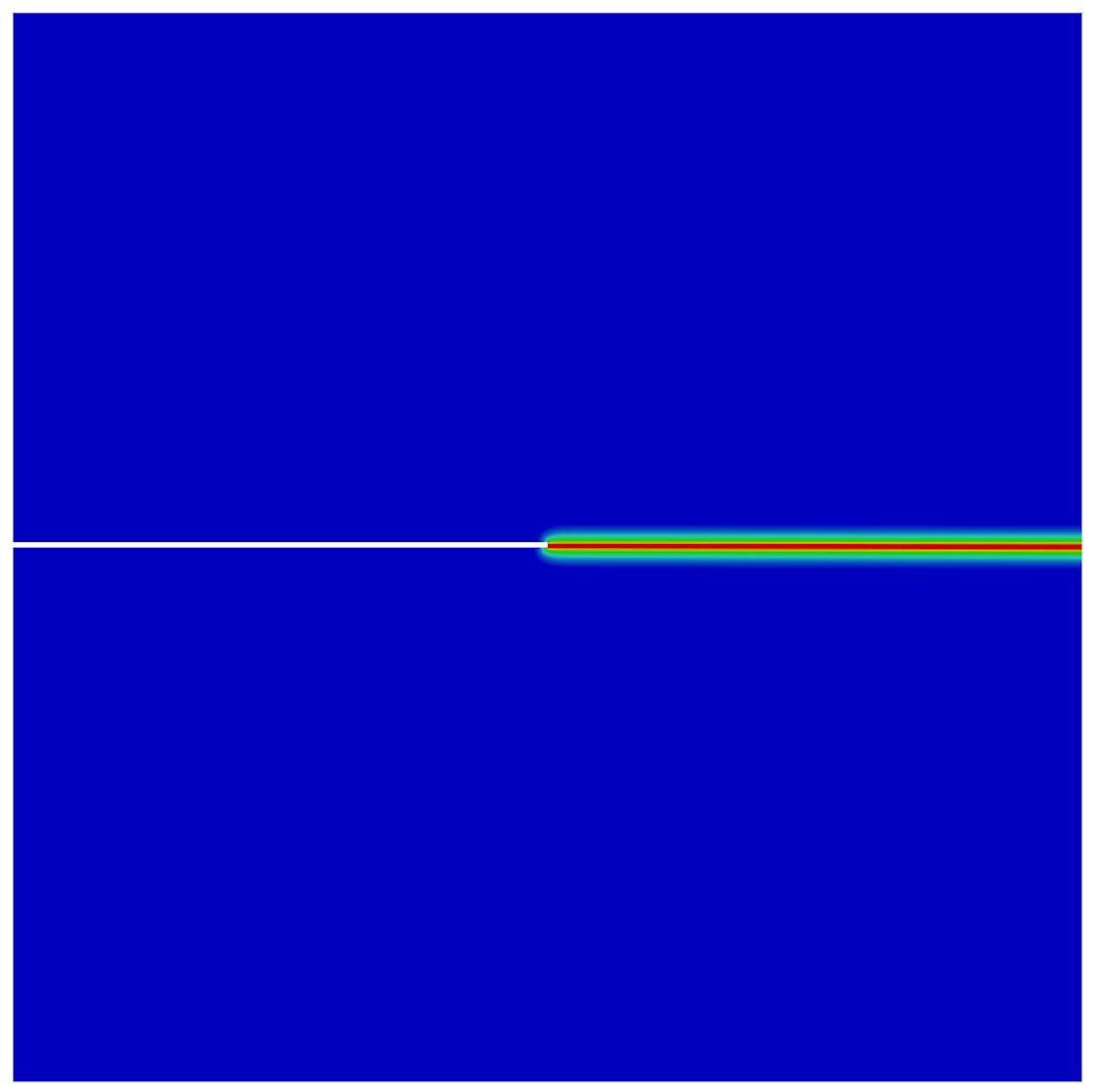}}
	\hspace{.5cm} 
	\subfloat{\includegraphics[width=2.0cm,height=4.0cm]{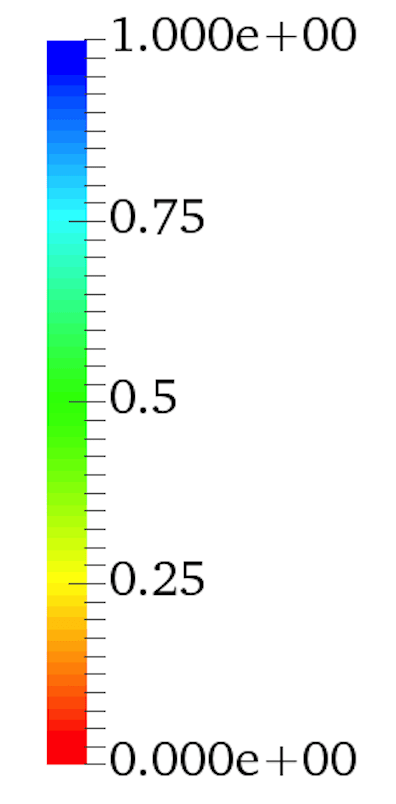}}
	\caption{The effect of the mesh size on the crack propagation
		in the SENT example (Example 1). The mesh sizes are (from the left) $h=1/20$, $h=1/40$, $h=1/80$, $h=1/160$, and $h=1/320$ (the reference). The effective parameters are chosen according to the prior values.}
	\label{fig:exam1_mesh1}
\end{figure}

\begin{figure}[ht!]
	\centering
	\includegraphics[width=10cm,height=6.5cm]{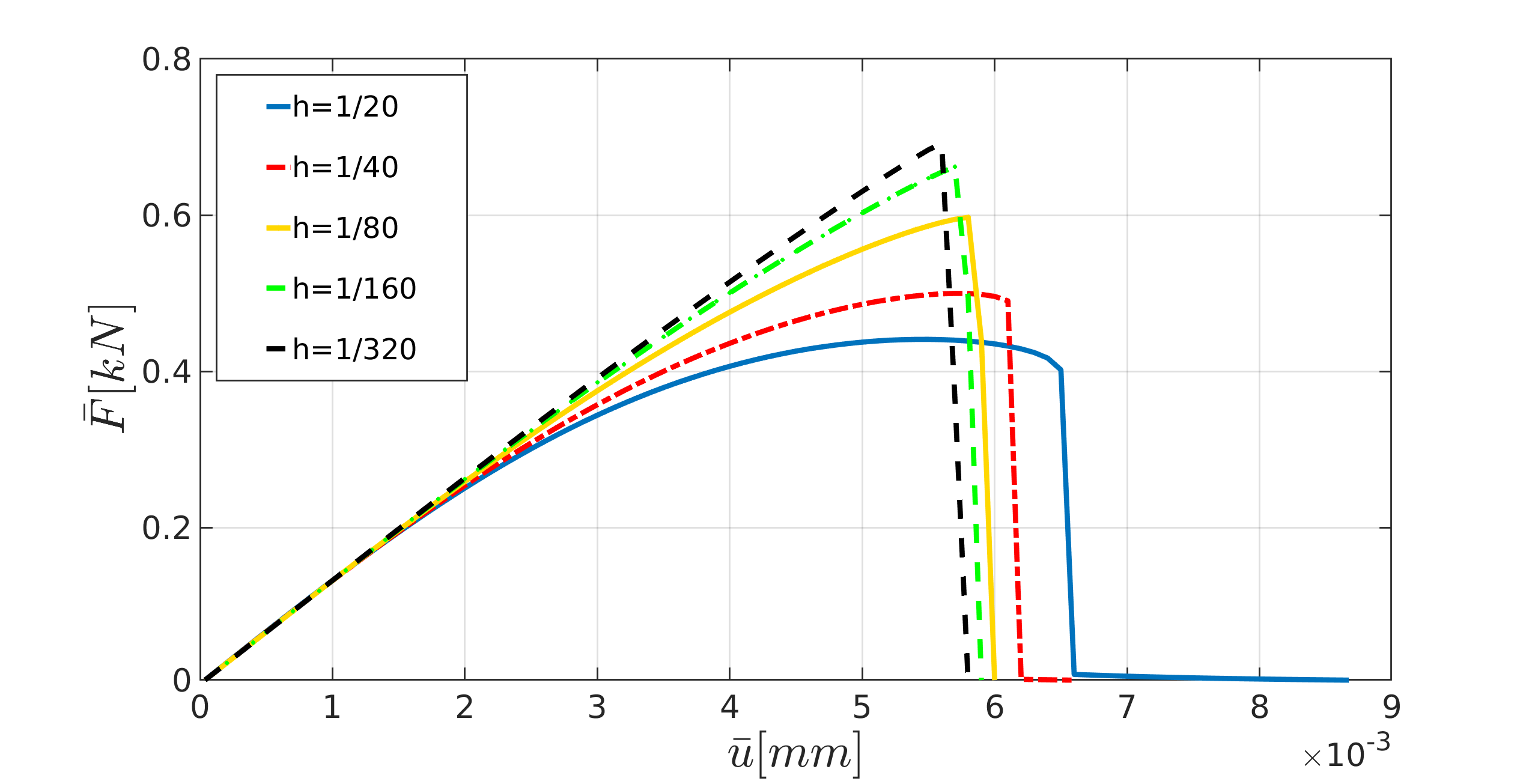}
	\caption{The load-displacement  curve in the SENT example
		(Example 1) for different mesh sizes, where the parameters
		are chosen according to the prior.}
	\label{fig:Ex1prior}
\end{figure}


As noted above, the phase-field solution depends on $h$ and $\ell$. A detailed computational 
analysis was for instance performed in \cite{HeWheWi15,HeiWi18_pamm}. 
In general, for smaller $h$ (and also smaller $\ell$) the crack path is better resolved, but 
leads to a much higher computational cost.

Figure \ref{fig:exam1_mesh1} illustrates the crack pattern using
different mesh sizes varying between $h=1/20$ and $h=1/320$. For these
mesh sizes, we show the load-displacement diagram in Figure
\ref{fig:Ex1prior} and the corresponding CPU time in Table \ref{time}.

\begin{figure}[!htb]
	\begin{tabular}{cc}
		\begin{minipage}{0.5\textwidth} 
			\includegraphics[width=10cm,height=6.5cm]{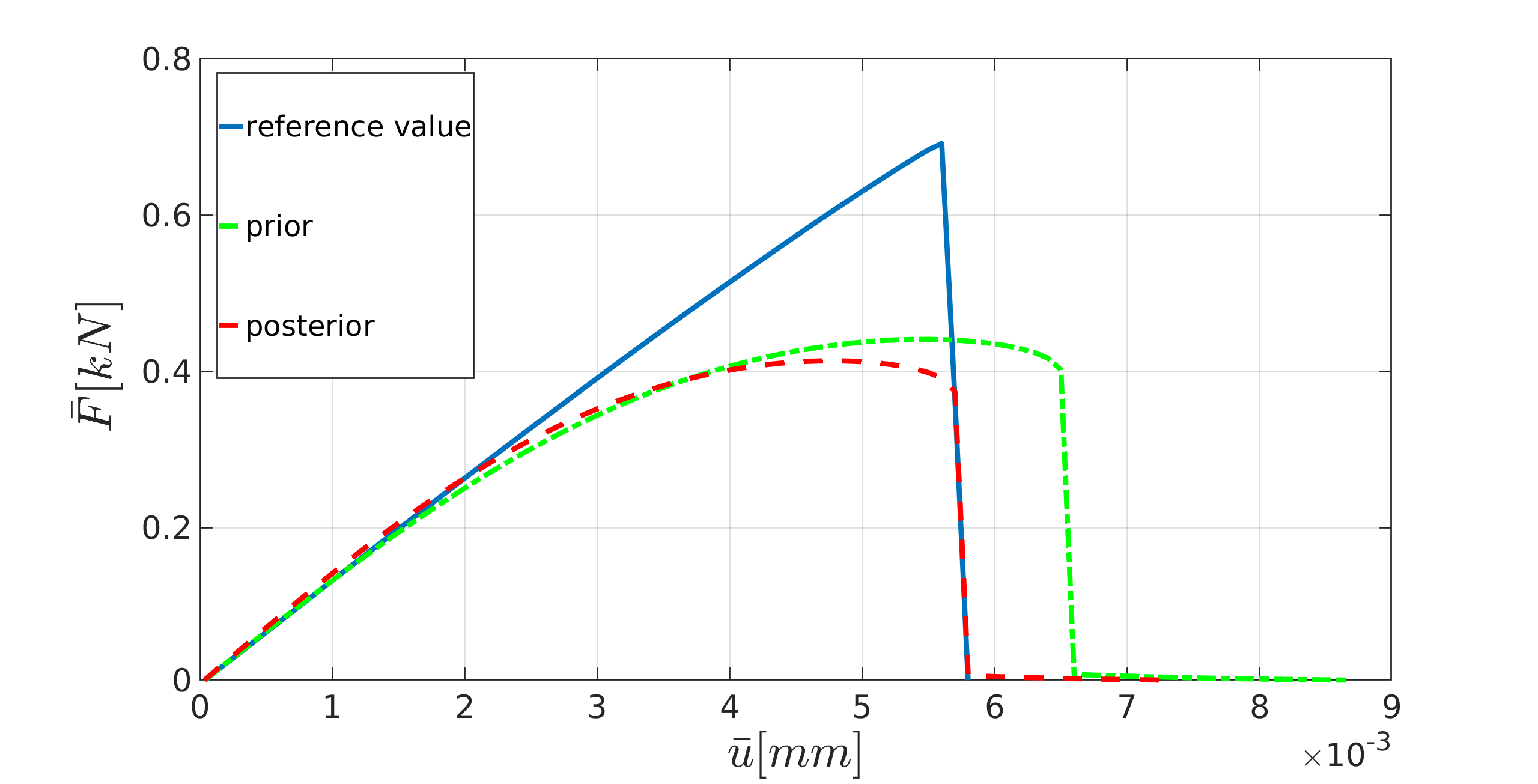} \end{minipage}& \begin{minipage}{0.1\textwidth} \includegraphics[width=5.5cm,height=3cm]{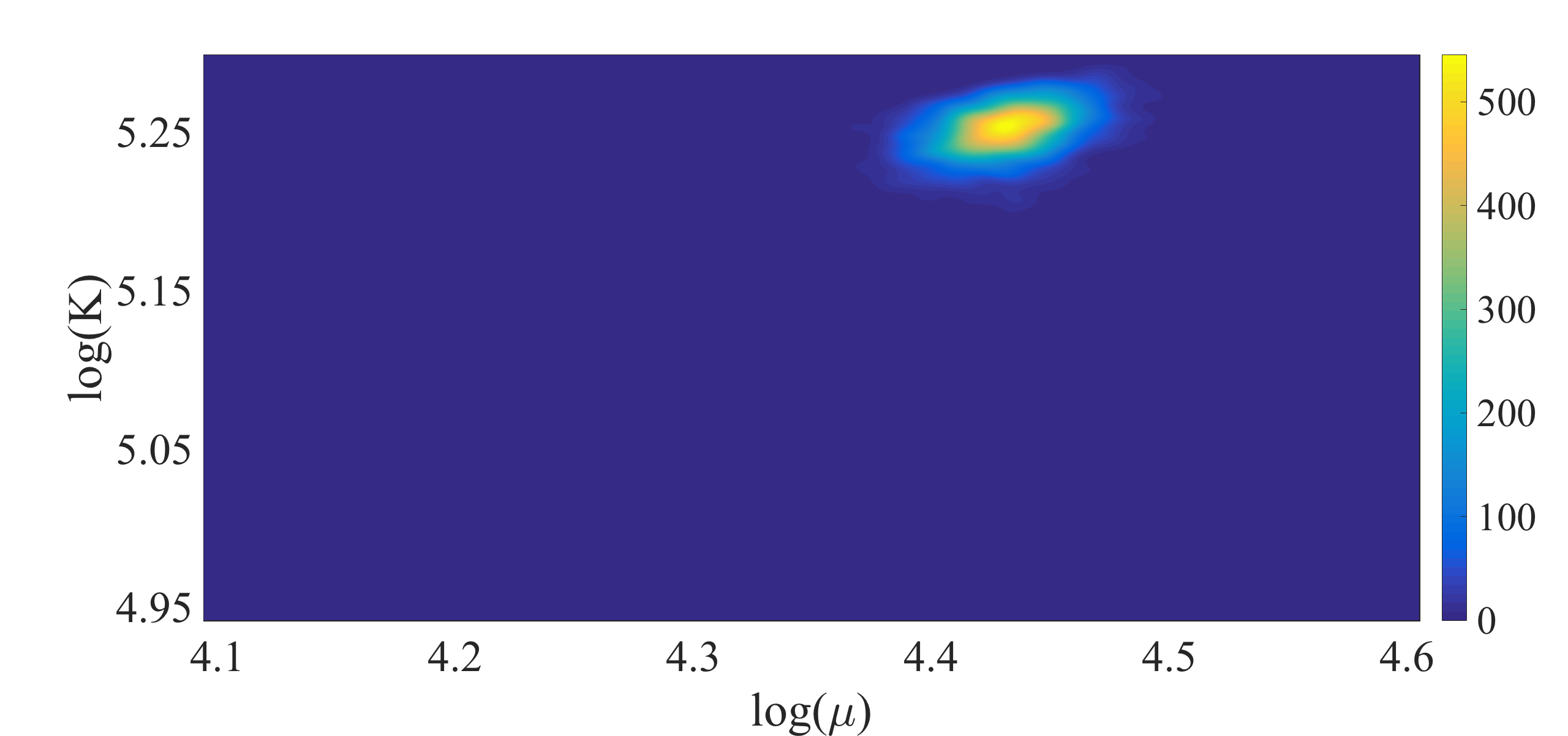} \\ \includegraphics[width=5.5cm,height=3cm]{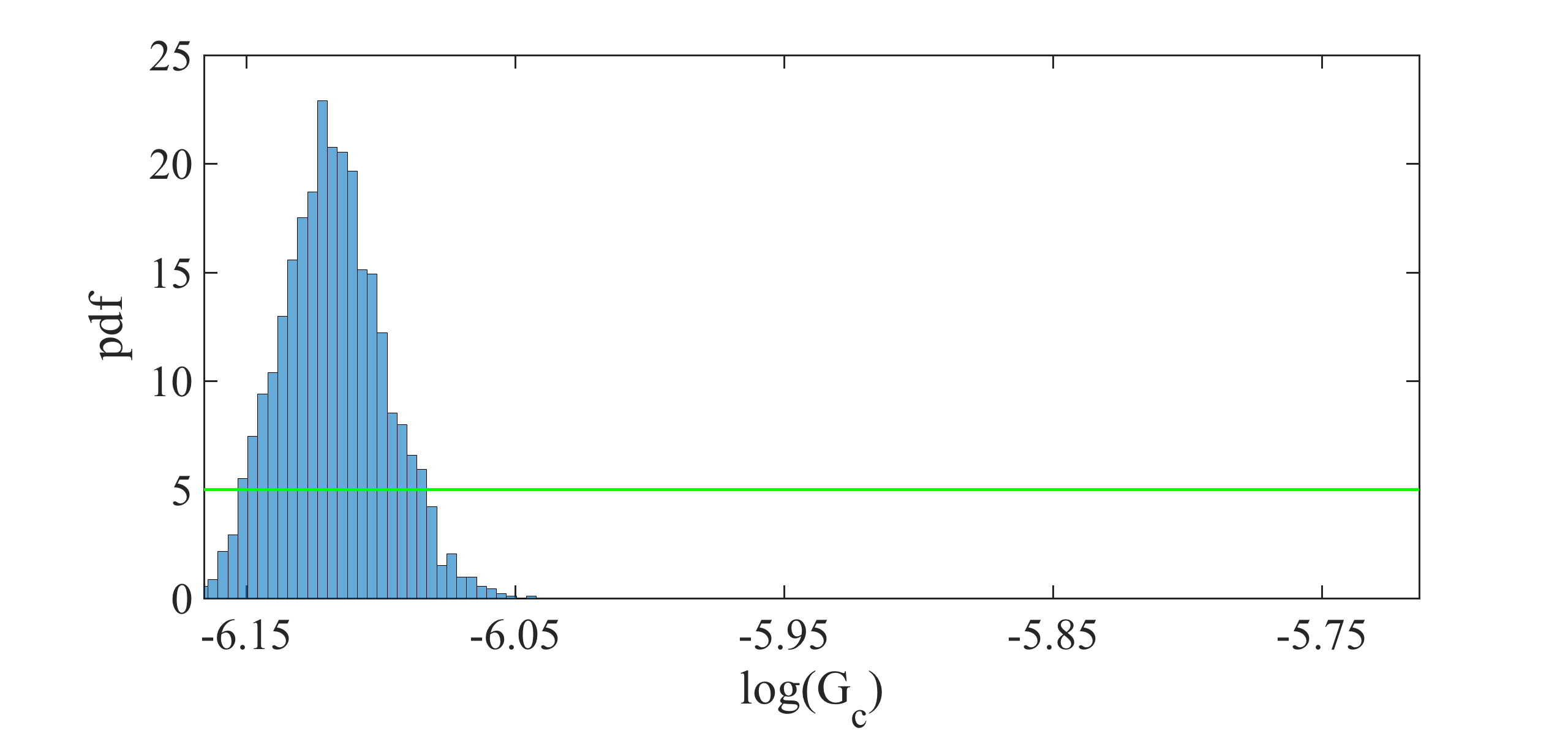} 
		\end{minipage}
	\end{tabular}
	\caption{The load-displacement curves for the reference, prior, and
		posterior values (left panel) for $h=1/20$ in the SENT
		example (Example 1). The \new{joint and} marginal posterior distributions of the effective parameters are shown in the right panel.}
	\label{fig:Ex1_h120}
\end{figure}

\begin{figure}[!htb]
	\begin{tabular}{cc}
		\begin{minipage}{0.5\textwidth} 
			\includegraphics[width=10cm,height=6.5cm]{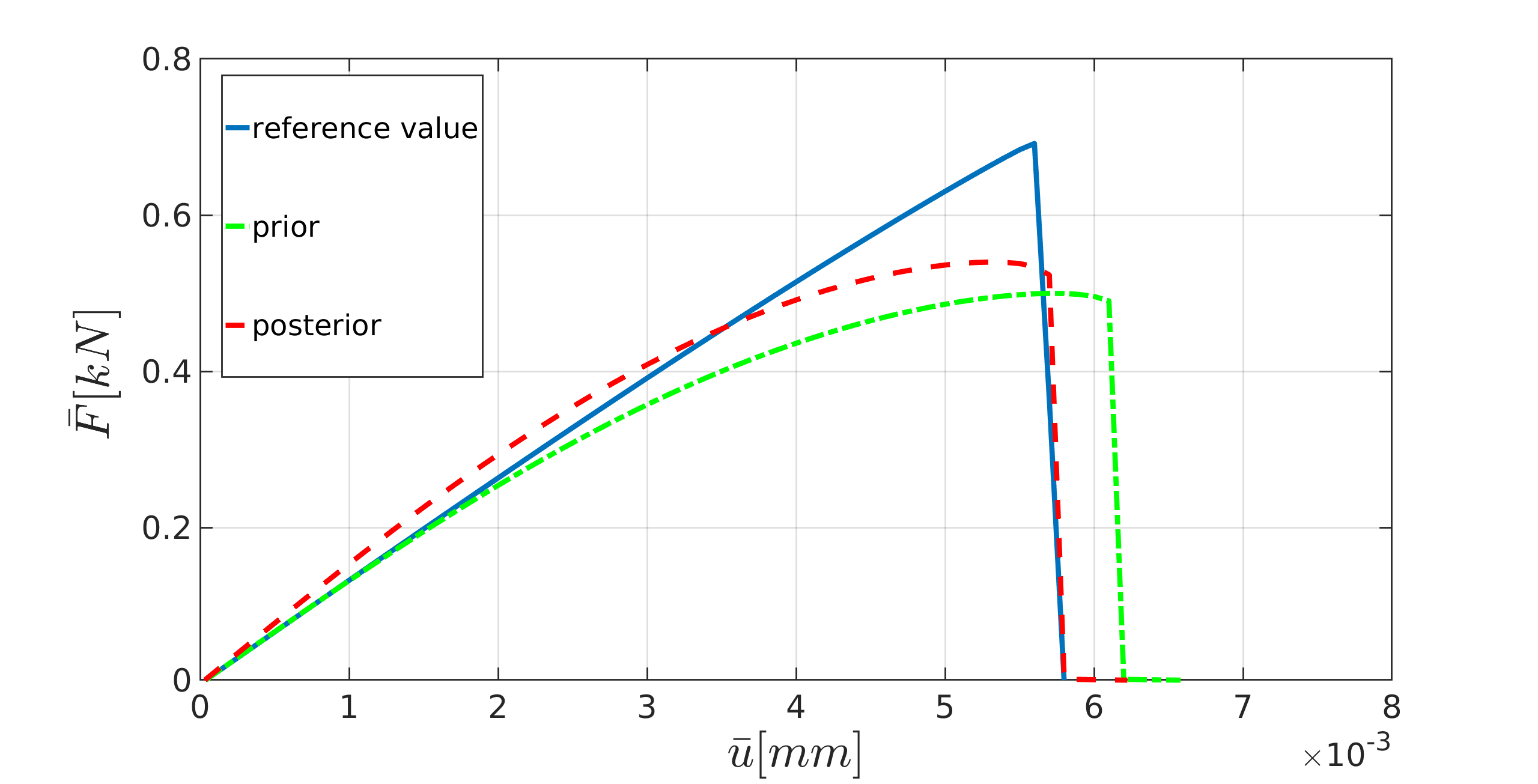} \end{minipage}& \begin{minipage}{0.1\textwidth} \includegraphics[width=5.5cm,height=3cm]{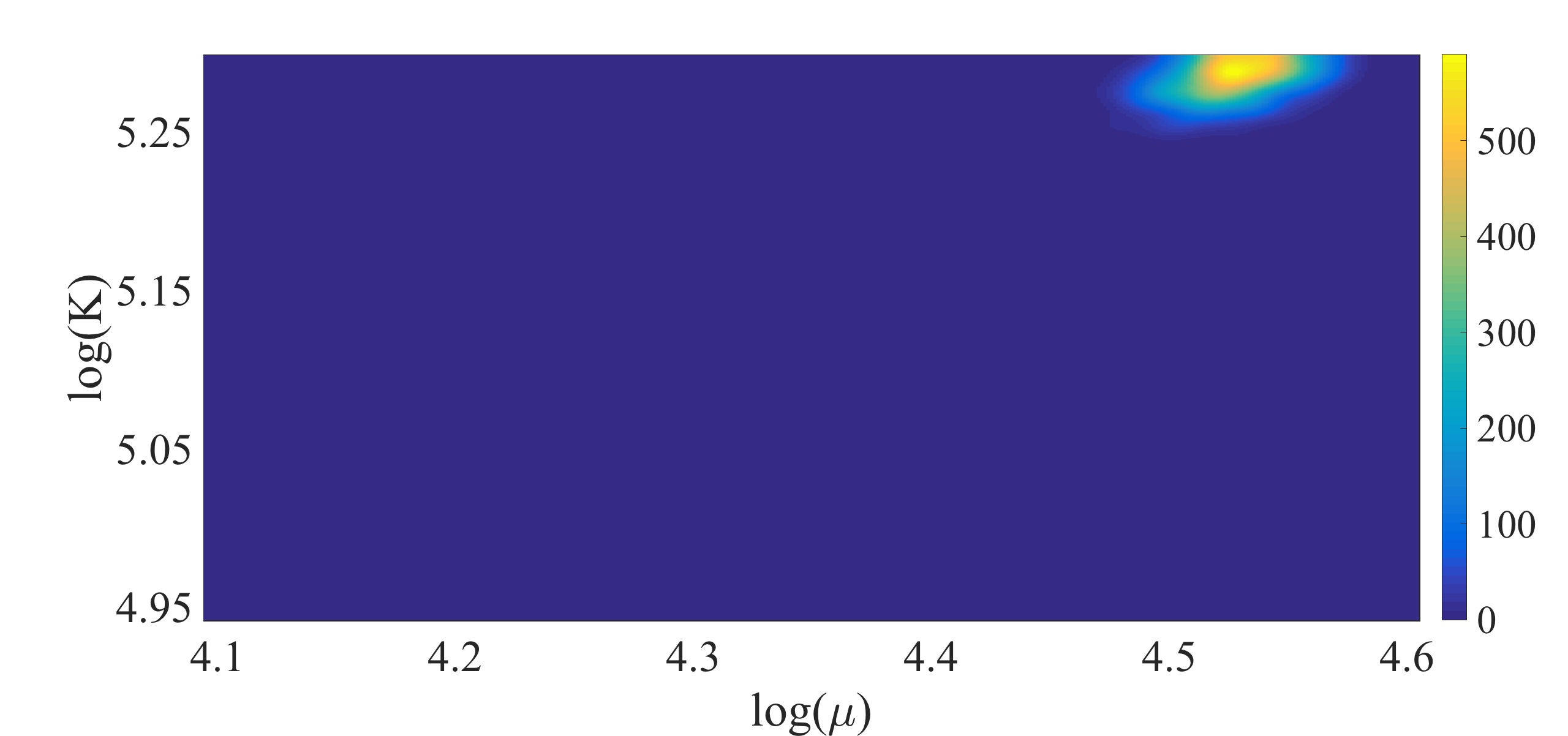} \\ \includegraphics[width=5.5cm,height=3cm]{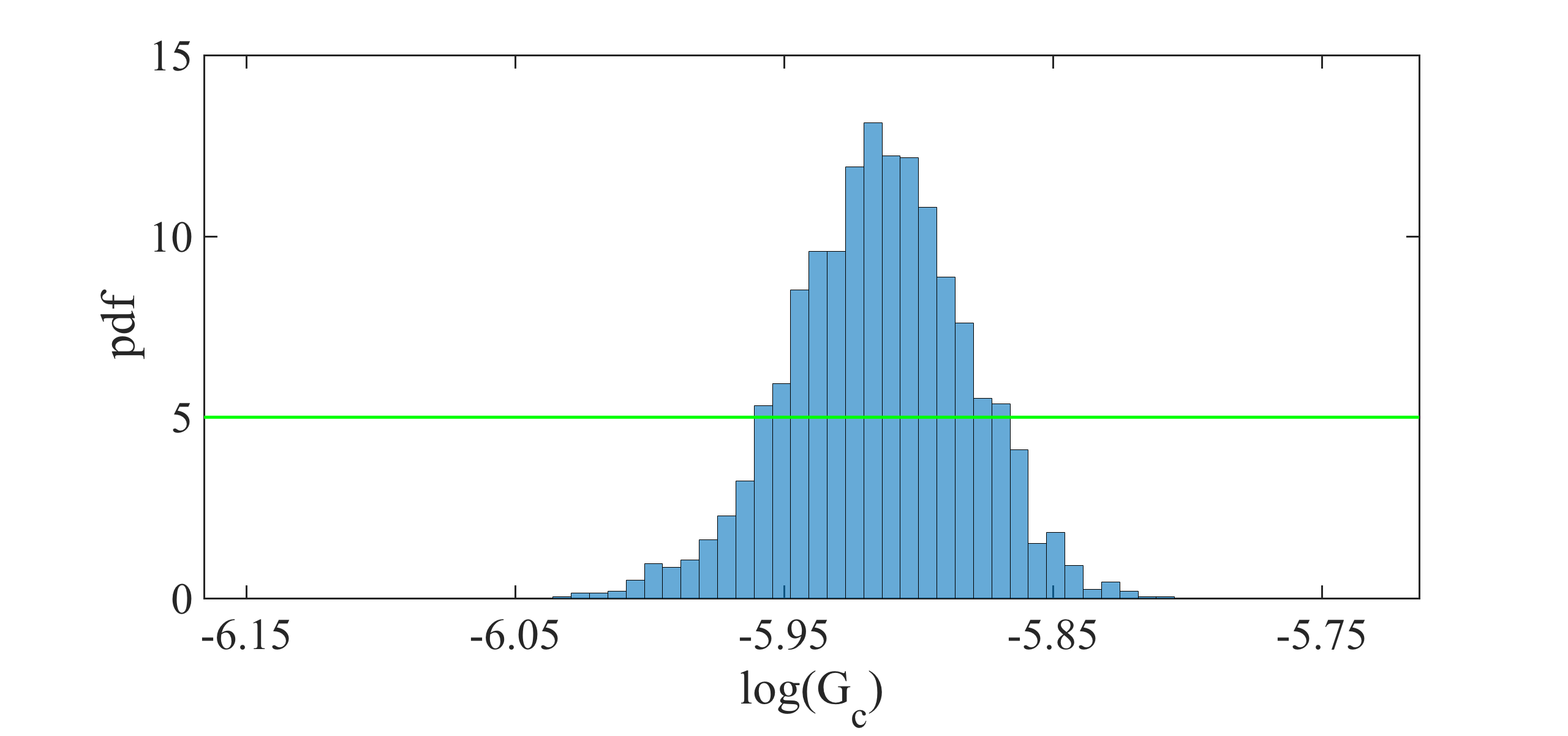}
		\end{minipage}
	\end{tabular}
	\caption{The load-displacement curves for the reference,
		prior, and posterior values (left panel) for $h=1/40$ in the
		SENT example (Example 1). The \new{joint and the} marginal posterior distributions of the effective parameters are shown in the right panel.}
	\label{fig:Ex1_h140}
\end{figure}

\begin{figure}[!htb]
	\begin{tabular}{cc}
		\begin{minipage}{0.5\textwidth} 
			\includegraphics[width=10cm,height=6.5cm]{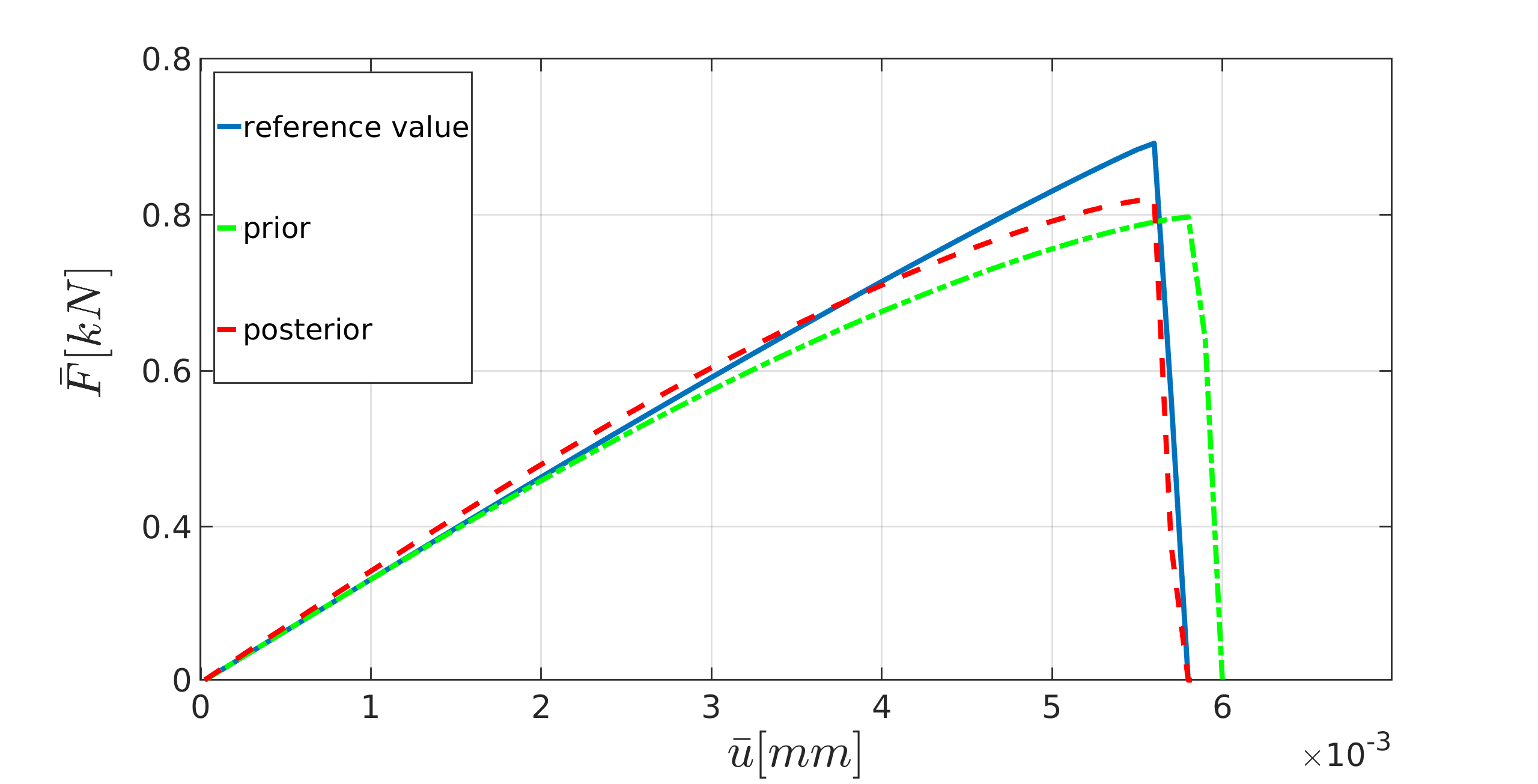} \end{minipage}& \begin{minipage}{0.1\textwidth} \includegraphics[width=5.5cm,height=3cm]{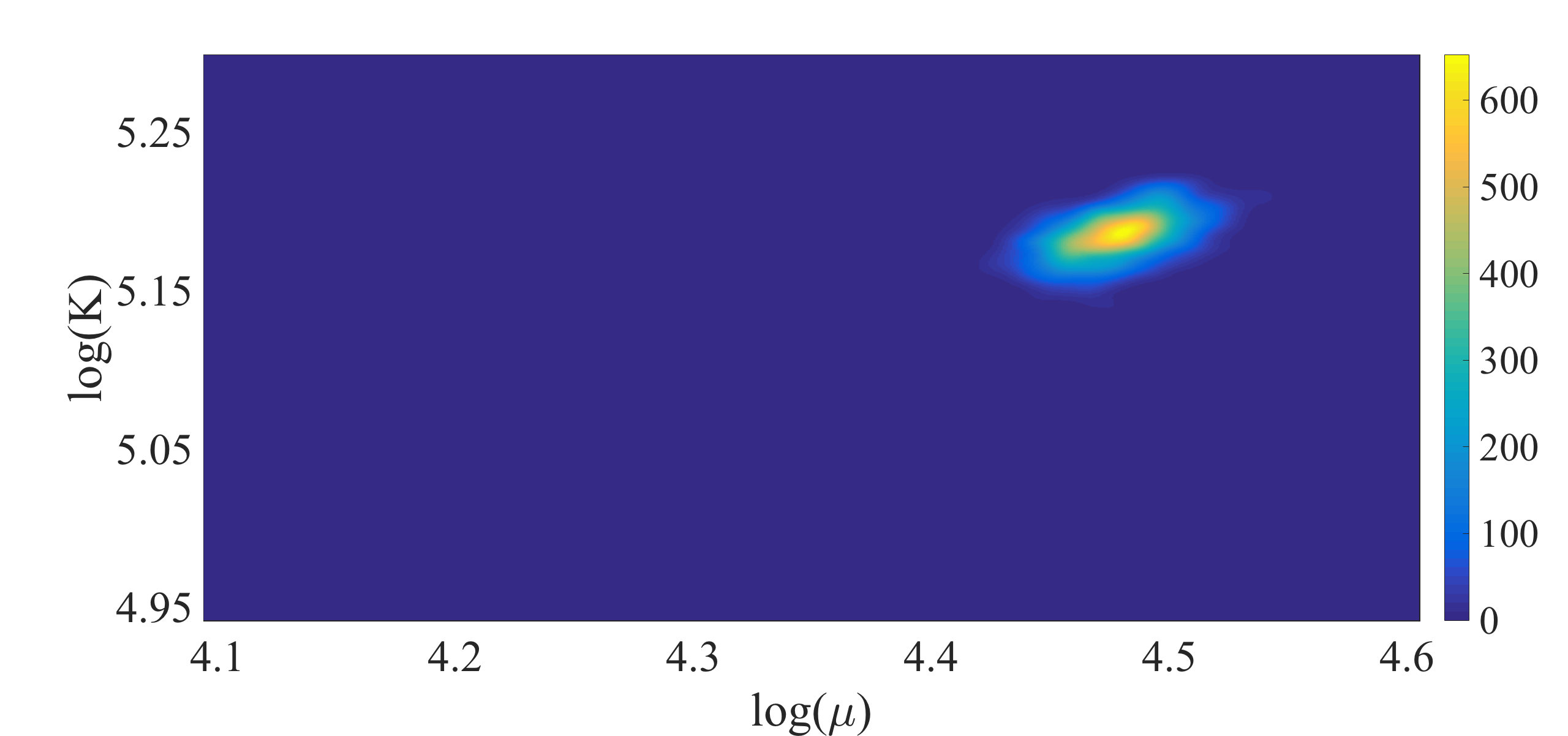} \\ \includegraphics[width=5.5cm,height=3cm]{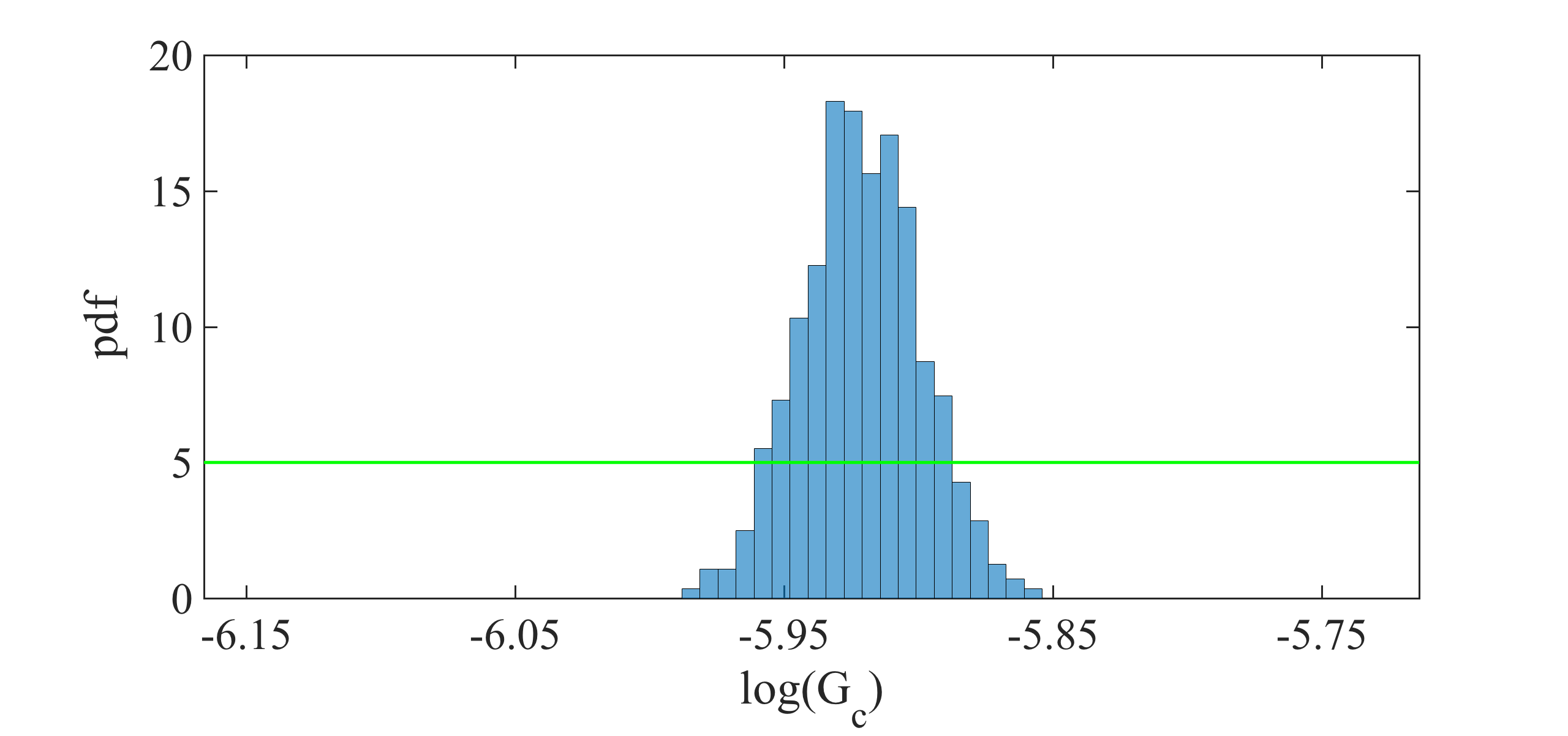} \end{minipage}
	\end{tabular}
	\caption{The load-displacement curves for the reference, prior, and
		posterior values (left panel) for $h=1/80$ in the
		SENT example (Example 1). The \new{joint and the} marginal posterior distributions of the effective parameters are shown in the right panel.}
	\label{fig:Ex1_h180}
\end{figure}

Here we strive to solve the problem using a coarse mesh and employ
MCMC to find parameters that make the solution more precise compared
with the reference value. Figure \ref{fig:Ex1_h120} shows the obtained
displacement with both prior and posterior distribution for
$h=1/20$. The efficiency of the Bayesian estimation is pointed out
here since the peak point and the failure point are estimated
precisely. The posterior distributions are shown in the right
panel as well. The estimation can also be performed for finer meshes:
Figure \ref{fig:Ex1_h140} and Figure \ref{fig:Ex1_h180} illustrate the
load-displacement curves for $h=1/40$ and $h=1/80$, respectively. In
both cases, in addition to the precise estimation of the
crack-initiation point and the material-failure point, the curve is
closer to the reference value. Again, the posterior distributions are
shown on the right panels. Finally, the mean values of the posterior
distributions in addition to their acceptance rates are indicated in
Table \ref{Exam1Table}.

\begin{table}[ht]
	\begin{tabular}{lllllll}
		\hline
		& $\mu$ &   rate (\%)	 & $\new{K}$&   rate (\%)   &$G_c$ &   rate (\%) \\
		\hline
		$h=1/20$     &  \new{83,9}   & \new{21} &  \new{190.4}    & \new{21} & \new{0.00221} & \new{23}     \\
		$h=1/40$     &  \new{92.4}   & \new{20} &  \new{197.2}    &  \new{20} & \new{0.00271} & \new{26}     \\
		$h=1/80$     &  \new{88.1}   & \new{27} &   \new{178.4}    &\new{27} & \new{0.00266} & \new{28.1}     \\
		\hline
	\end{tabular}
	\caption{The mean of thel posterior distributions of
		$\mu$, $\new{K}$, and $G_c$ in the SENT example (Example 1) for $h=1/20$, $h=1/40$, and $h=1/80$. All units are in $\mathrm{kN/mm^2}$.}
	\label{Exam1Table}
\end{table}

\begin{table}[ht]
	\begin{tabular}{lllllll}
		\hline
		&mesh size	& $h=1/20$  &   $h=1/40$ 	 & $h=1/80$ &  $h=1/160$   &$h=1/320$    \\
		\hline
		&	CPU time [s]    &  10   & 23 &  1428    & 4810 & 15\,854     \\
		\hline
	\end{tabular}
	\caption{The elapsed CPU time for the estimation of the load-displacement diagram (with the reference values) until the failure point in SENT.}
	\label{time}
\end{table}

\subsection{Example 2. Double edge notch tension (DENT) test}

\begin{figure}[ht!]
	\centering
	\includegraphics[width=15cm,height=6cm]{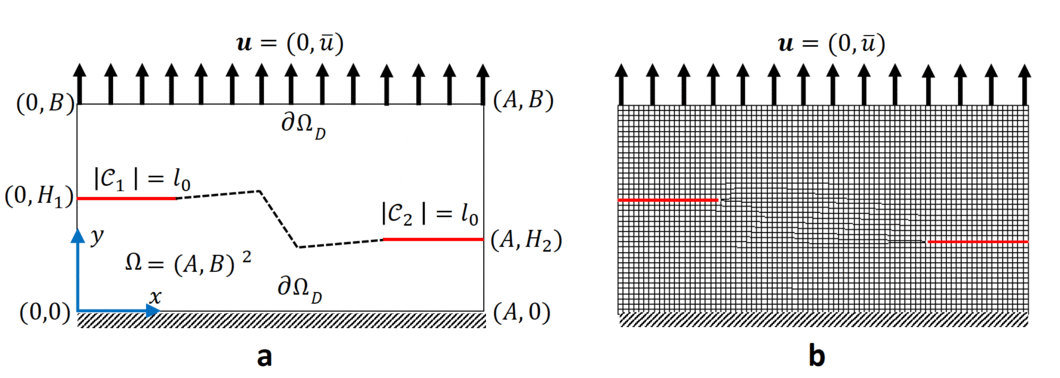}
	\caption{Schematic diagram for the DENT example  (left) and its corresponding mesh with $h=1/80$ (right).}
	\label{schematics2}
\end{figure}

This numerical example is a fracture process that occurs through the
coalescence and merging of two cracks in the domain. We consider the
tension test with a double notch located on the left and right
edge. The specimen is fixed on the bottom. We have traction-free
conditions on both sides. A non-homogeneous Dirichlet condition is
applied to the top-edge. The domain has a predefined two-notch located in
the left and right edge in the body as shown in Figure
\ref{schematics2}a. We set $A:=\unit{20}{mm}$ and $B:=\unit{10}{mm}$
hence $\Omega=\unit{(20,10)^2}{mm^2}$. For the double-edge-notches,
let $H_1:=\unit{5.5}{mm}$ and $H_2:=\unit{3.5}{mm}$ with the
predefined crack length of $l_0:=\unit{5}{mm}$ (Figure
\ref{schematics2}a). This numerical example is computed by imposing a
monotonic displacement $\bar{u}=1\times10^{-4}$ at the top surface of
the specimen in a vertical direction. The finite element
discretization that uses $h=1/80$ is indicated in Figure
\ref{schematics2}b.

\begin{figure}[ht]
	\centering
	\subfloat{\includegraphics[width=8cm,height=5cm]{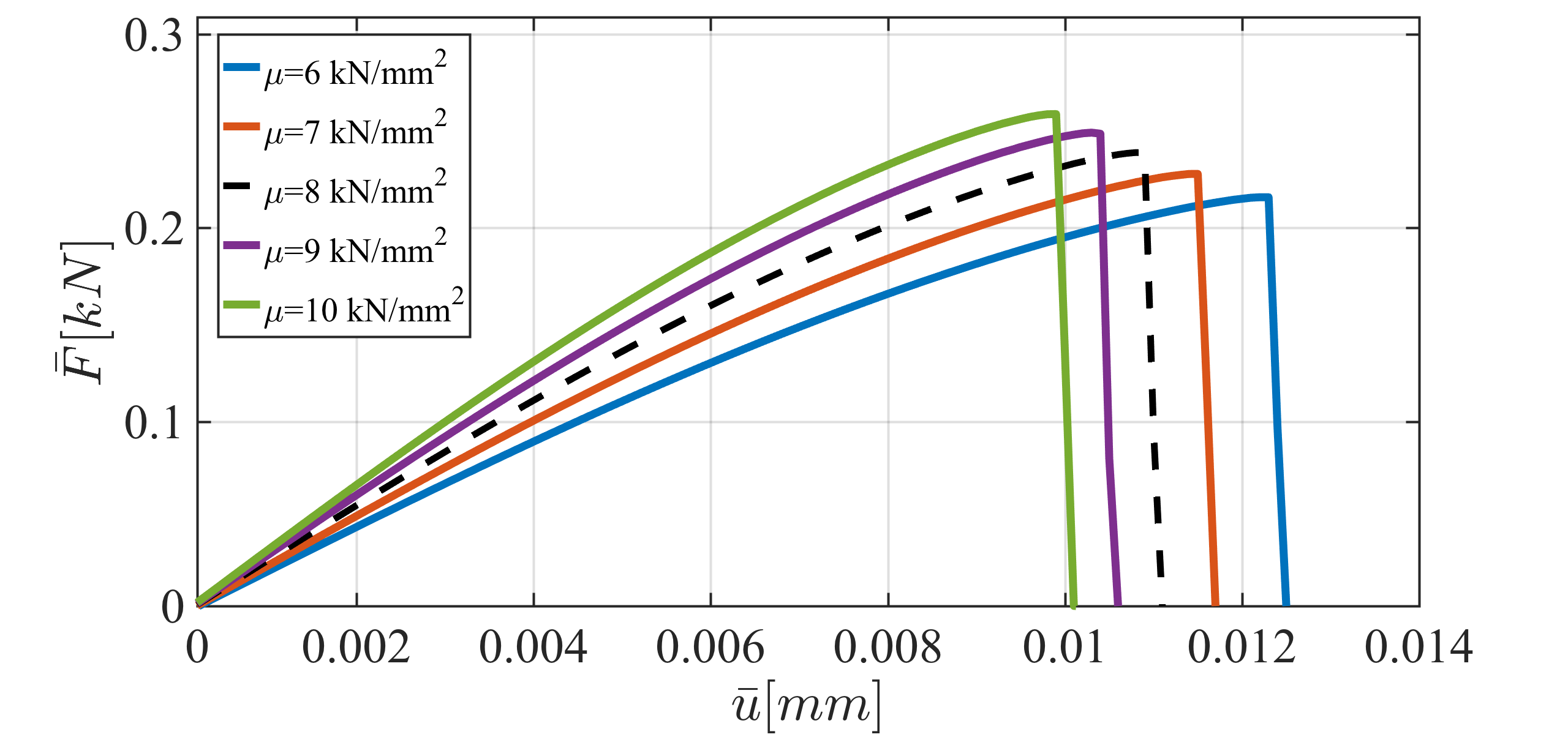}}%
	\hfill 
	\subfloat{\includegraphics[width=8cm,height=5cm]{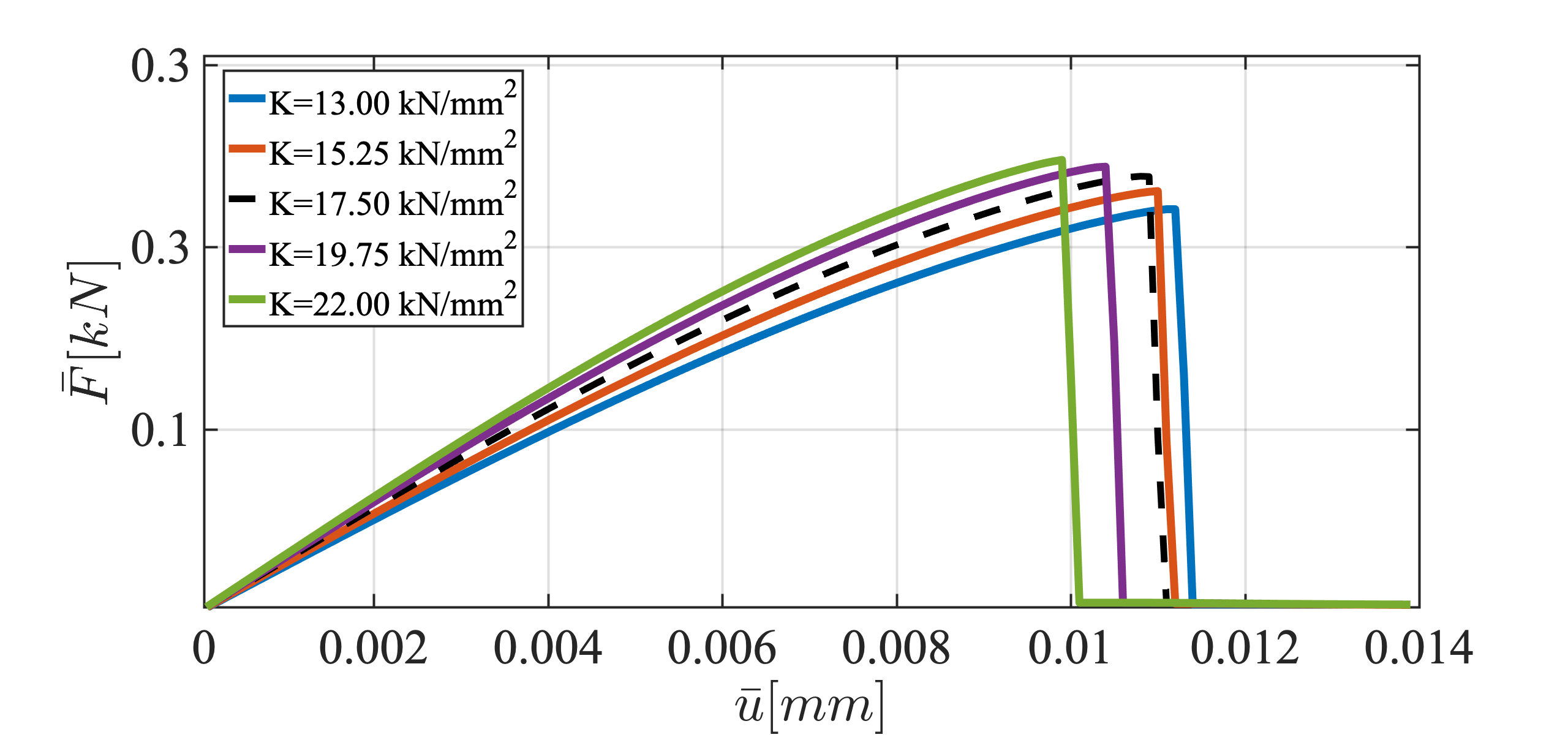}}%
	\newline
	\subfloat{\includegraphics[width=8cm,height=5cm]{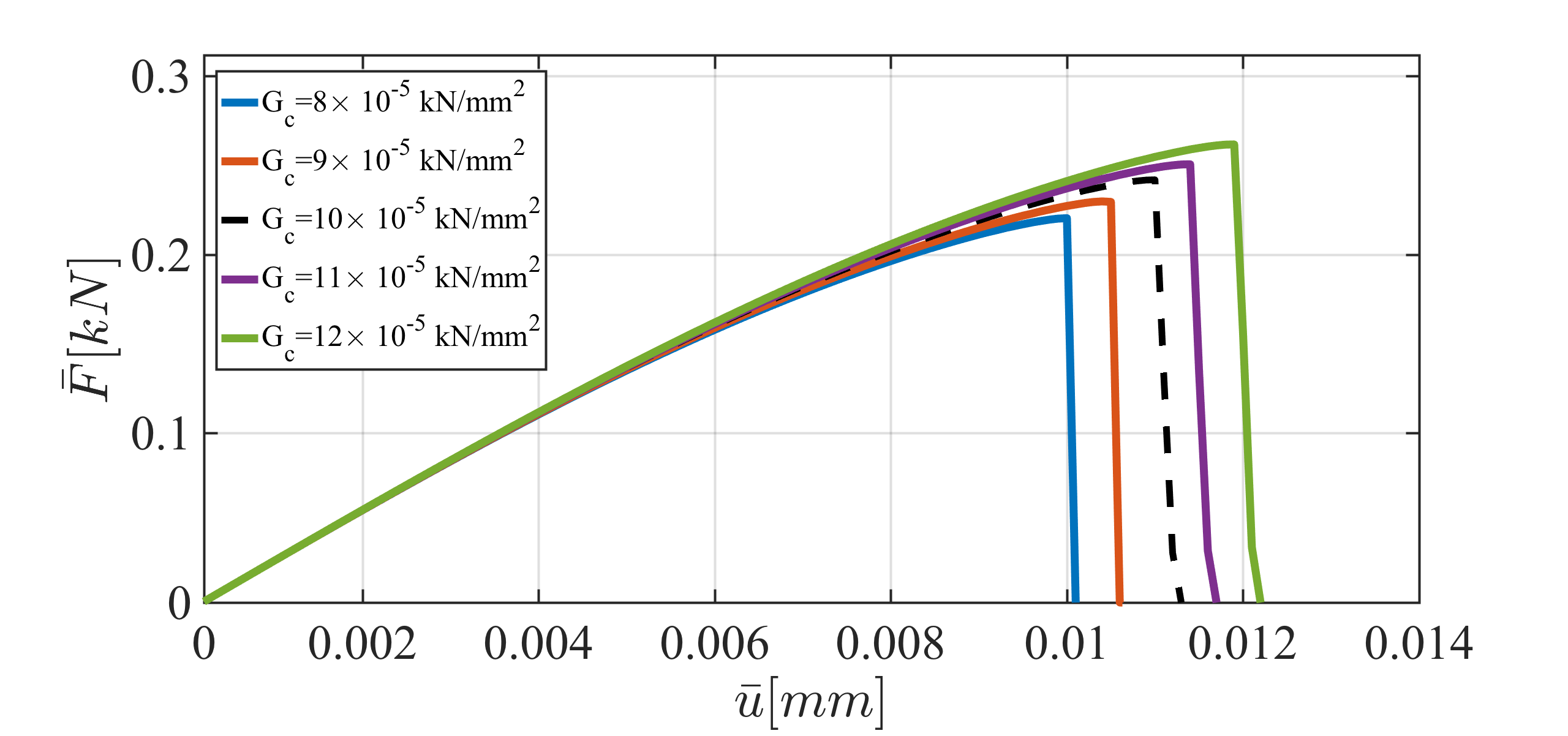}}%
	\caption{The load-displacement curve for different values of
		$\mu$  (top left), \new{$K$} (top right) and $G_c$ (bottom)
		for the DENT example.}
	\label{fig:exam2_parameters}
\end{figure}

\begin{figure}[h!]
	\subfloat{\includegraphics[width=7cm,height=4.5cm]{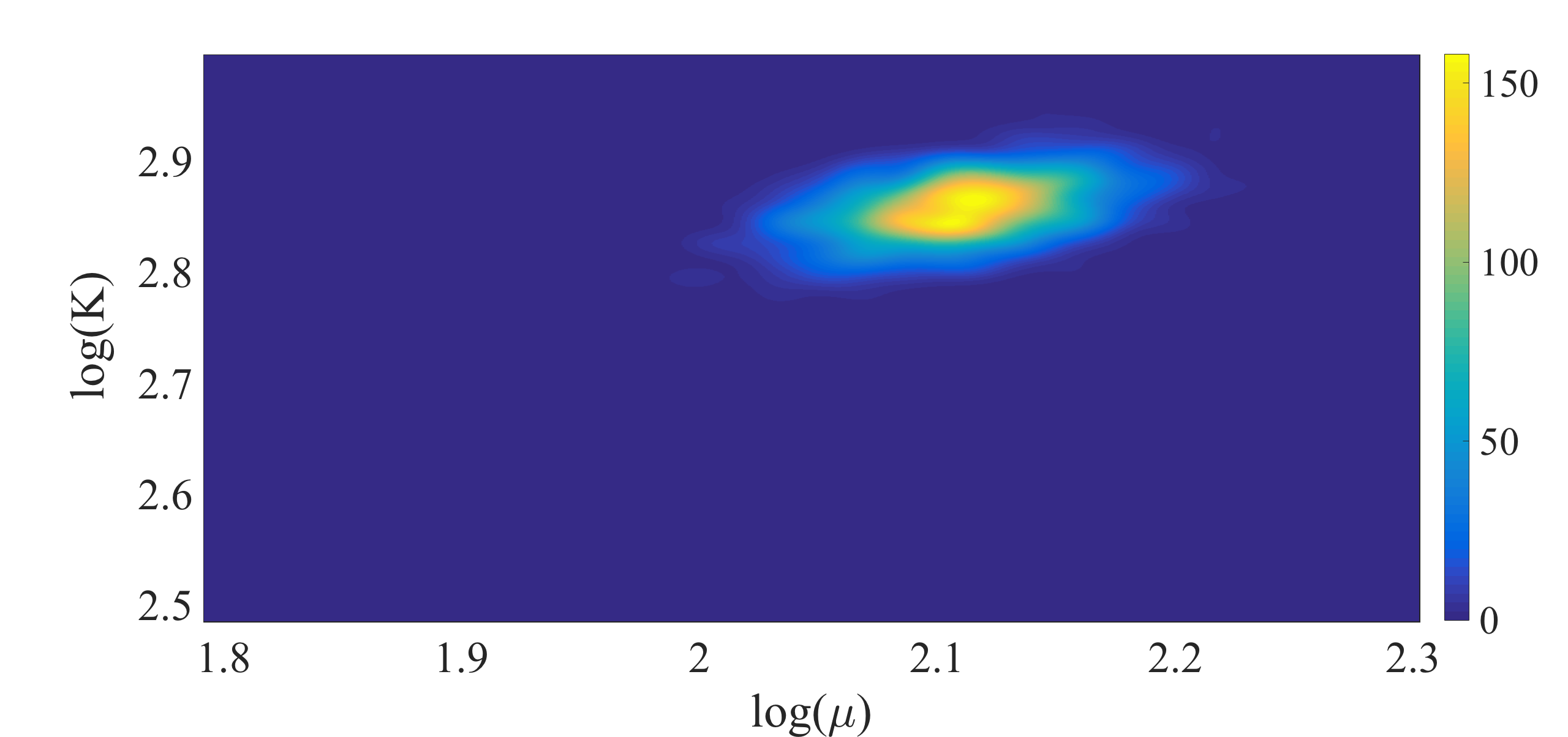}}%
	\subfloat{\includegraphics[width=7cm,height=4.5cm]{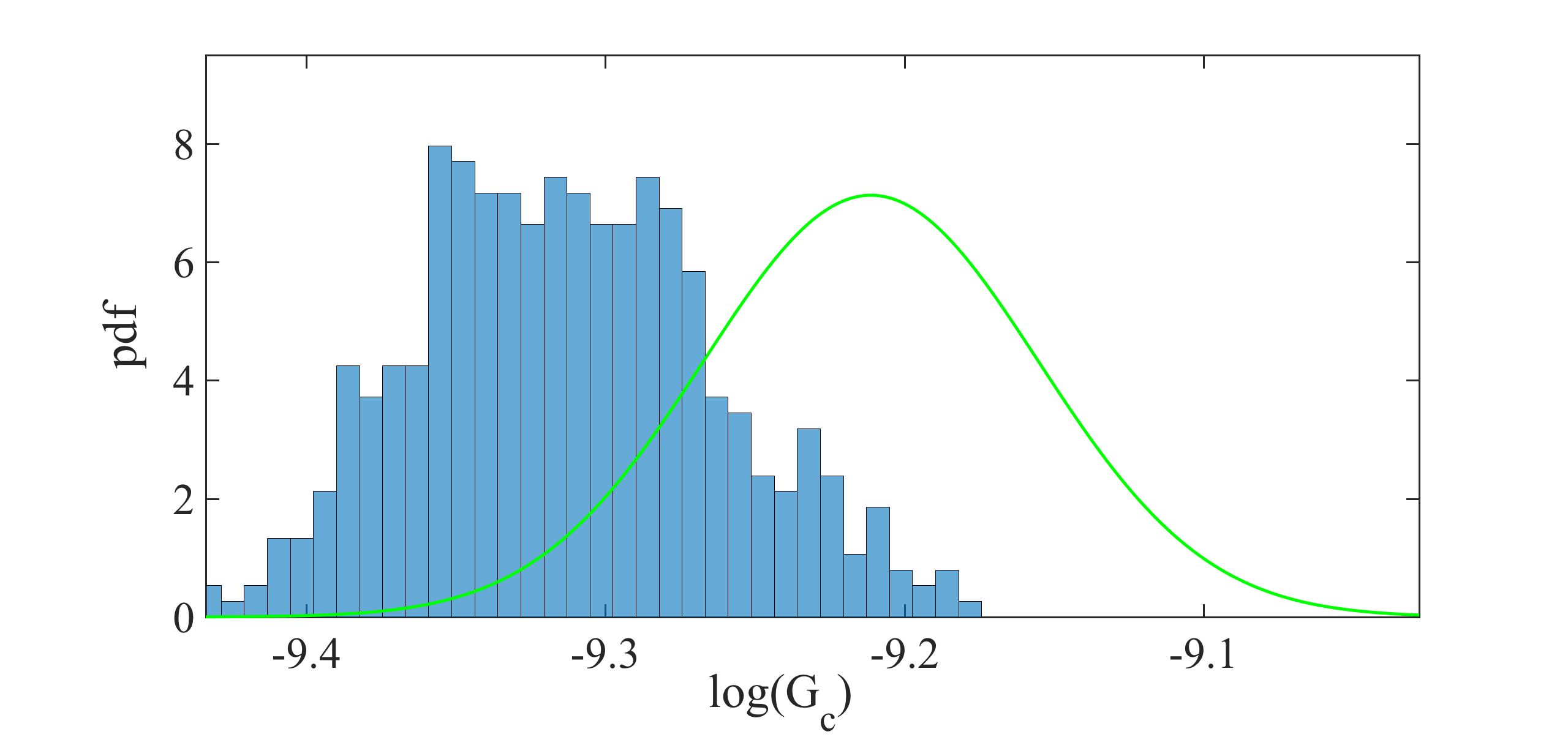}}
	\newline
	\subfloat{\includegraphics[width=7cm,height=4.5cm]{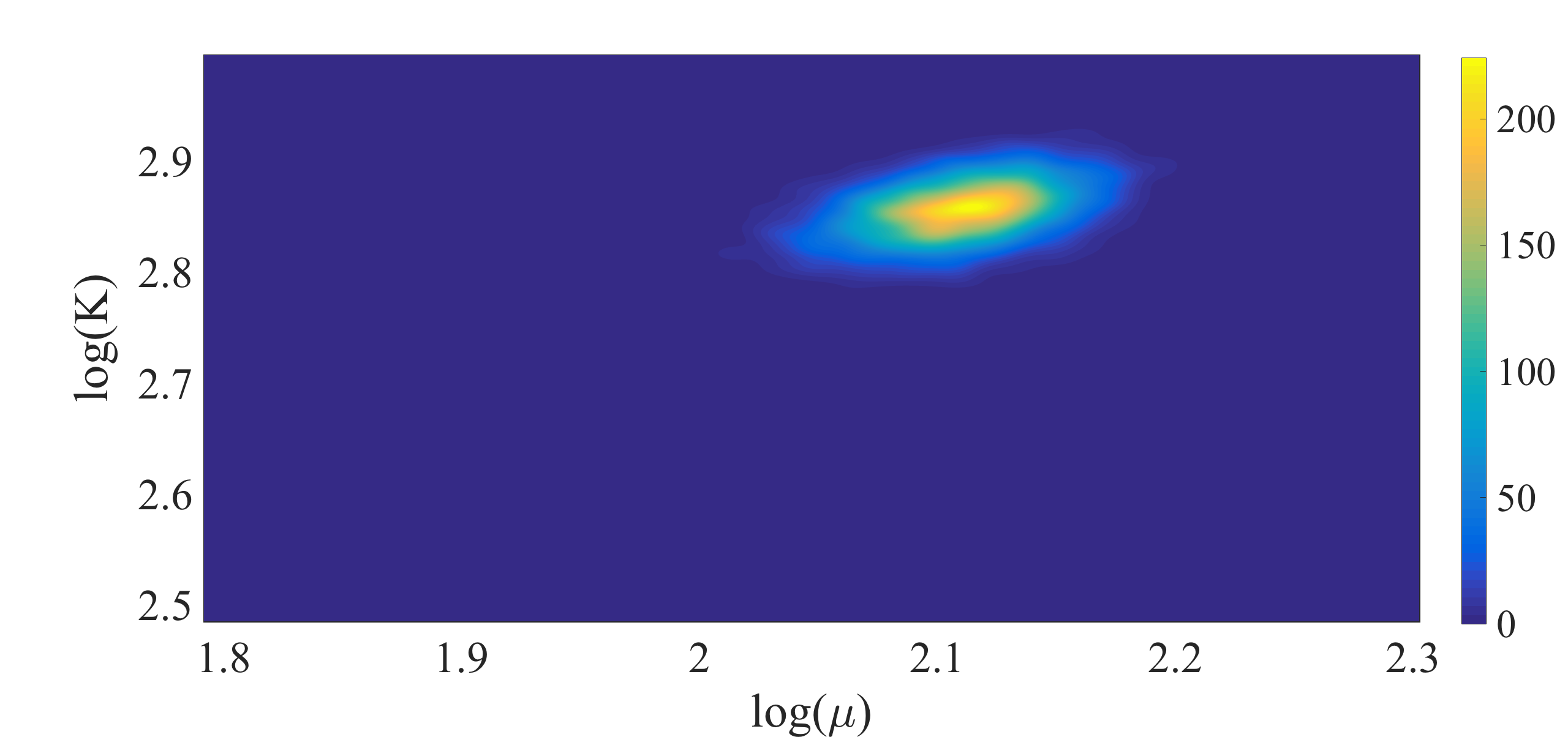}}%
	\subfloat{\includegraphics[width=7cm,height=4.5cm]{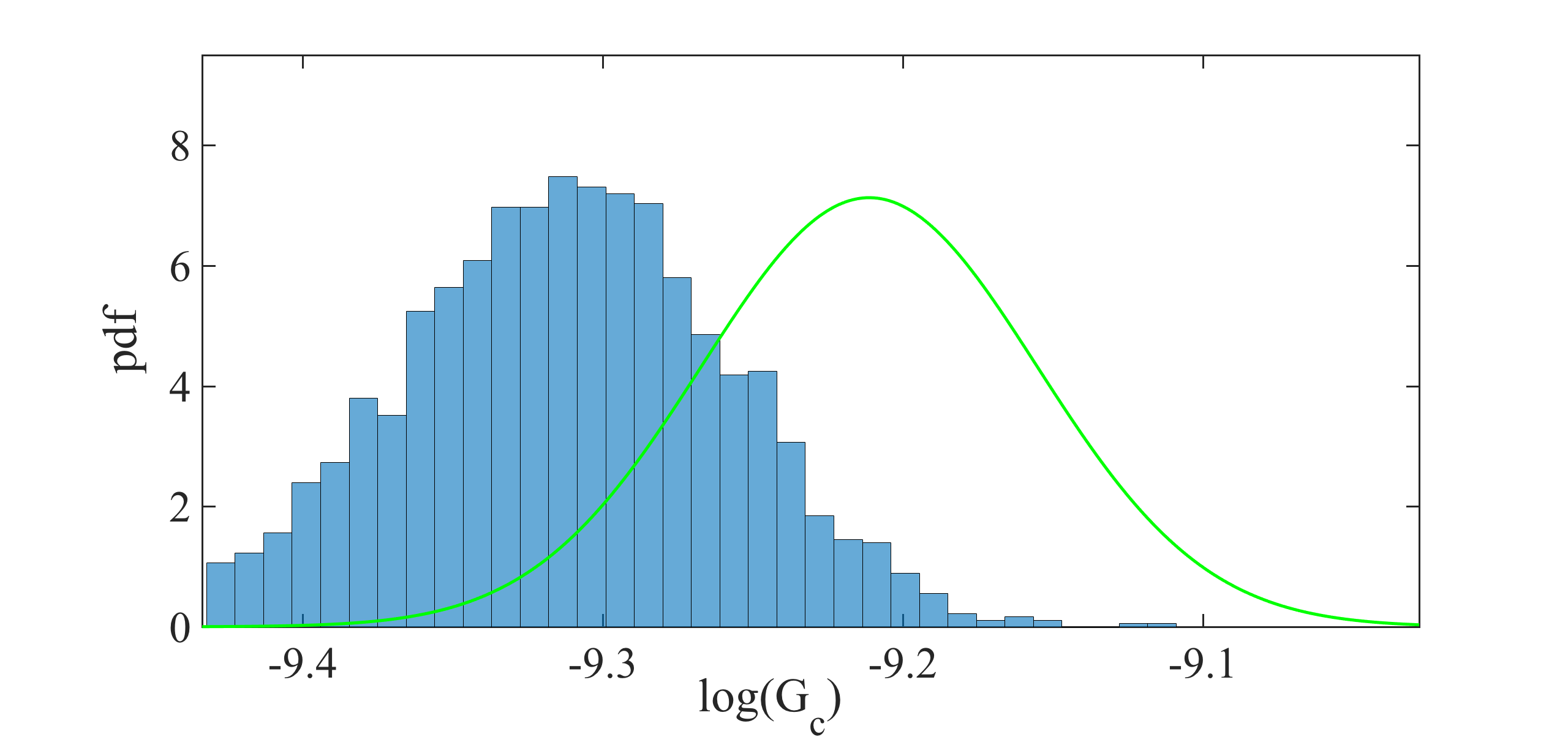}} 
	\hspace{-2cm}
	\newline
	\subfloat{\includegraphics[width=7cm,height=4.5cm]{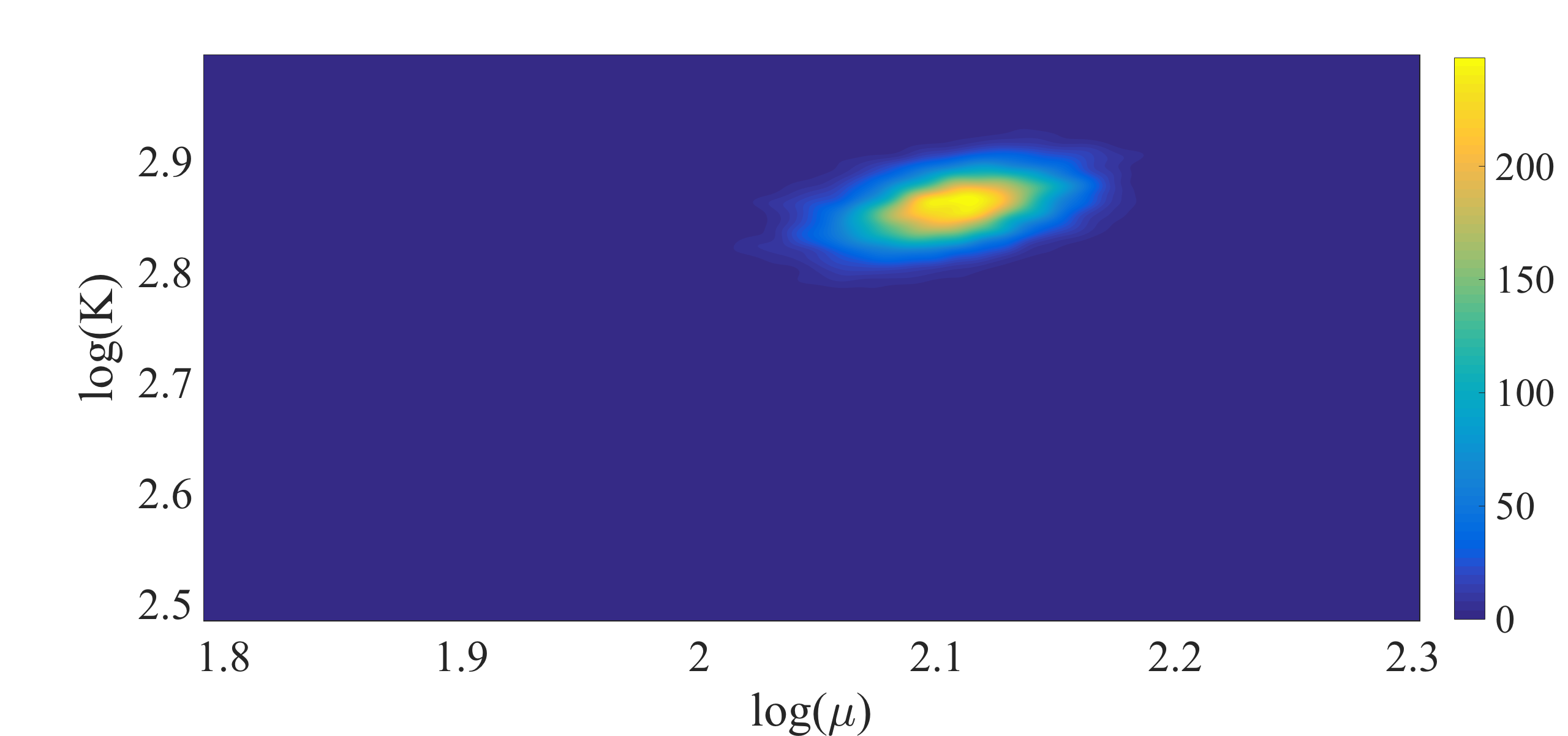}}%
	\subfloat{\includegraphics[width=7cm,height=4.5cm]{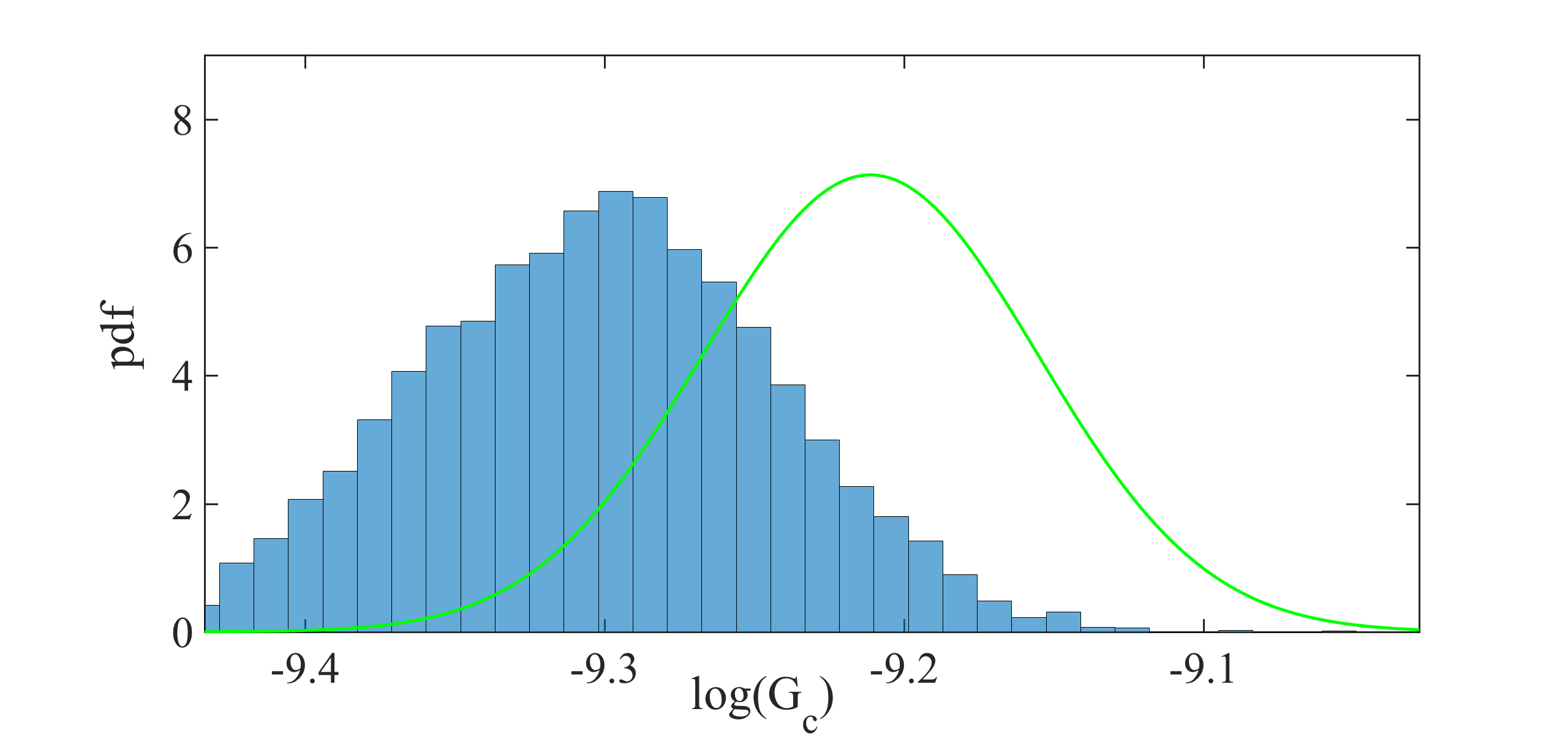}} 
	\caption{\new{Left: the joint probability density of  the elasticity parameters. Right:} the  
		prior (green line) and the posterior (histogram) of $G_c$  for the DENT example (Example 2). For the posterior distribution, we used $3\,000$ samples (the first row), $15\,000$ samples (the second row), and $50\,000$ samples (the third row).}
	\label{fig:exam2_hist}
\end{figure}

According to the truncated KL-expansion, for the \new{bulk} modulus
$K$, Eq.~\ref{KL2} gives the mean value of \new{$\bar{K}=23.58$}
and the standard deviation of \new{$\sigma_{K}=0.28$}. Therefore, the
parameter varies between $\unit{10}{kN/mm^2}$ and
$\unit{14}{kN/mm^2}$. For the shear modulus, the expectation of $
\bar{\mu}=22.8$ and the standard deviation of $\sigma_\mu=0.23$ leads
to the variation range
$(\unit{6}{kN/mm^2},\unit{10}{kN/mm^2})$. Similarly, by using a
KL-expansion for $G_c$, we obtained the variation range between
$\unit{8\times 10^{-5}}{kN/mm^2}$ and $\unit{12\times
	10^{-5}}{kN/mm^2}$.  Figure \ref{fig:exam2_parameters} illustrates
the effect of their different values \new{including shear and bulk modulus and $G_c$} on the curve.

We assumed the \new{uniform} proposal distribution, namely the normal distribution
\begin{align}
\mathcal{K}(\theta\rightarrow \theta^*):=\frac{1}{\sqrt{2\pi \sigma^2}}\exp\left(-\frac{(\theta-\theta^*)^2}{2\sigma^2} \right).
\end{align}
\new{As we aforementioned, the random field can be represented using the KL-expansion. In this example \eqref{KL} is employed to parameterize the elasticity and energy rate parameters. The random perturbations are imposed on the $\xi$ coefficients in the KL expansion. According to the proposal, the mean of the KL-expansion is updated.}

Here we plan to study the effect of the number of samples $N$ on the
posterior distribution. Figure \ref{fig:exam2_hist} shows the \new{joint} distributions of ($\new{K}$, $\mu$), and \new{the marginal distribution of} $G_c$  using $N=3\,000$, $N=12\,000$, and $N=50\,000$. The
calculations are done with $h=1/80$, and $h=1/320$ is used as the
reference.  As shown, with a larger number of samples, the
distribution is close to a normal distribution. Table \ref{Table:Ex2}
points out the mean values in addition to the acceptance rate of all influential parameters. 


\begin{figure}[t!]
	\centering
	\hspace{-0.0cm}\subfloat{\includegraphics[width=5.85cm,height=4cm]{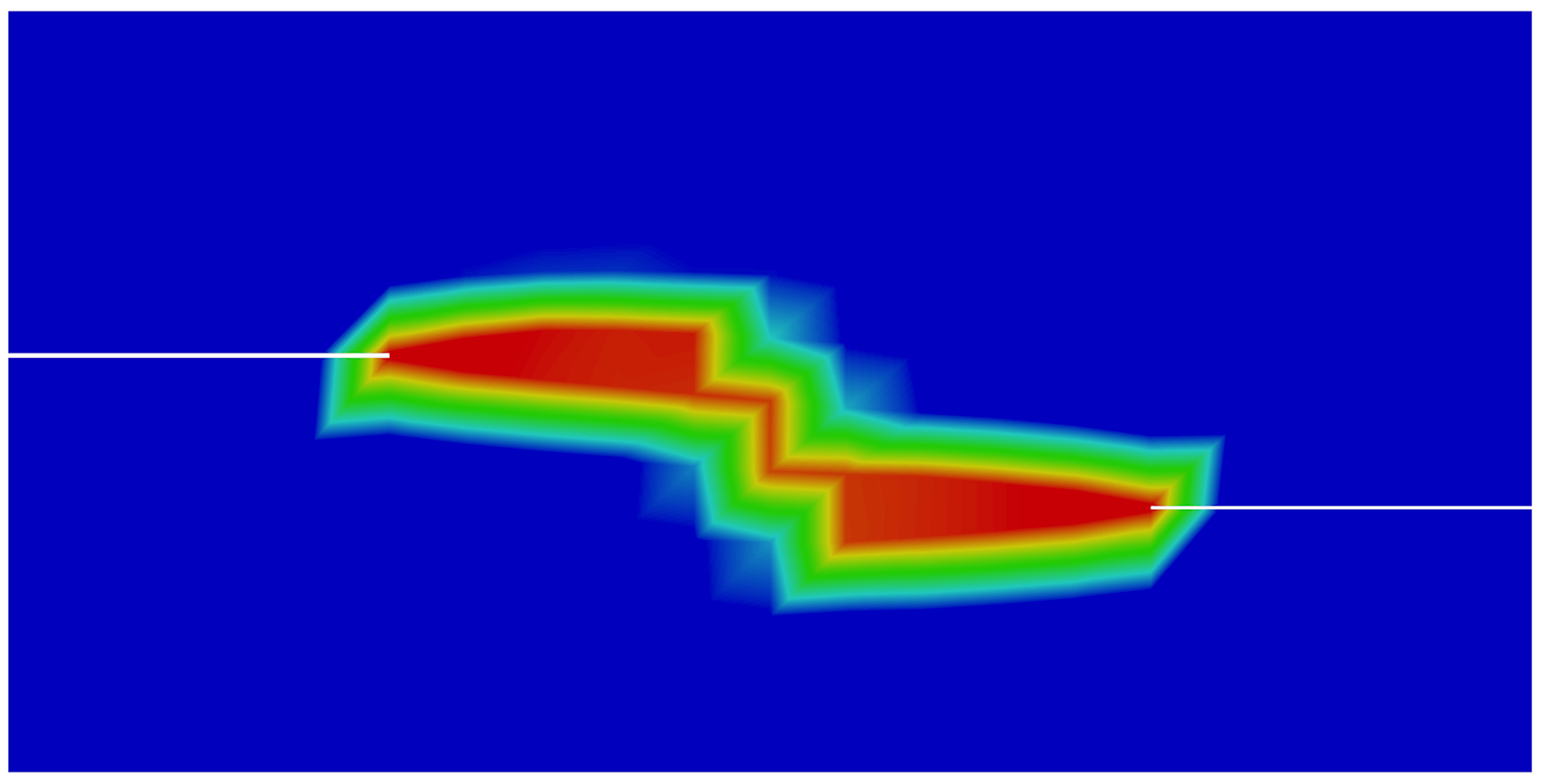}}\hspace{1cm}
	\subfloat{\includegraphics[width=5.85cm,height=4cm]{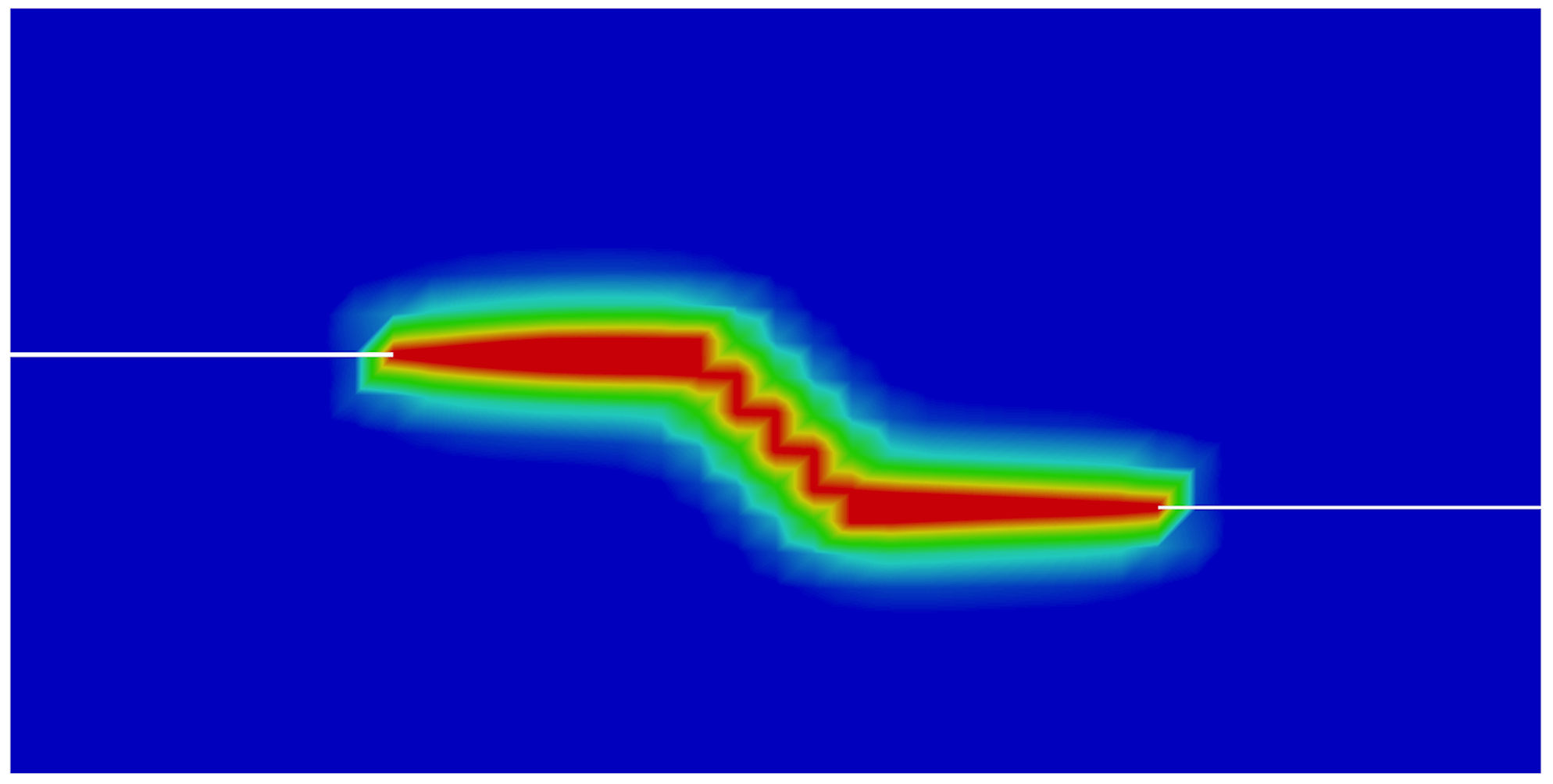}}\newline\\
	\subfloat{\includegraphics[width=5.85cm,height=4cm]{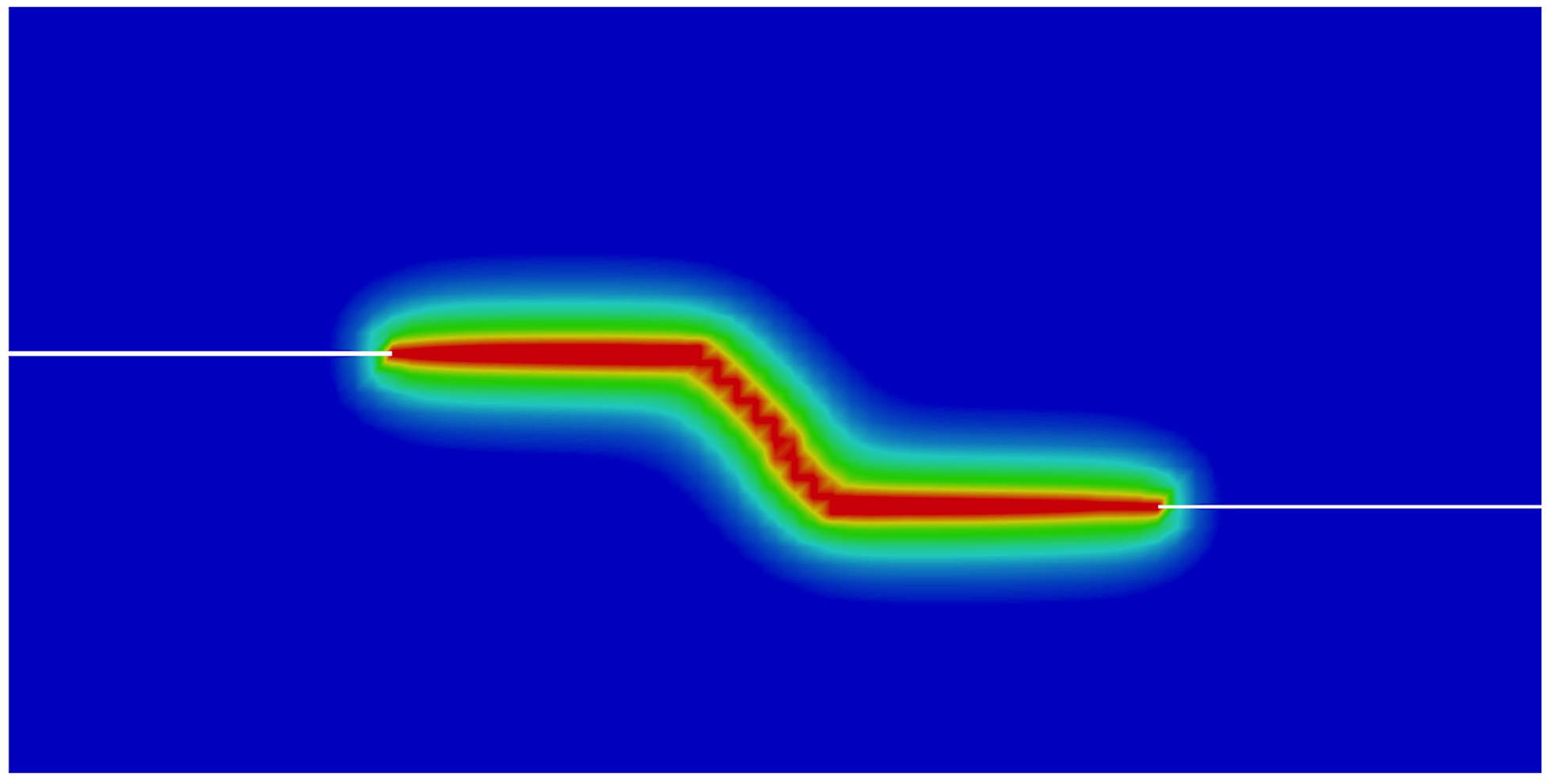}}\hspace{1cm}
	\subfloat{\includegraphics[width=5.85cm,height=4cm]{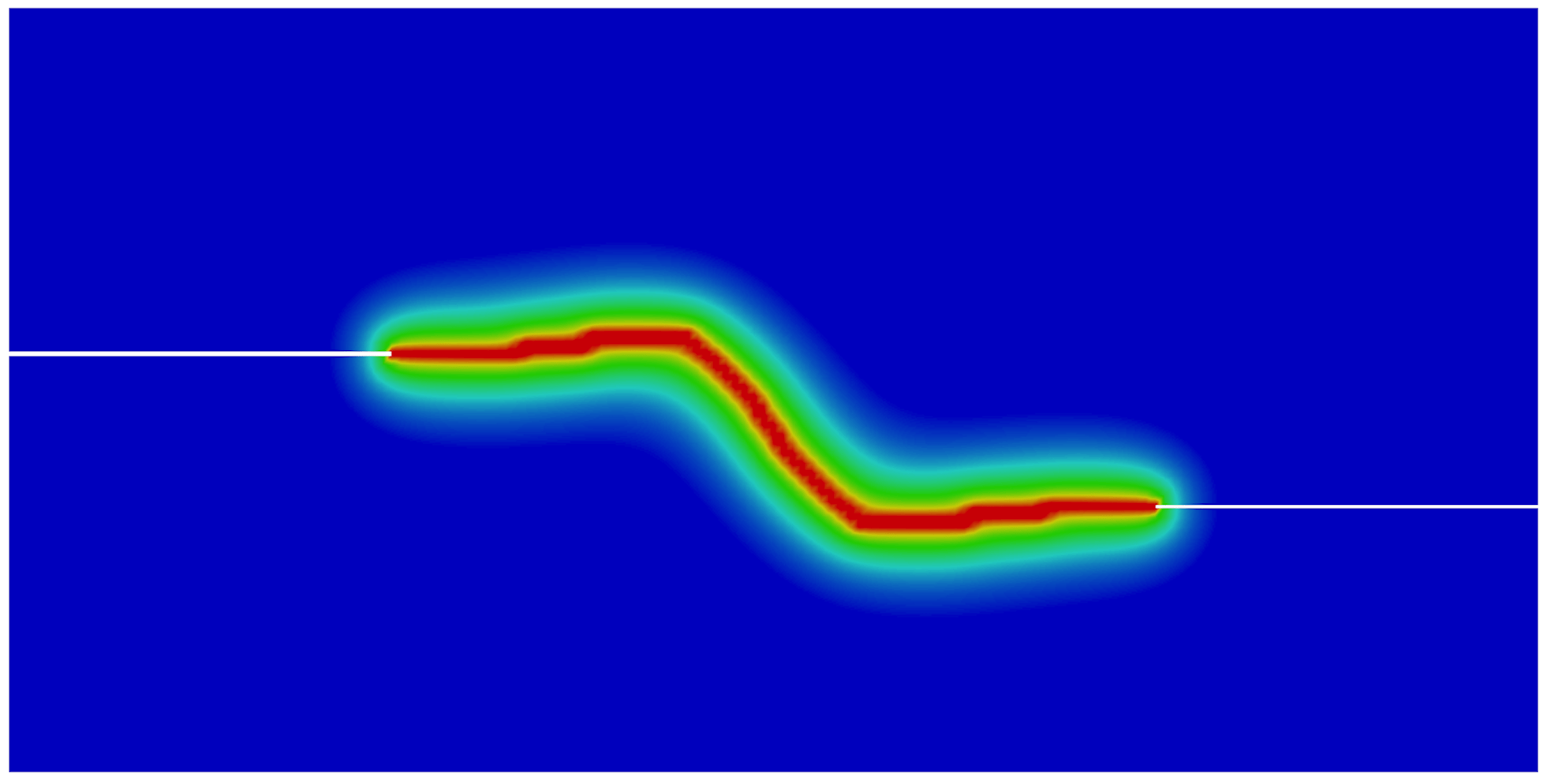}}\newline\\
	\hspace{0.9cm}\subfloat{\includegraphics[width=5.85cm,height=4cm]{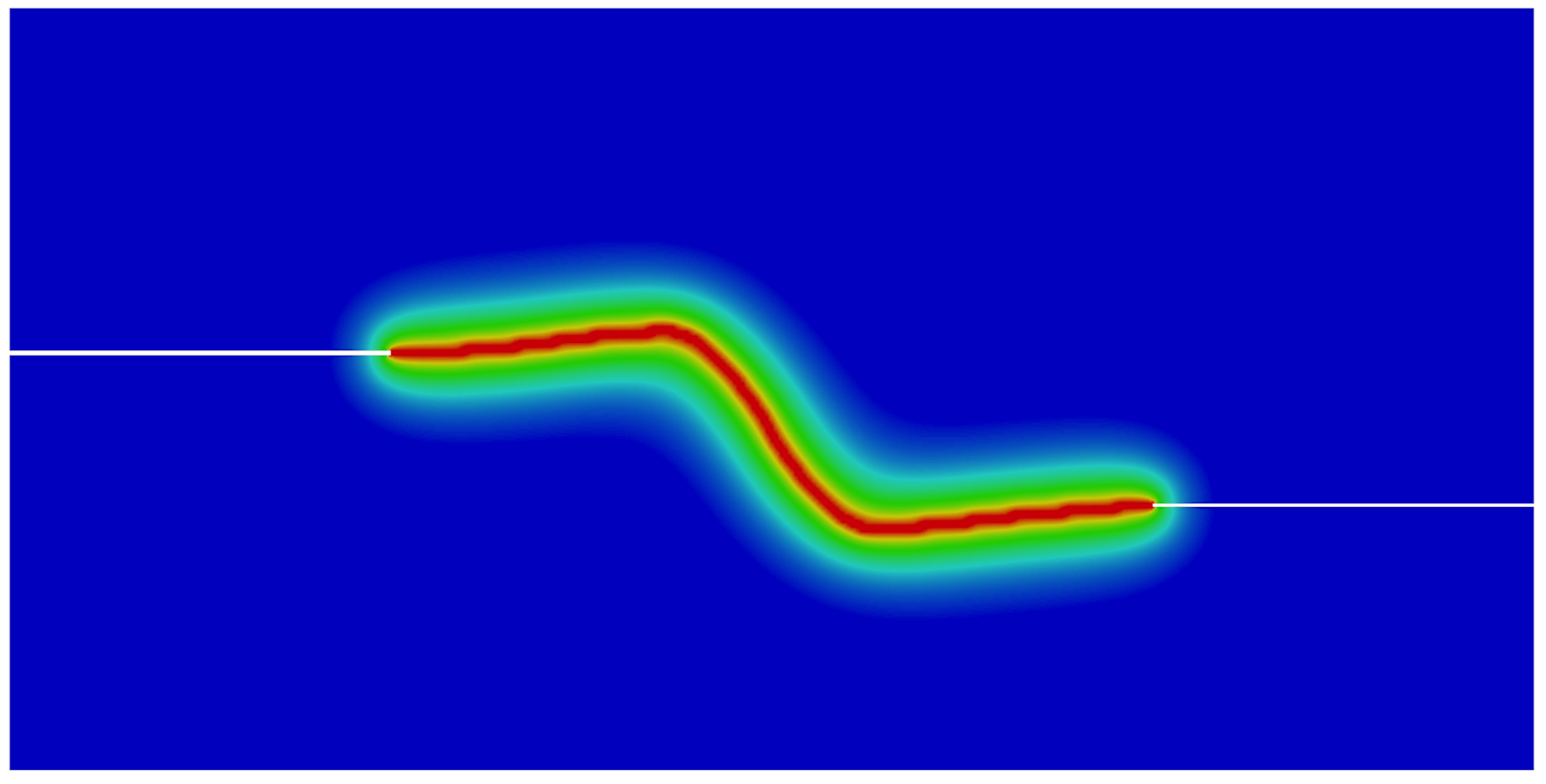}}\hspace{1cm} 
	\subfloat{\includegraphics[width=5.85cm,height=4cm]{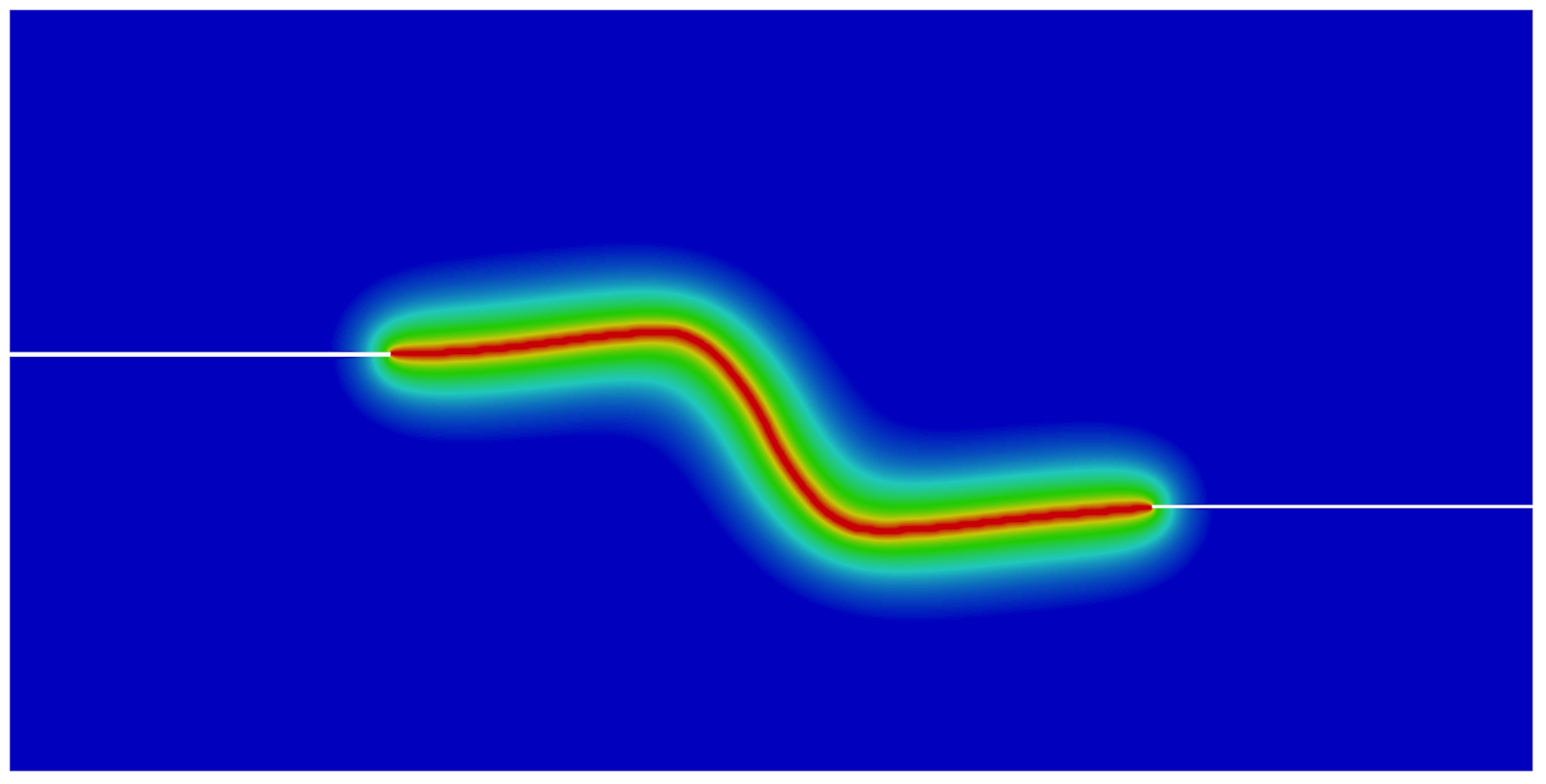}}\hspace{0.2cm} 
	\subfloat{\includegraphics[width=2cm,height=4cm]{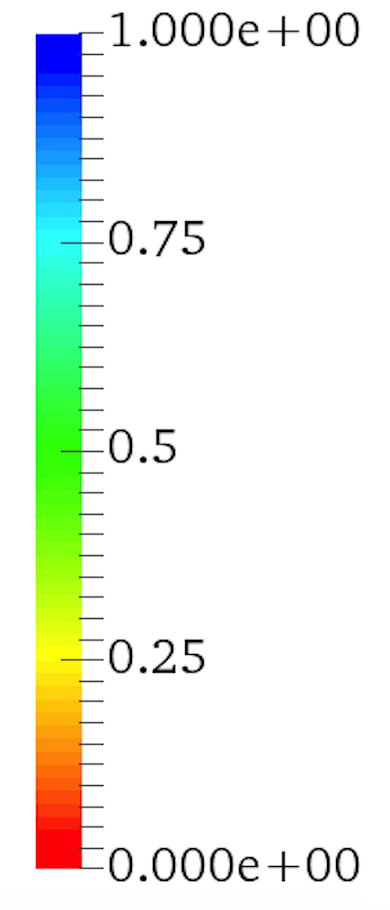}}  
	\caption{The effect of the mesh size on the crack propagation
		in the DENT example (Example 2). The mesh sizes are (from  the left) $h=1/10$, $h=1/20$, $h=1/40$, $h=1/80$, $h=1/160$, and $h=1/320$ (the reference).}
	\label{fig:exam2_mesh1}
\end{figure}

\begin{table}
	\begin{tabular}{lllllll}
		\hline
		& $N_\text{samples}$=3\,000&& $N_\text{samples}$=15\,000 && $N_\text{samples}$=50\,000\\
		\cline{2-3}\cline{4-5}\cline{6-7}
		& mean ($\mathrm{kN/mm^2}$) &   rate (\%)	 & mean ($\mathrm{kN/mm^2}$)&   rate (\%)   & mean ($\mathrm{kN/mm^2}$) &   rate (\%) \\
		\hline
		$\mu$     &   \new{8.26}  & \new{31.0} &  \new{8.23}    & \new{31.5} & \new{8.35} & \new{36.2}     \\
		$\new{K}$   &  \new{17.4}  & \new{31.0}  & 	\new{17.42}     &\new{31.5}  & \new{17.12} & \new{36.2}    \\
		$G_c$       & \new{$9.02\times 10^{-5}$}     & \new{15.2}   & 	\new{$9.1\times 10^{-5}$}    & \new{15.8} & \new{$9.68\times 10^{-5}$} & \new{16.7}   \\
		\hline
	\end{tabular}
	\caption{The mean value and the acceptance rate of the
		posterior distributions of $\mu$, \new{$K$}, and $G_c$ with
		$N_\text{samples}=3\,000$, $N_\text{samples}=12\,000$, and
		$N_\text{samples}=50\,000$ in the DENT example (Example 2). All units are in $\mathrm{kN/mm^2}$.}
	\label{Table:Ex2}
\end{table}

\begin{figure}[t!]
	\centering
	\subfloat{\includegraphics[width=8.2cm,height=5.5cm]{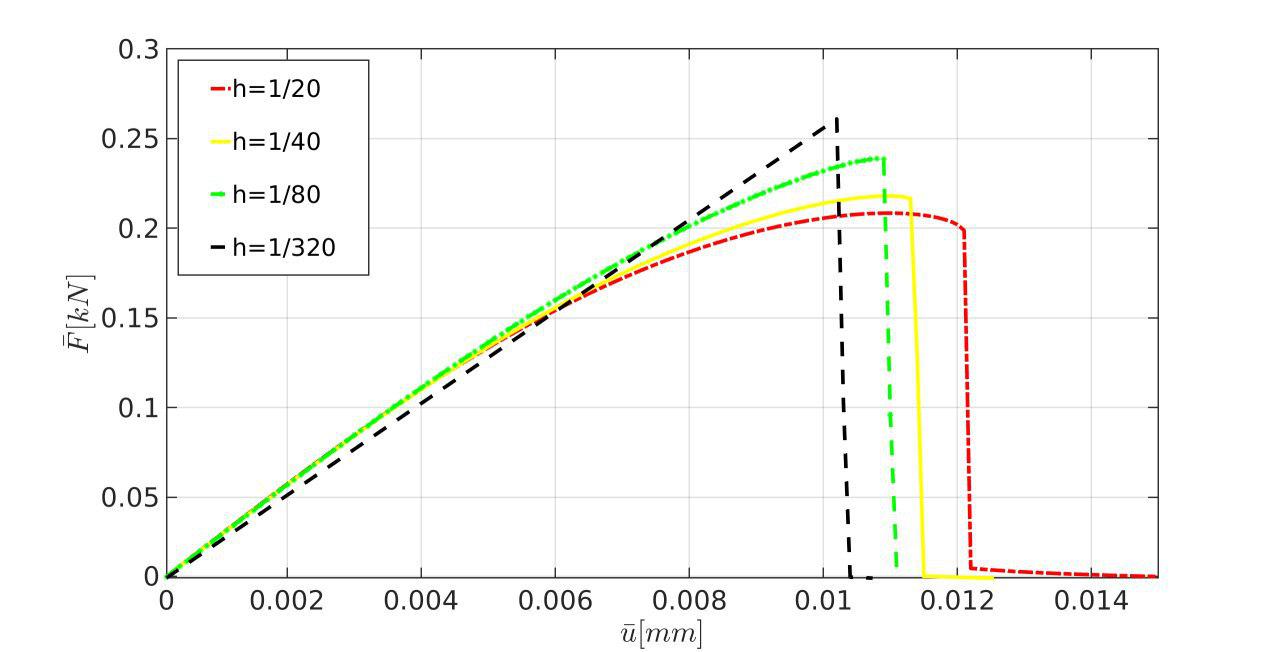}}%
	\hfill 
	\subfloat{\includegraphics[width=8.2cm,height=5.5cm]{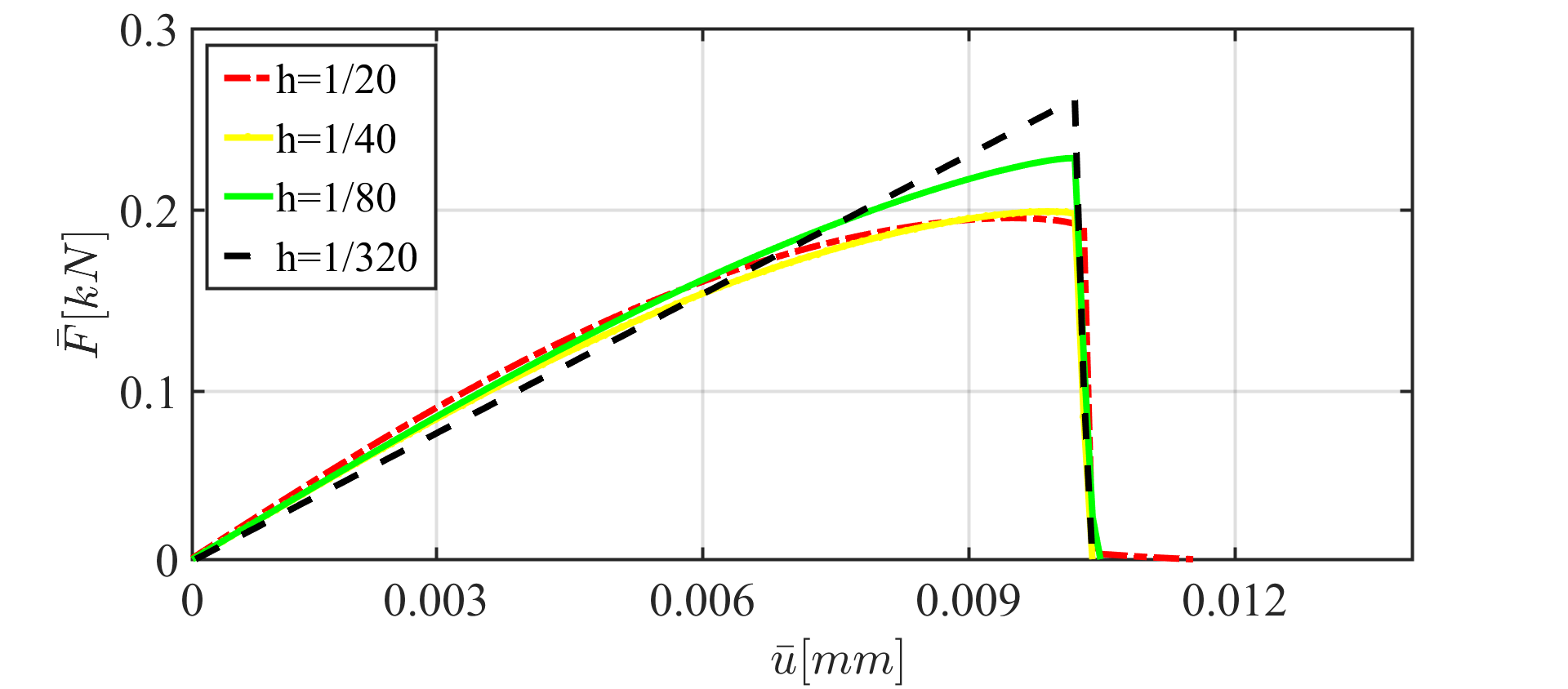}}%
	\caption{The load-displacement curve of DENT (Example 2) with
		different mesh sizes. Here the parameters are chosen
		according to the prior (left) and posterior (right) distributions.}
	\label{Ex2:pror_psot}
\end{figure}

As the next step, we use different mesh sizes for the Bayesian 
inversion using 15\,000 samples. Figure \ref{fig:exam2_mesh1} shows the crack pattern using
different mesh sizes changing from $h=1/10$ to $h=1/320$. Finer meshes
lead to a smoother and more reliable pattern.  Figure
\ref{Ex2:pror_psot} depicts the load-displacement diagram using the
prior values. With coarse meshes, the curve is significantly different
from the reference including crack initiation. Using Bayesian
inversion (see Figure \ref{fig:exam2_hist}) enables us to predict the
crack propagation and initiation more precisely. As the figure shows,
even for the coarsest mesh (compare $h=1/10$ to $h=1/320$) the peak
and fracture points are estimated precisely. For finer meshes (e.g.,
$h=1/80$) the diagram is adjusted tangibly compared to the reference
value. Finally, a summary of the mean values (of posterior
distributions) and their respective acceptance rate is given in
Table \ref{Ex2_Table}.

\begin{table}
	\begin{tabular}{lllllll}
		\hline
		& $\mu$ &   rate (\%)	 & $\new{K}$&   rate (\%)   &$G_c$ &   rate (\%) \\
		\hline
		$h=1/20$     &  \new{8.90}   & \new{33.3} &  \new{18.08}    & \new{33.3} & \new{$8.15\times 10^{-5}$} & \new{13}     \\
		$h=1/40$     &  \new{8.15}   & \new{37.5} &  \new{17.61}    & \new{37.5} &  \new{$8.23\times 10^{-5}$} &  \new{15.2}     \\
		$h=1/80$     &  \new{8.35}   & \new{38} &  \new{17.4}    &  \new{38} &  \new{$9.10\times 10^{-5}$} & \new{26.4}     \\
		\hline
	\end{tabular}
	\caption{The mean of the posterior distributions of $\mu$,
		\new{$K$}, and $G_c$ for different mesh sizes in the DENT example. The units are in $\mathrm{kN/mm^2}$}
	\label{Ex2_Table}
\end{table}

The significant advantage of the developed Bayesian inversion is a
significant computational cost reduction. As shown, for SENT and DENT,
by using Bayesian inference for coarser meshes, the estimated
load-displacement curve is very close to the reference values. We
should note that the needed CPU time for $h=1/320$ is approximately
4.5 hours; however the solution with $h=1/80$ is obtained in less than
10 minutes.  This fact pronounces the computational efficiency
provided by Bayesian inversion, i.e., obtaining a relatively  precise
solution in spite of using much coarser meshes.

\subsection{Example 3. Tension test with two voids}
\label{voids}

\begin{figure}[ht!]
	\centering
	\includegraphics[width=17cm,height=8cm]{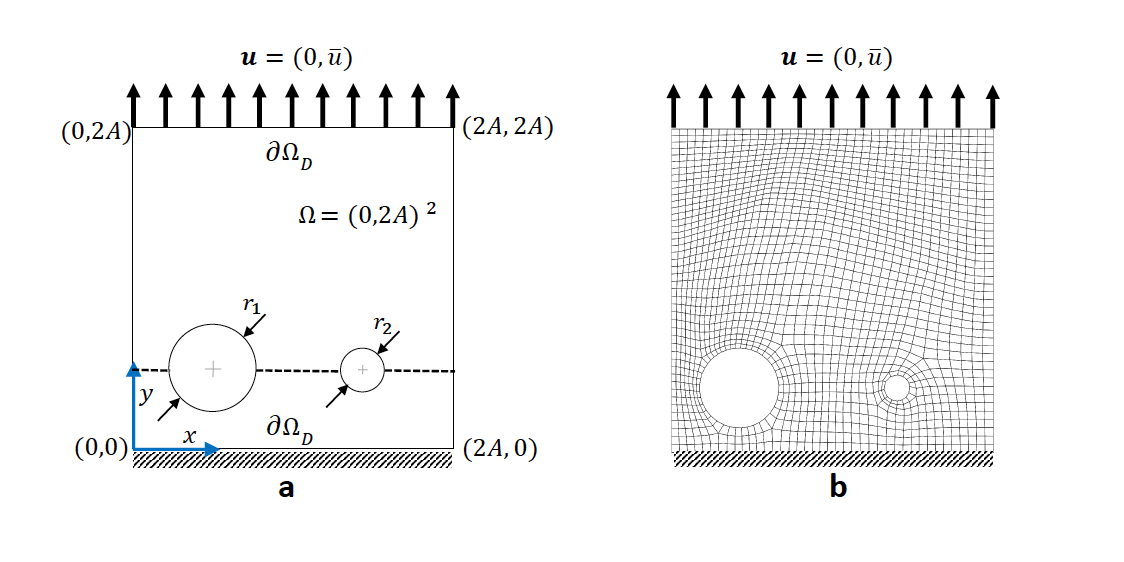}
	\vspace{-1cm}
	\caption{Schematic of SENT with voids (Example 3) (left) and its corresponding meshes with $h=1/40$ (right).}
	\label{schematics3}
\end{figure}

Here we consider the tension test where two voids are located in the
domain as a more complicated example.  The voids are used to weaken
the material and to lead to crack nucleation/initiation without an
initial singularity (i.e., a pre-existing crack).  The specimen is
fixed on the bottom. We have traction-free conditions on both sides. A
non-homogeneous Dirichlet condition is applied to the top. Domain
includes a predefined two voids in the body, as depicted in Figure
\ref{schematics3}a. We set $A=\unit{0.5}{mm}$ hence
$\Omega=\unit{(0,1)^2}{mm^2}$. The radius of left void is $r_1:=0.247
$ with the center $c_1:=(0.21,0.197)$. The radius of the right void is
$r_2:=0.0806$ with the center $c_2:=(0.7,0.197)$. This numerical
example is computed by imposing a monotonic displacement
$\bar{u}:=1\times10^{-4}$ at the top surface of the specimen in
vertical direction.  The finite-element discretization corresponding
$h=1/40$ is shown in Figure \ref{schematics3}b.

\begin{figure}[h!]
	\centering
	\subfloat{\includegraphics[width=8cm,height=5cm]{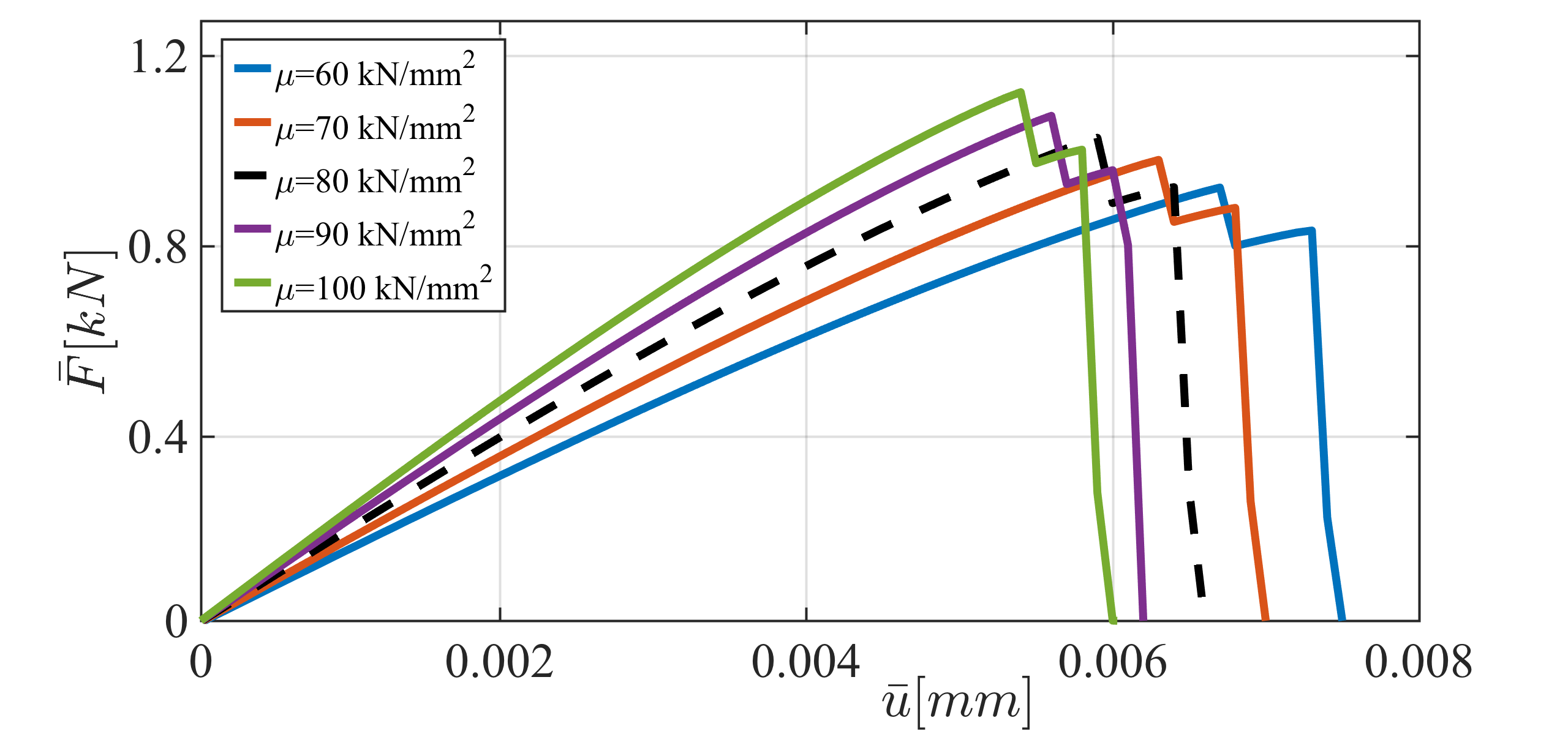}}%
	\hfill 
	\subfloat{\includegraphics[width=8cm,height=5cm]{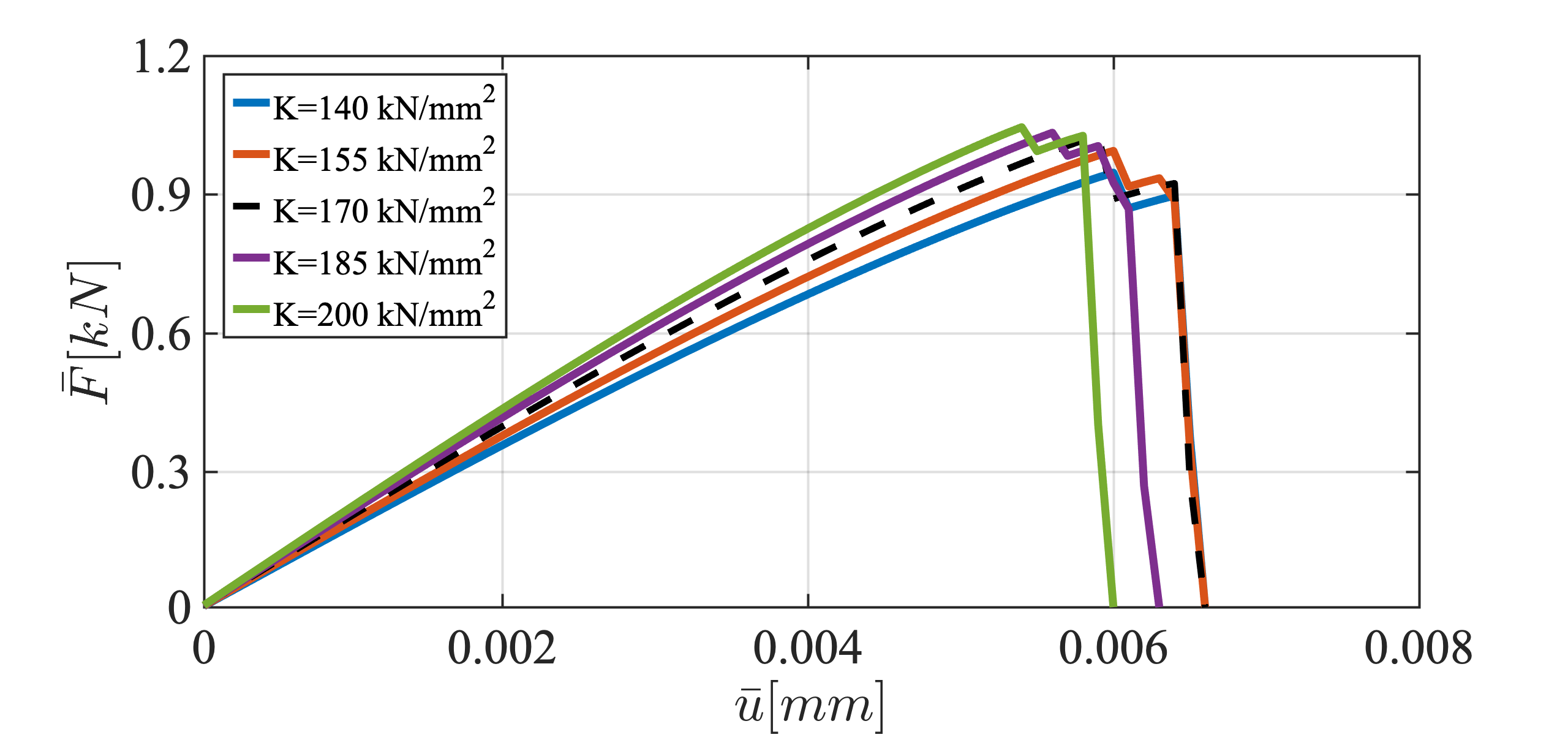}}%
	\newline
	\subfloat{\includegraphics[width=8cm,height=5cm]{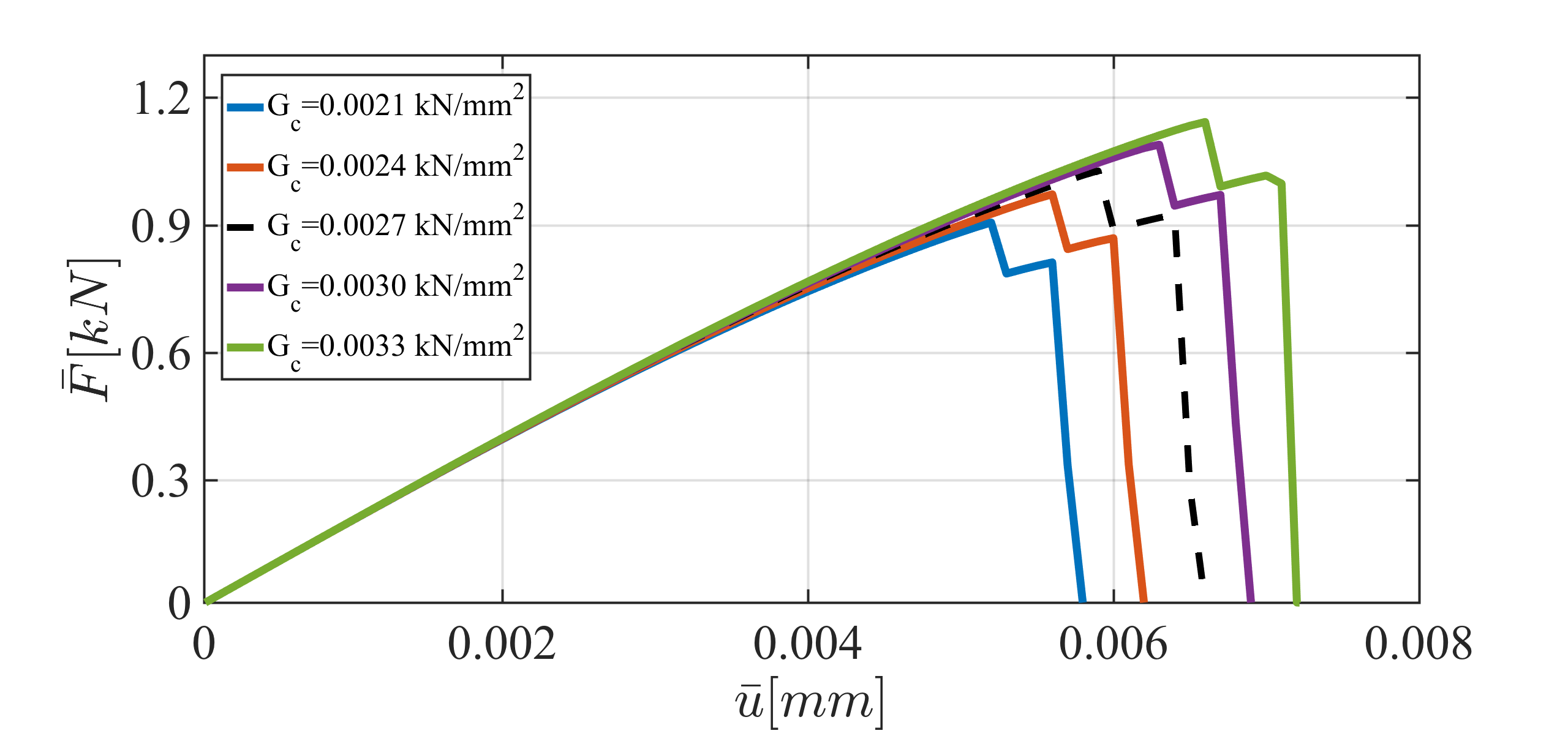}}%
	\caption{The load-displacement curve for different values of
		$\mu$ (top left), \new{$K$} (top right) and $G_c$ (bottom) for the SENT voids example (Example 3).}
	\label{fig:exam3_parameters}
\end{figure}

\begin{figure}[h!]
	\centering
	\subfloat{\includegraphics[width=7cm,height=4.5cm]{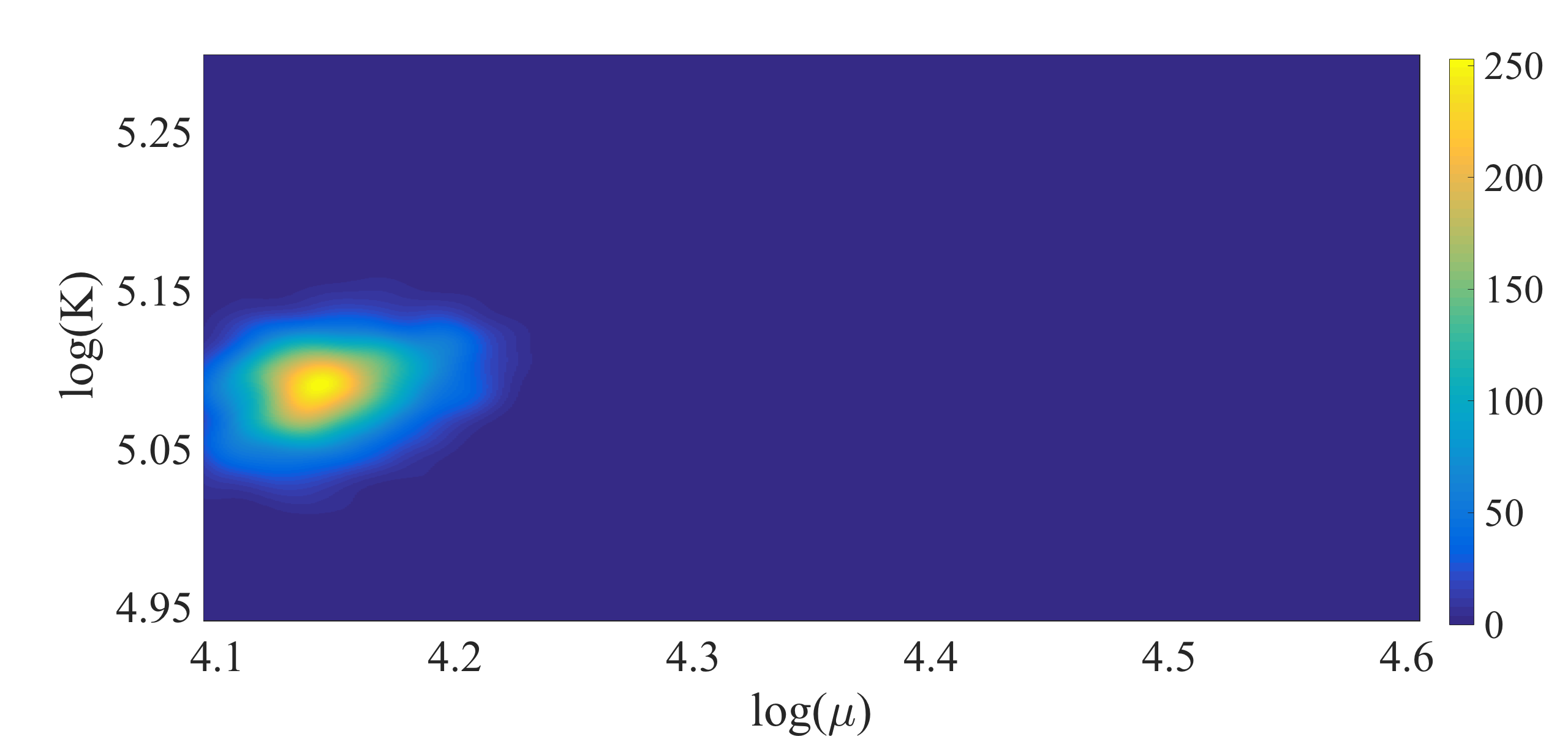}}%
	\subfloat{\includegraphics[width=7cm,height=4.5cm]{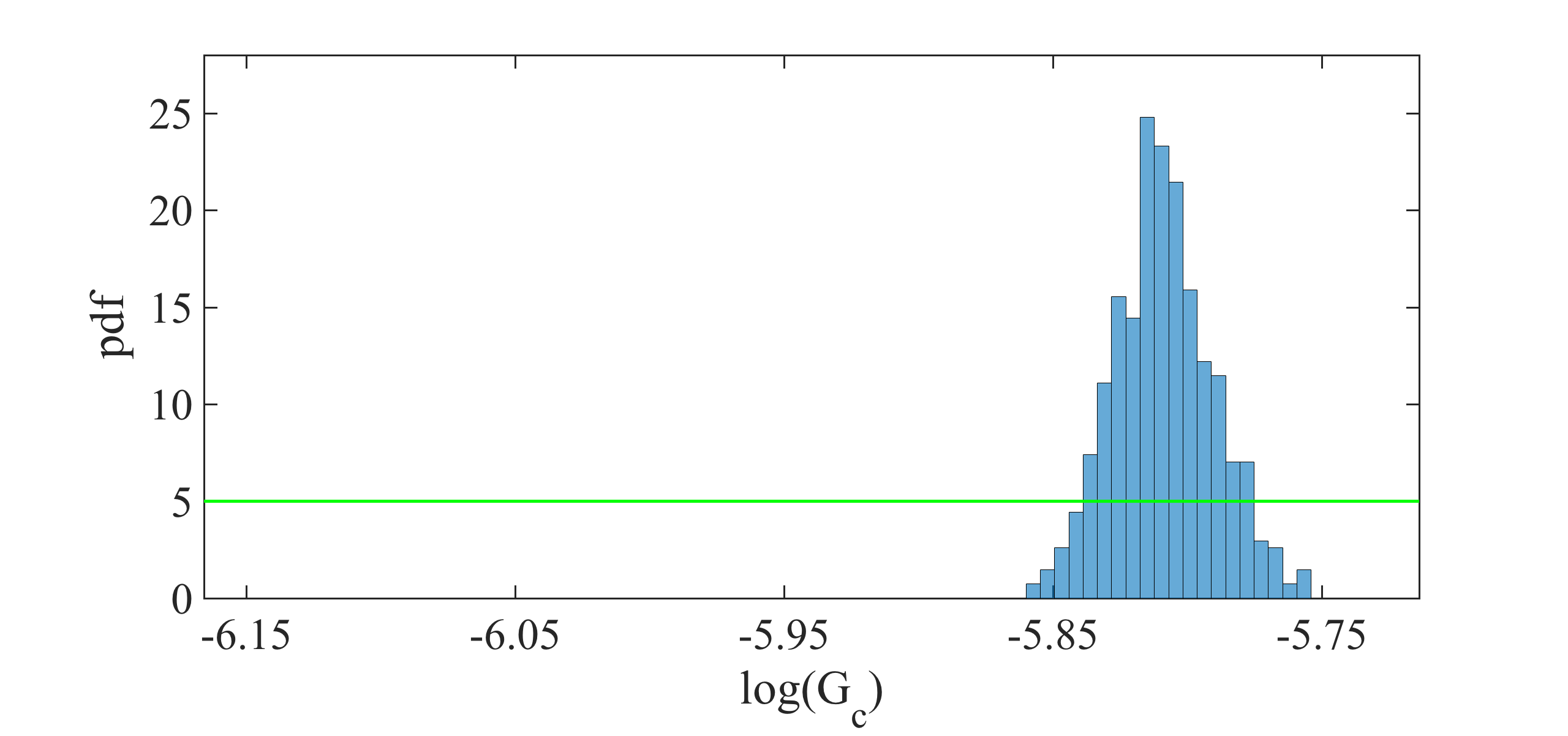}}
	
	\caption{\new{Left: the joint probability density of the elasticity parameters. Right:} the 
		prior (green line), and the posterior (histogram) distribution of  $G_c$ for SENT with voids (Example 3).}
	\label{fig:exam3_histogram}
\end{figure}

Due to the resemblance to the first example (SENT) we use the same range of parameters \new{  and the variables are again spatially constant random variables.}
The load-displacement curves obtained from different values of $\mu$, $\new{K}$, and $G_c$ are illustrated in Figure \ref{fig:exam3_parameters}.

This numerical example includes two voids results in multi-stage crack propagation. Hence, in the load-displacement curve, two peak points exist to demonstrate multi-stage crack propagation, see Figure \ref{fig:exam3_parameters}.

Figure \ref{fig:exam3_histogram} shows the proposal distribution where
a uniform prior distribution is used for Bayesian inversion with 10\,000 samples. Here we
use $h=1/160$ as the reference solution and $h=1/80$ is employed to
estimate the parameters.  In
summary, the mean values are \new{$\mu=\unit{63.1}{KN/mm^2}$},
\new{$K=\unit{162.2}{KN/mm^2}$}, and \new{$G_c=\unit{0.0031}{KN/mm^2}$}, and
the acceptance rates are \new{28\% (the elasticity parameters) and 21\% (the critical energy rate)}.


\begin{figure}[ht!]
	\centering
	\includegraphics[width=10cm,height=6.25cm]{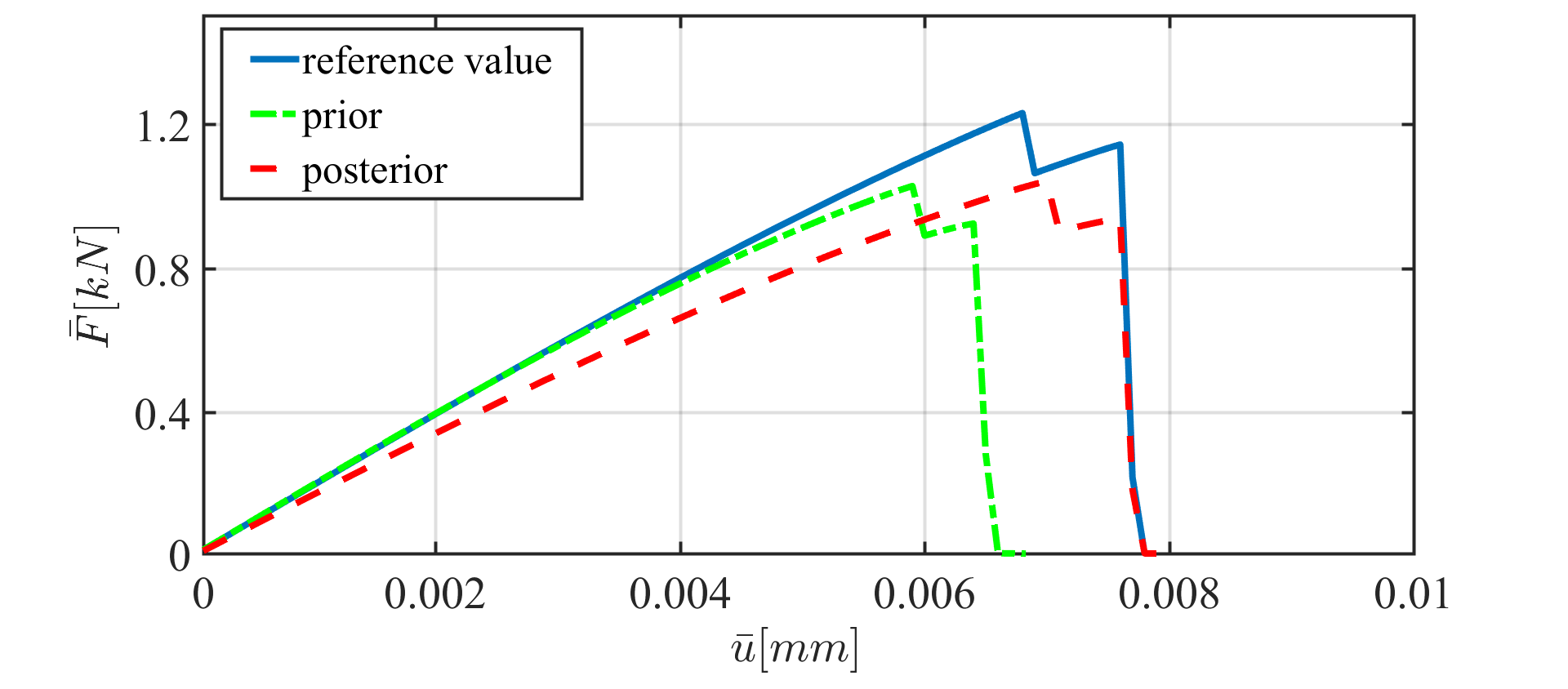}
	\caption{The load-displacement diagram of SENT with voids
		(Example 3). The parameters are the mean values
		($\mu=\unit{80}{kN/mm^2}$,~\new{$K=\unit{170}{kN/mm^2}$},
		and $G_c=\unit{2.7\times10^{-3}}{kN/mm^2}$) obtained by the three-dimensional Bayesian inference.}
	\label{Ex3}
\end{figure}

\begin{figure}[t!]
	\centering
	\subfloat{\includegraphics[width=5cm,height=5cm]{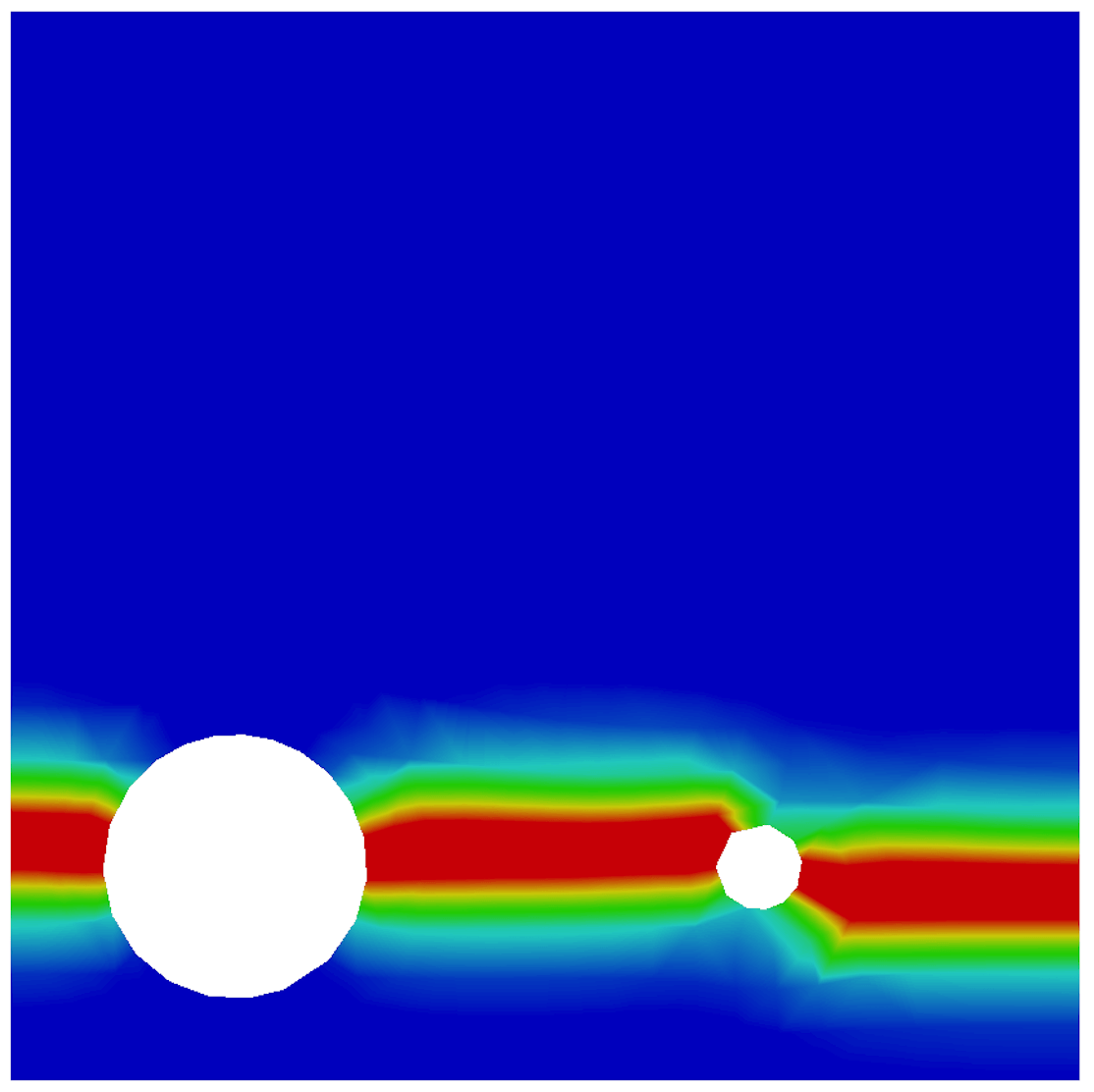}} \hspace{0.5cm}
	\subfloat{\includegraphics[width=5cm,height=5cm]{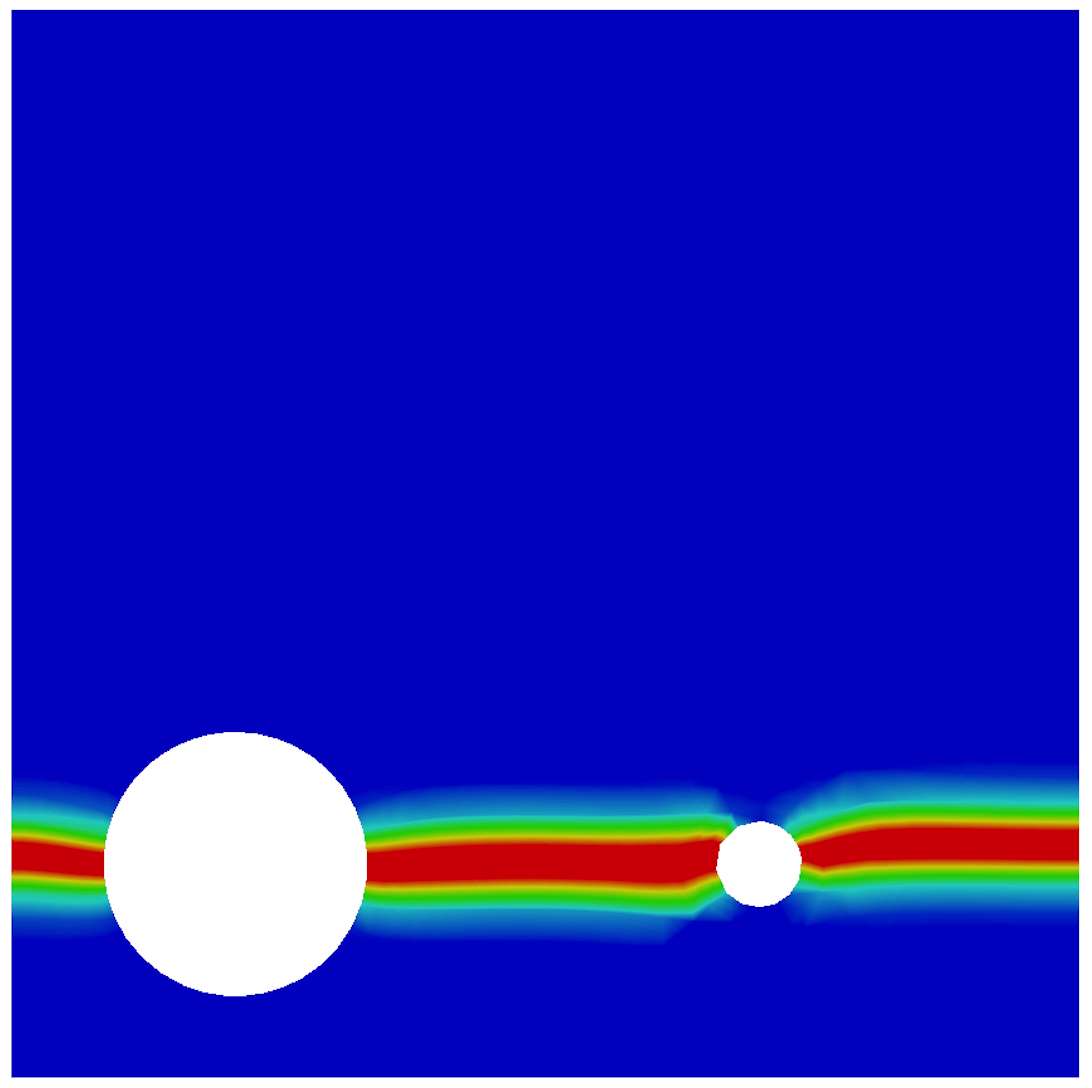}}\hspace{0.5cm}
	\subfloat{\includegraphics[width=5cm,height=5cm]{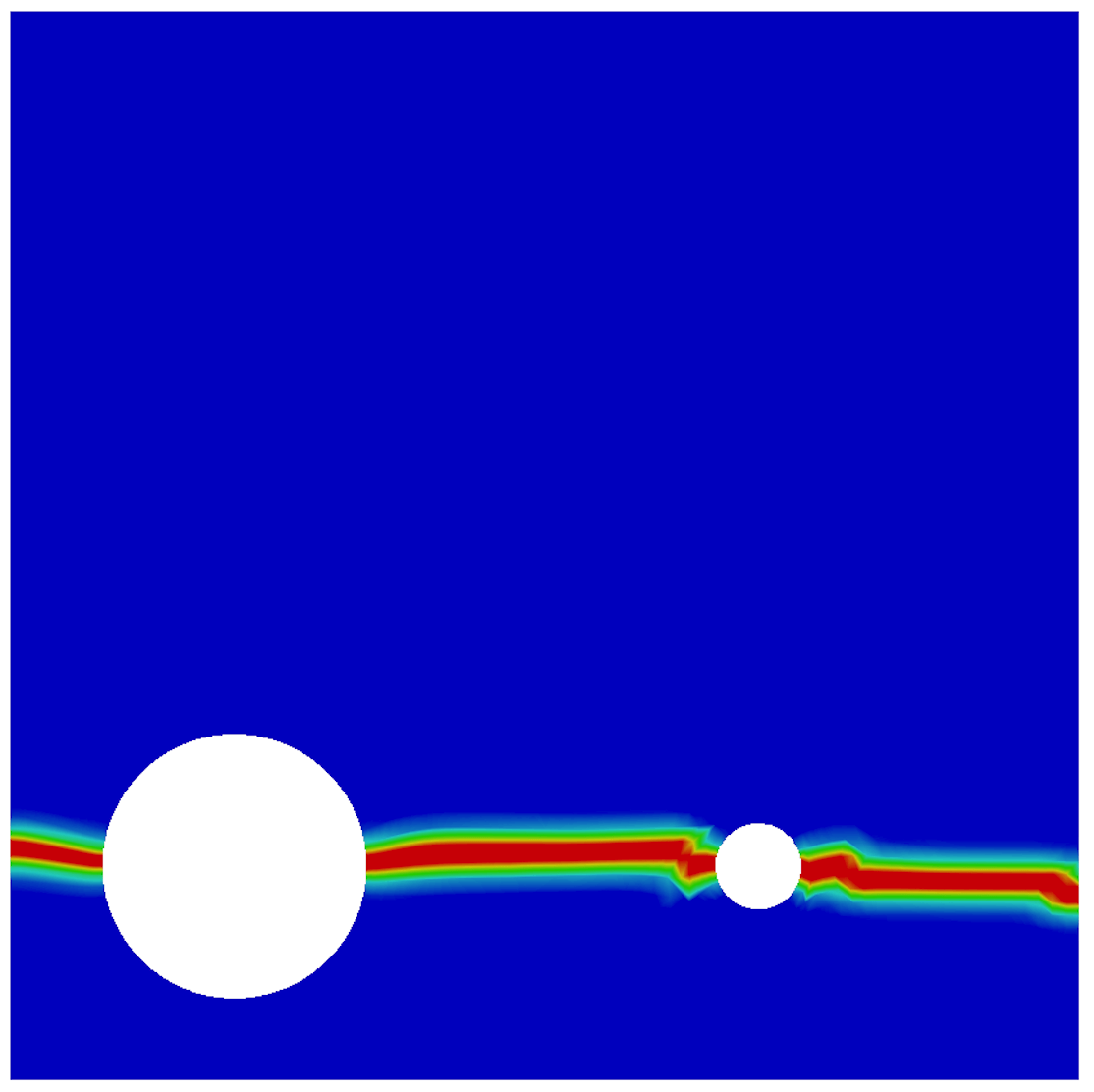}}\newline\\
	\subfloat{\includegraphics[width=5cm,height=5cm]{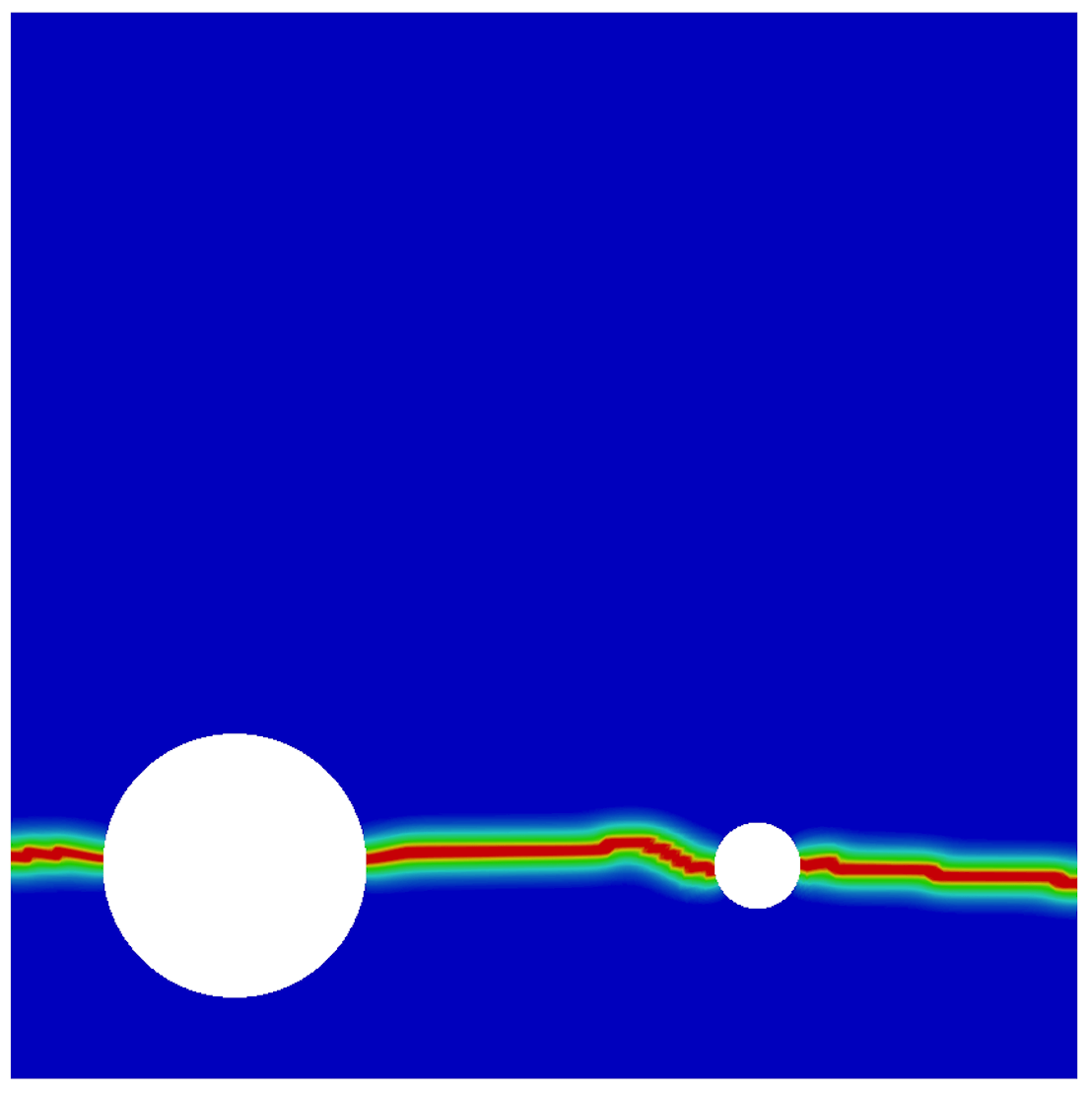}}\hspace{0.5cm}
	\subfloat{\includegraphics[width=5cm,height=5cm]{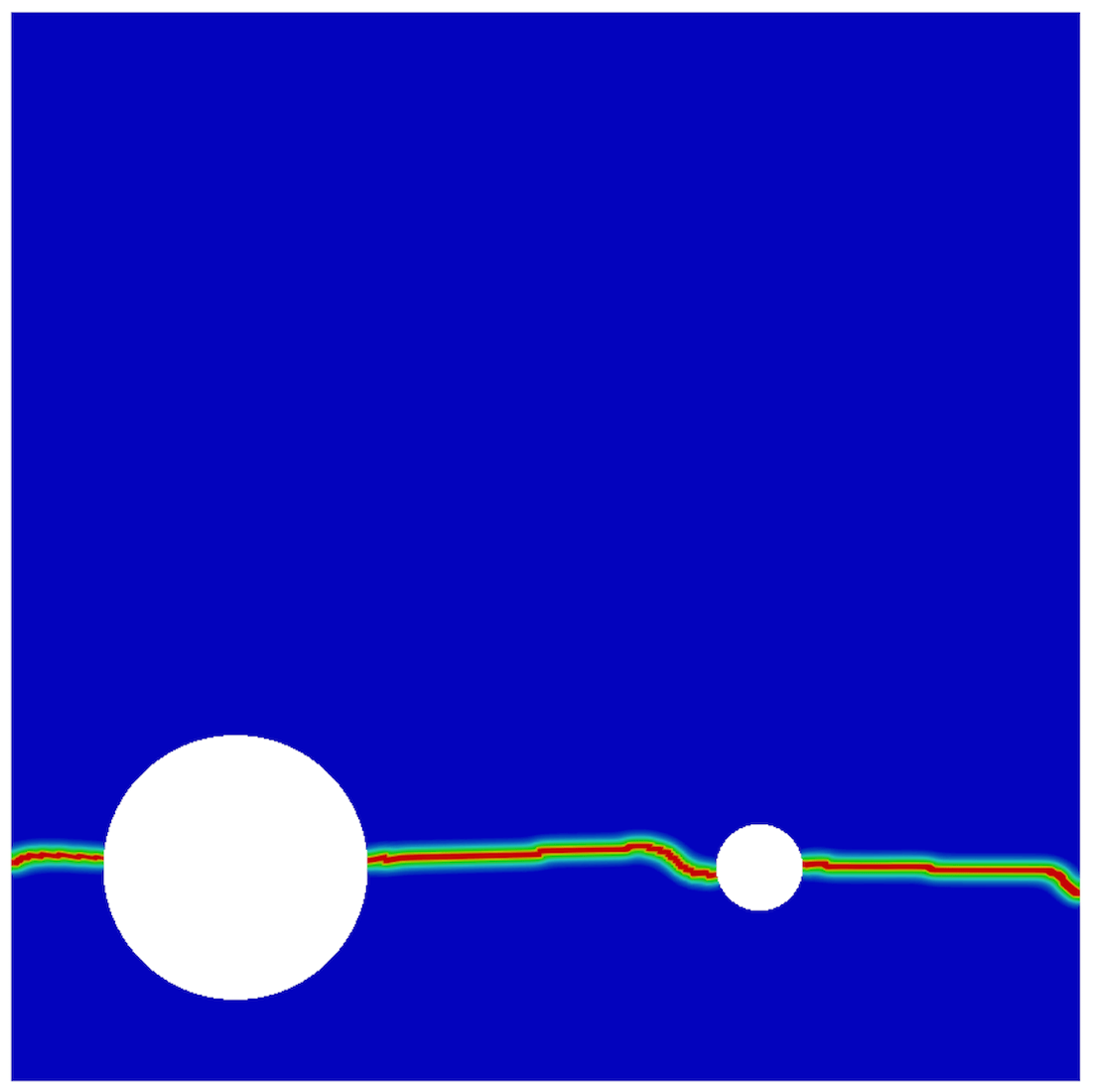}}\hspace{0.5cm}
	\subfloat{\includegraphics[width=2cm,height=4.0cm]{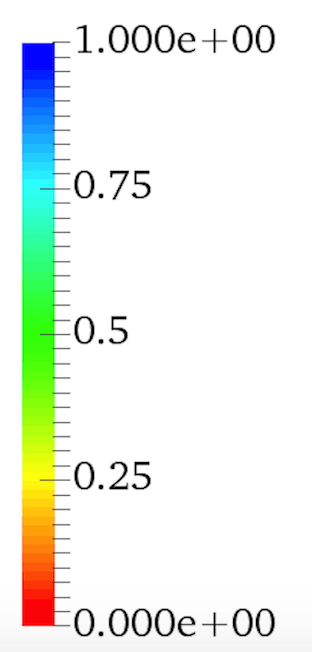}}%
	\caption{The effect of the mesh size on the crack
		propagation. The mesh sizes are (from the left) $h=1/10$, $h=1/20$, $h=1/40$, $h=1/80$, and $h=1/160$ (the reference).}
	\label{fig:exam3_mesh1}
\end{figure}

We solve the forward model with the mean values of the estimated
parameters. As Figure \ref{Ex3} shows, the difference between the
prior distribution and the reference solution is significantly large. By using Bayesian inversion, we could compensate this difference; crack initiation and material failure points are estimated precisely. Although multidimensional Bayesian inversion increases the computational costs (CPU time), the estimated solution is closer to the reference value.

\begin{figure}[t!]
	\centering
	\subfloat{\includegraphics[width=8cm,height=5cm]{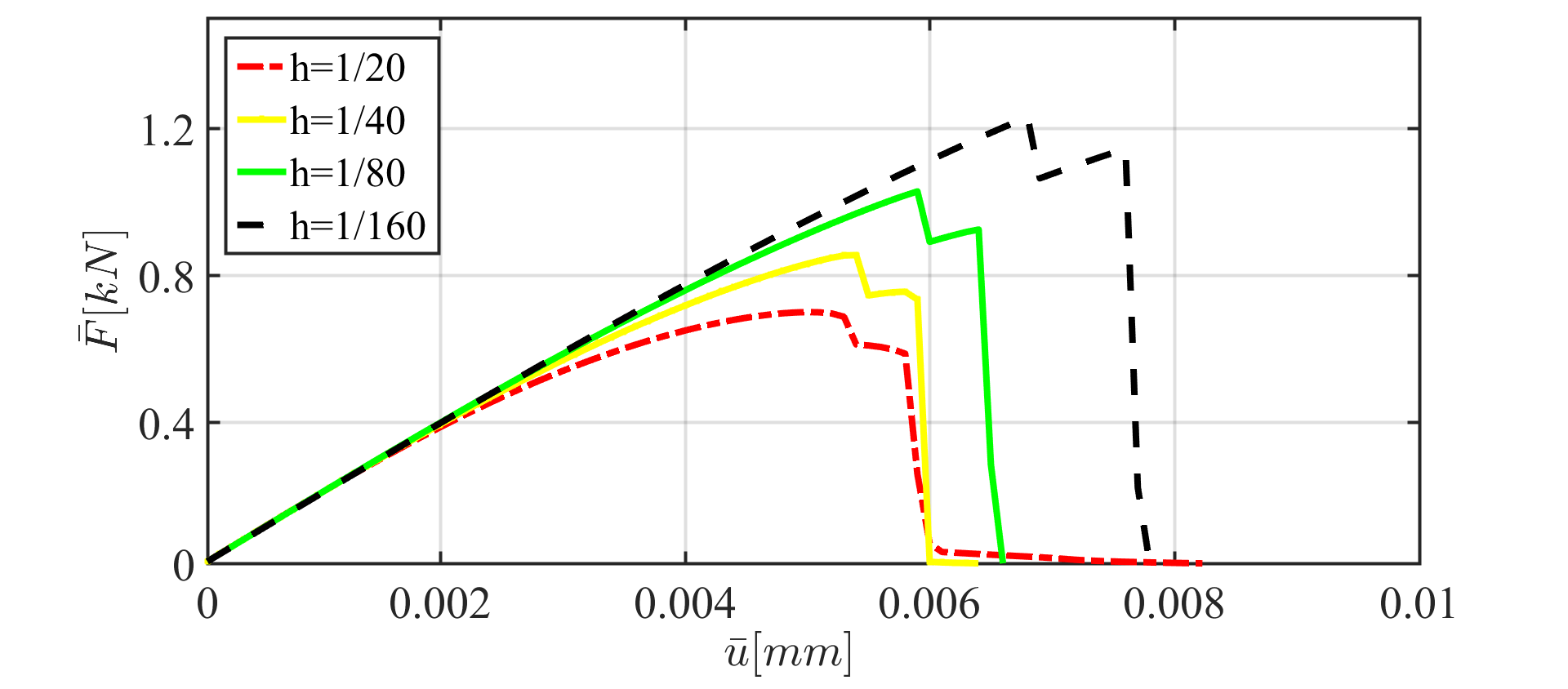}}%
	\hfill 
	\subfloat{\includegraphics[width=8cm,height=5cm]{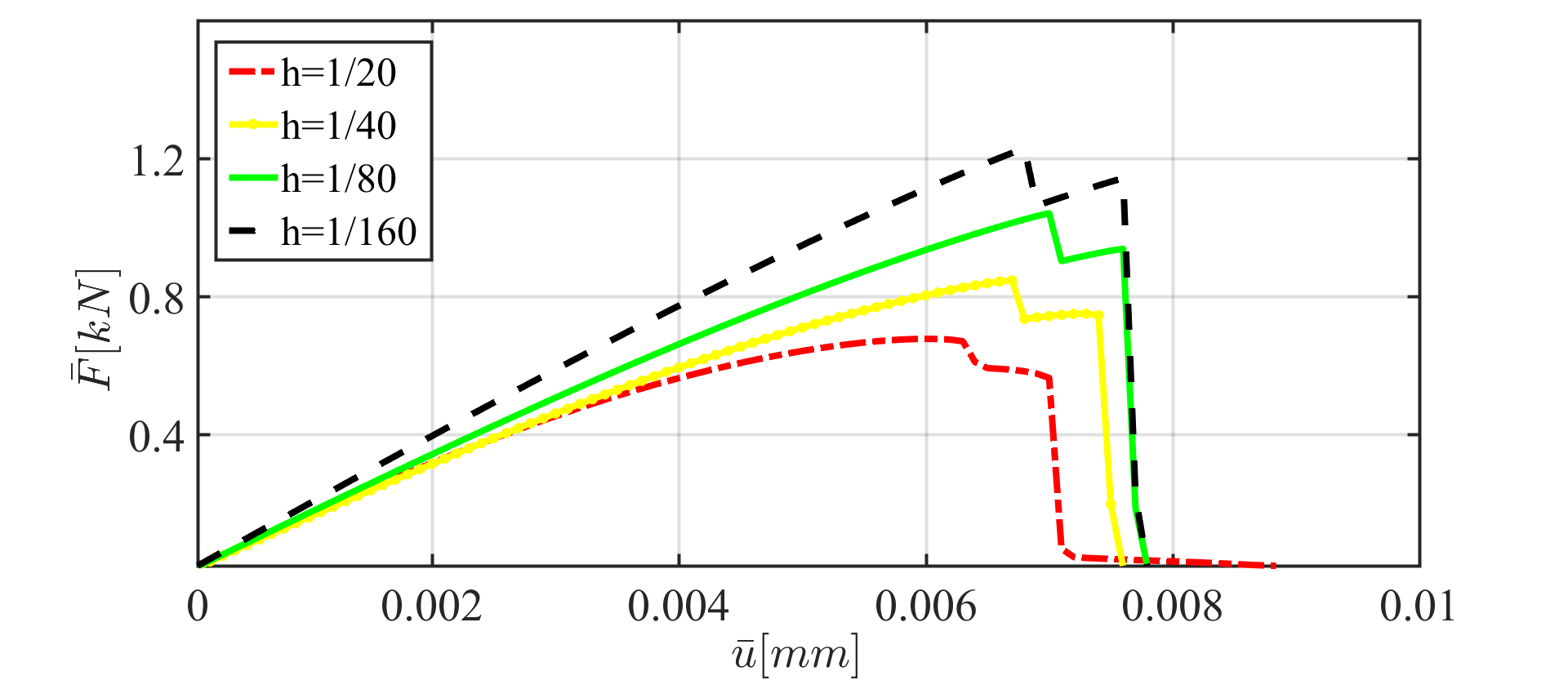}}%
	\caption{The load-displacement curve for different mesh sizes for
		the tension test with voids. The effective parameters are chosen
		according to the prior (left) and posterior (right) distributions.}
	\label{fig:Ex3}
\end{figure} 

Finally, we show the crack patterns obtained by different meshes
varying from $h=1/10$ to $h=1/160$. We use Bayesian inference to
estimate the unknown parameters with the three-dimensional approach
for $h=1/20$ and $h=1/40$. As Figure \ref{fig:Ex3} illustrates,
although the solution based on the posterior distribution is more
precise (i.e., a better estimation of crack initiation and the
fracture point) compared to the one based on the prior distribution,
there is still a difference compared to the reference value. These
results conform to Figure \ref{fig:exam3_mesh1}, since the estimated
crack pattern is considerably larger than the reality.


\section{Conclusions and future works}
\label{conclusions}
In this work, we proposed a Bayesian approach to estimate material parameters 
for propagating fractures in elastic solids. For the fracture model, we adopted a phase-field approach. For the parameter estimation, 
we employed a Bayesian framework. We studied three phase-field
fracture settings, and in each one, 
\new{bulk and shear modulus} as well as
the critical elastic energy release rate were estimated with respect to a reference solution. 


\blue{The developed Bayesian framework enabled us to provide useful knowledge about 
	unknown parameters. }
 By using Bayesian inversion, we could estimate the
load-displacement curve precisely even with coarse meshes. For
instance, in the first example (SENT), the diagram for $h=1/320$ and
$h=1/80$ are essentially same, although a noticeable CPU time
reduction is achieved. Interestingly, using even coarser meshes, the
crack initiations and material fracture times can be estimated very
well in all examples.

\new{As one future application, the Bayesian approach will be used in multiscale
	problems to study crack propagation in heterogeneous materials, e.g.,
	in composites. Due to their complexities, Bayesian inference will be
	employed to estimate material properties when the fiber-reinforced
	structures have a random distribution.
}

\section{Acknowledgments} 

T.~Wick and N.~Noii have been financially supported by the German Research Foundation, 
Priority Program 1748 (DFG SPP 1748) in the 
subproject \textit{Structure Preserving Adaptive Enriched Galerkin Methods for Pressure-Driven 3D Fracture Phase-Field Models} 
with the project No.~392587580. A.~Khodadadian and C.~Heitzinger acknowledge financial support by 
FWF (Austrian Science Fund) START Project no.\ Y660 \textit{PDE Models for Nanotechnology}. M.~Parvizi has been supported by FWF \textit{Project no.\ P28367-N35}. Furthermore, the authors appreciate the useful comments given by the anonymous reviewers.

\end{document}